\documentclass[11pt,a4paper]{article}
\usepackage[top=3cm, bottom=3cm, left=3cm, right=3cm]{geometry}

\usepackage{booktabs}
\usepackage{color}
\usepackage{latexsym} 
\usepackage{amsmath}               
\usepackage{amssymb}               
\usepackage{amsfonts}              
\usepackage{amsthm}                
\usepackage{tcolorbox}
\usepackage{times}
\usepackage[utf8]{inputenc}
\usepackage[T1]{fontenc}
\usepackage{enumitem}
\usepackage{xcolor}
\usepackage{thmtools}
\usepackage{thm-restate}
\usepackage[square,numbers]{natbib}
\usepackage[colorlinks=true,allcolors=blue]{hyperref}
\usepackage{cleveref}
\usepackage{url}            
\usepackage{booktabs} 
\usepackage{nicefrac}       
\usepackage{microtype}      
\usepackage[american]{babel}
\usepackage{thmtools}
\usepackage{booktabs}
\usepackage{algorithm}
\usepackage{algorithmic}
\usepackage[tikz]{bclogo}
\usepackage{multicol,tabularx,capt-of}
\usepackage{hhline}
\usepackage{multirow}

\renewenvironment{proof}[1][\proofname]{{\bfseries #1.}}{}

\makeatother

\newcommand{\eqals}[1]{\begin{align*}#1\end{align*}}
\newcommand{\eqal}[1]{\begin{align}#1\end{align}}
\newcommand{\bpr}{\begin{proof}}
\newcommand{\epr}{\end{proof}}
\newcommand{\be}{\begin{equation}}
\newcommand{\ee}{\end{equation}}

 \newtheorem{example}{Example}
  \newtheorem{theorem}{Theorem}
  \newtheorem{lemma}{Lemma}
  \newtheorem{proposition}{Proposition}
  \newtheorem{remark}{Remark}
  \newtheorem{corollary}{Corollary}
  \newtheorem{definition}{Definition}

\newcommand{\bd}{\begin{definition}}
\newcommand{\ed}{\end{definition}}

\newcommand{\bi}{\begin{itemize}}
\newcommand{\ei}{\end{itemize}}

\newtheorem{ass}{Assumption}
\newcommand{\ba}{\begin{ass}}
\newcommand{\ea}{\end{ass}}

\newtheorem{mtd}{Method}
\newcommand{\bmtd}{\begin{mtd}}
\newcommand{\emtd}{\end{mtd}}

\newcommand{\bre}{\begin{restatable}}
\newcommand{\ere}{\end{restatable}}

\newcommand{\br}{\begin{remark}}
\newcommand{\er}{\end{remark}}

\newcommand{\bp}{\begin{proposition}}
\newcommand{\ep}{\end{proposition}}

\newcommand{\blm}{\begin{lemma}}
\newcommand{\elm}{\end{lemma}}

\newcommand{\bt}{\begin{theorem}}
\newcommand{\et}{\end{theorem}}

\newcommand{\bcor}{\begin{corollary}}
\newcommand{\ecor}{\end{corollary}}

\newcommand{\bex}{\begin{example}}
\newcommand{\eex}{\end{example}}

\crefname{proposition}{Proposition}{Propositions}

\crefname{definition}{Definition}{Definitions}

\crefname{ass}{Assumption}{Assumptions}
\crefname{equation}{Eq.}{Eqs.}
\crefname{figure}{Fig.}{Figs.}
\crefname{table}{Table}{Tables}
\crefname{section}{Sec.}{Secs.}

\crefname{theorem}{Thm.}{Thms.}

\crefname{lemma}{Lemma}{Lemmas}

\crefname{corollary}{Cor.}{Cors.}

\crefname{example}{Example}{Examples}

\crefname{appendix}{Appendix}{Appendixes}

\crefname{remark}{Remark}{Remark}

\newcommand{\ka}{\kappa}

\newcommand{\Kfinal}{\mathsf{K}}
\newcommand{\tHess}{\widetilde{\mathbf{H}}}

\newcommand{\Ab}{\mathbf{A}}
\newcommand{\tAb}{\widetilde{\mathbf{A}}}

\newcommand{\Pj}{\mathbf{P}}
\newcommand{\Id}{\mathbf{I}}

\newcommand{\xla}{x_\la^{\star}}

\newcommand{\Deff}{\mathsf{df}}

\newcommand{\Exp}[1]{\mathbb{E}\left[#1 \right]}

\newcommand{\Expb}[2]{\mathbb{E}_{#1}\left[#2\right]}

\newcommand{\clas}{\mathcal{P}}

\newcommand{\Bos}{\mathsf{B}_1^\star}
\newcommand{\Bts}{\mathsf{B}_2^\star}

\newcommand{\Qs}{\mathsf{Q}^\star}
\newcommand{\Lc}{\mathsf{L}}

\newcommand{\Qc}{\mathsf{Q}}

\newcommand{\rla}{\mathsf{r}_\la}

\newcommand{\nm}[1]{\mathsf{t}(#1)}
\newcommand{\tn}{\mathsf{t}}

\newcommand{\Rad}{R}

\newcommand{\xlap}{x_{\Pj,\la}^*}
\newcommand{\cp}{\mathsf{C}_{\Pj}}
\newcommand{\srce}{\mathsf{s}}
\newcommand{\rpla}{\mathsf{r}_{\Pj,\la}}

\newcommand{\X}{\mathcal{X}}
\newcommand{\Y}{\mathcal{Y}}
\newcommand{\Z}{\mathcal{Z}}
\newcommand{\R}{\mathbb{R}}

\newcommand{\N}{\mathbb{N}}

\renewcommand{\L}{L}

\newcommand{\Hess}{\mathbf{H}}
\newcommand{\nd}{\nu}
\newcommand{\ndt}{\widetilde{\nu}}
\newcommand{\ns}{\Delta}
\newcommand{\nsa}{\widetilde{\Delta}}
\newcommand{\cc}{\mathsf{c}}

\newcommand{\radd}{\mathsf{r}}
\newcommand{\dik}{\mathsf{D}}
\newcommand{\hf}{\widehat{f}}
\newcommand{\hxla}{\widehat{x}_\la^\star}
\newcommand{\hG}{\widehat{\mathcal{G}}}

\newcommand{\hL}{\widehat{L}}
\newcommand{\hHess}{\widehat{\mathbf{H}}}

\newcommand{\hnm}[1]{\widehat{\mathsf{t}}(#1)}
\newcommand{\hrla}{\widehat{\mathsf{r}}_{\lambda}}
\newcommand{\hnd}{\widehat{\nu}}

\newcommand{\bias}{\mathsf{b}}

\newcommand{\la}{\lambda}
\newcommand{\hh}{{\mathcal{H}}}

\DeclareMathOperator{\supp}{supp}
\DeclareMathOperator{\lspan}{span}

\DeclareMathOperator{\Tr}{Tr}

\DeclareMathOperator{\poly}{poly}
\DeclareMathOperator{\diag}{diag}

\newcommand{\dikin}{\underline{\phi}}

\newcommand{\dikins}{\overline{\phi}}

\newcommand{\bvar}{\widehat{\mathsf{v}}}

\DeclareMathOperator*{\argmin}{arg\,min}

\newcommand{\G}{\mathcal{G}}
\newcommand{\xo}{x^\star}
\newcommand{\hsrce}{\widehat{s}}
\newcommand{\htn}{\widehat{\mathsf{t}}}
\newcommand{\hcp}{\widehat{\mathsf{C}}_{\Pj}}
\newcommand{\hradd}{\widehat{\mathsf{r}}}
\newcommand{\Cone}{\mathsf{C}_1}
\newcommand{\Ctwo}{\mathsf{C}_2}
\newcommand{\Cfone}{\mathsf{K}_1}
\newcommand{\Cftwo}{\mathsf{K}_2}
\newcommand{\Cfthree}{\mathsf{K}_3}
\newcommand{\xmu}{x_\mu^\star}
\newcommand{\tmu}{\widetilde{\mu}}
\newcommand{\err}{\rho}
\DeclareMathOperator{\lso}{LinApprox}

\newcommand{\era}{\epsilon_0}
\newcommand{\erat}{\widetilde{\epsilon}_0}
\newcommand{\terr}{\tilde{\rho}}
\newcommand{\Bb}{\mathbf{B}}
\newcommand{\fstar}{f^\star}
\newcommand{\hAb}{\widehat{\Ab}}
\newcommand{\Ny}{\mathcal{N}}
\newcommand{\Nyi}{\mathcal{N}_{\infty}}

\newcommand{\Gb}{\mathbf{G}}
\newcommand{\bts}{\mathsf{b}_2^\star}
\newcommand{\btns}{\mathsf{b}_2}
\newcommand{\Cov}{\mathbf{\Sigma}}
\DeclareMathOperator{\cond}{Cond}
\newcommand{\btb}{\overline{\mathsf{b}}_2}
\newcommand{\Nyb}{\overline{\mathcal{N}}}
\newcommand{\Nyib}{\overline{\mathcal{N}}_{\infty}}
\newcommand{\Rell}{\mathsf{R}_{\ell}}

\newcommand{\pthm}{ \widehat{\Gamma}^M}
\newcommand{\csla}{\mathsf{c}_{\texttt{samp}}}

\renewcommand{\leq}{\leqslant}
\renewcommand{\geq}{\geqslant}

\usepackage{authblk}

\author{\bfseries
Ulysse Marteau-Ferey}
\author{\bfseries
Francis Bach}
\author{\bfseries
Alessandro Rudi
}
\affil{
INRIA - D{\'e}partement d'Informatique de l'{\'E}cole Normale Sup{\'e}rieure \\
PSL Research University\\
Paris, France
}

\begin{document}

\title{Globally Convergent Newton Methods for Ill-conditioned Generalized Self-concordant Losses}

\date{}


\maketitle

\begin{abstract}
In this paper, we study large-scale convex optimization algorithms based on the Newton method applied to regularized generalized self-concordant losses, which include logistic regression and softmax regression.  We first prove that our new simple scheme based on a sequence of problems with decreasing regularization parameters is provably globally convergent, that this convergence is linear with a constant factor which scales only logarithmically with the condition number. In the parametric setting, we obtain an algorithm with the same scaling than regular first-order methods but with an improved behavior, in particular in ill-conditioned problems. Second, in the non-parametric machine learning setting, we provide an explicit algorithm combining  the previous scheme with Nystr\"om projection techniques, and prove that it achieves optimal generalization bounds with a time complexity of order $O(n\Deff_\la)$, a memory complexity of order $O(\Deff_\la^2)$ and {\em no dependence on the condition number}, generalizing the results known for least-squares regression. Here $n$ is the number of observations and $\Deff_\la$ is the associated degrees of freedom. In particular, this is the first large-scale algorithm to solve logistic and softmax regressions in the non-parametric setting with large condition numbers and theoretical guarantees.
\end{abstract}

\section{Introduction}
Minimization algorithms constitute a crucial algorithmic part of many machine learning methods, with   algorithms available for a variety of situations~\cite{bottou2018optimization}. In this paper, we focus on \emph{finite sum} problems of the form
\[\min_{x \in \hh}{f_\la(x) = f(x) + \frac{\la}{2}\|x\|^2}, \mbox{ with }  f(x) = \frac{1}{n}\sum_{i=1}^n{f_i(x)},\] 
where $\hh$ is a Euclidean or a Hilbert space, and each function is convex and smooth. The running-time of minimization algorithms classically depends on the number of functions $n$, the explicit (for Euclidean spaces) or implicit (for Hilbert spaces) dimension $d$ of the search space, and the condition number of the problem, which is upper bounded by   $\ka = L/\lambda$, where $L$ characterizes the smoothness of the functions~$f_i$, and $\lambda$ the regularization parameter.
    
In the last few years, there has been a strong focus on problems with large $n$ and $d$, leading to \emph{first-order} (i.e., gradient-based) stochastic algorithms, culminating in a sequence of linearly convergent algorithms whose running time is favorable in $n$ and $d$, but scale at best in $\sqrt{\ka}$~\cite{defazio2014saga,lin2015universal,defazio2016simple,Allen2017}.  
However, modern problems lead to objective functions with very large condition numbers, i.e., in many learning problems, the regularization parameter that is optimal for test predictive performance may be so small that the scaling above in $\sqrt{\ka}$ is not practical anymore (see examples in Sect.~\ref{sec:exp}).
    
These ill-conditioned problems are good candidates for \emph{second-order methods}  (i.e., that use the Hessians of the objective functions) such as Newton method. These methods are traditionally discarded within machine learning for several reasons: (1) they are usually adapted to high precision results which are not necessary for generalization to unseen data for machine learning problems~\cite{bottou2008tradeoffs}, (2) computing the Newton step $\ns_\la(x) = \nabla^2 f_\la(x)^{-1}\nabla f_\la(x)$ requires to form the Hessian and solve the associated linear system, leading to complexity which is at least quadratic in $d$, and thus prohibitive for large $d$, and (3) the global convergence properties are not applicable, unless the function is very special, i.e., self-concordant~\cite{nemirovskii1994interior} (which includes only few classical learning problems),  so they often are only shown to converge in a  small area around the optimal $x$.

In this paper, we argue that the three reasons above for not using Newton method can be circumvented to obtain competitive algorithms: (1) high {absolute} precisions are indeed not needed for machine learning, but faced with strongly ill-conditioned problems, even a low-precision solution requires second-order schemes; (2) many approximate Newton steps have been designed for approximating the solution of the associated large linear system \cite{Monta15,RoostaKhorasaniMahoney19,pilanci2017newton,bollapragada2018exact}; (3) we propose a novel second-order method which is globally convergent and which is based on performing approximate Newton methods for a certain class of so-called \emph{generalized self-concordant functions} which includes logistic regression~\cite{bach2010self}. For these functions, the conditioning of the problem is also characterized by a more \emph{local} quantity: $\ka_\ell = \Rad^2/\la$, where $\Rad$ characterizes the local evolution of Hessians. This leads to second-order algorithms which are competitive with first-order algorithms for well-conditioned problems, while being superior for ill-conditioned problems which are common in practice.

\paragraph{Contributions.} We make the following contributions:
\begin{itemize}

\item[$(a)$] We build a global second-order method for the minimization of $f_\la$, which relies only on computing approximate Newton steps of the functions $f_\mu,\mu \geq \la$. The number of such steps will be of order $O(c \log \ka_\ell + \log \frac{1}{\epsilon})$ where $\epsilon$ is the desired precision, and $c$ is an explicit constant. In the parametric setting ($\hh = \R^d$), $c$ can be as bad as $\sqrt{\ka_\ell}$ in the worst-case but much smaller in theory and practice. Moreover in the non-parametric/kernel machine learning setting ($\hh$ infinite dimensional), $c$ does not depend on the local condition number~$\ka_\ell$.

\item[$(b)$] Together with the appropriate quadratic solver to compute approximate Newton steps, we obtain an algorithm with the same scaling as regular first-order methods but with an improved behavior, in particular in ill-conditioned problems. Indeed, this algorithm matches the performance of the best quadratic solvers but covers any generalized self-concordant function, up to logarithmic terms.

\item[$(c)$] In the non-parametric/kernel machine learning setting we provide an explicit algorithm combining the previous scheme with Nystr\"om projections techniques. We prove that it achieves optimal generalization bounds with $O(n\Deff_\la)$ in time and $O(\Deff_\la^2)$ in memory, where~$n$ is the number of observations and $\Deff_\la$ is the associated degrees of freedom. In particular, this is the first large-scale algorithm to solve logistic and softmax regression in the non-parametric setting with large condition numbers and theoretical guarantees. 
\end{itemize}

\subsection{Comparison to related work}\label{sec:comparison}
We consider two cases for $\hh$ and the functions $f_i$ that are common  in machine learning:  $\hh = \R^d$ with linear (in the parameter) models with explicit feature maps, and  $\hh$ infinite-dimensional, corresponding in machine learning to learning with kernels~\cite{shawe2004kernel}. Moreover in this section we first consider the quadratic case, for example the squared loss in machine learning (i.e., $f_i(x) = \frac{1}{2}(x^\top z_i - y_i)^2$ for some $z_i \in \hh, y_i \in \R$). We first need to introduce the Hessian of the problem, for any $\la >0$, define 
\[ \Hess(x):= \nabla^2 f(x) , \qquad \Hess_\la(x):= \nabla^2 f_\la(x) = \Hess(x) + \la \Id,\]
in particular we denote by $\Hess$ (and analogously $\Hess_\la$) the Hessian at optimum (which in case of squared loss corresponds to the covariance matrix of the inputs). 

\paragraph{Quadratic problems and $\hh = \R^d$ (ridge regression).}
The problem then consists in solving a (ill-conditioned) positive semi-definite symmetric linear system of dimension $d \times d$. Methods based on {\em randomized linear algebra}, {\em sketching} and suitable {\em subsampling} \cite{drineas2011faster,drineas2012fast,boutsidis2013improved} are able to find the solution with precision $\epsilon$ in time that is $O((n d + \min(n,d)^3)\log(L/\la \epsilon))$, so  essentially independently of the condition number, because of the logarithmic complexity in $\lambda$.

\paragraph{Quadratic problems and $\hh$ infinite-dimensional (kernel ridge regression).} Here the problem corresponds to solving a (ill-conditioned) infinite-dimensional linear system in a reproducing kernel Hilbert space \cite{shawe2004kernel}. Since however the sum defining $f$ is finite, the problem can be projected on a subspace of dimension at most $n$~\cite{aronszajn1950theory}, leading to a linear system of dimension $n \times n$. Solving it with the techniques above would lead to a complexity of the order $O(n^2)$, which is not feasible on massive learning problems (e.g., $n \approx 10^7$). Interestingly these problems are usually approximately low-rank, with the rank represented by the so called {\em effective-dimension} $\Deff_\la$ \cite{devito}, counting essentially the eigenvalues of the problem larger than $\la$, 
\eqal{\label{eq:eff-dim-text}
\Deff_\la = \Tr(\Hess \Hess_\la^{-1}).
}
Note that $\Deff_\la$ is bounded by $\min\{ n, L/\la\}$ and in many cases $\Deff_\la \ll \min(n, L/\la)$.  Using suitable projection techniques, like {\em Nystr\"om}~\cite{williams2001using} or {\em random features}~\cite{rahimi2008random} it is possible to further reduce the problem to dimension $\Deff_\la$, for a total cost to find the solution of $O(n \Deff_\la^2)$. Finally recent methods~\cite{Rudi17}, combining suitable projection methods with refined preconditioning techniques, are able to find the solution with precision compatible with the optimal statistical learning error \cite{devito} in time that is $O(n \Deff_\la \log(L/\lambda))$, so being essentially independent of the condition number of the problem.

\paragraph{Convex problems and explicit features (logistic regression).} When the loss function is {\em self-concordant} it is possible to leverage the fast techniques for linear systems in approximate Newton algorithms \cite{pilanci2017newton} (see more in \cref{sec:background}), to achieve the solution in essentially $O(n d + \min(n,d)^3)$ time, modulo logarithmic terms. However only few loss functions of interest are self-concordant, in particular the widely used logistic and soft-max losses are not self-concordant, but {\em generalized-self-concordant}~\cite{bach2010self}. In such cases we need to use (accelerated/stochastic) first order optimization methods to enter in the quadratic convergence region of Newton methods \cite{Agarwal2017}, which leads to a solution in $O(d n + d \sqrt{ n L/\la} + \min(n,d)^3)$ time, which does not present any improvement on a simple accelerated first-order method.  Globally convergent second-order methods have also been proposed to solve such problems \cite{jaggi18}, but the number of Newton steps needed being bounded only by $L/\la$, they lead to a solution in $O(L/\la~(nd +\min(n,d)^3))$. With $\la$ that could be as small as $10^{-12}$ in modern machine learning problems, this makes both these kind of approaches expensive from a computational viewpoint for ill-conditioned problems. For such problems, with our new global second-order scheme, the algorithm we propose achieves instead a complexity of essentially $O((n d + \min(n,d)^3)\log(\Rad^2/\la\epsilon))$ (see \cref{thm:main-thm}).

\paragraph{Convex problems and $\hh$ infinite-dimensional (kernel logistic regression).}
Analogously to the case above, it is not possible to use Newton methods profitably as global optimizers on losses that are not self-concordant as we see in \cref{sec:global_cv_scheme}. In such cases by combining projecting techniques developped in \cref{sec:kernels} and accelerated first-order optimization methods, it is possible to find a solution in $O(n\Deff_\la + \Deff_\la \sqrt{ n  L/\la})$ time. This can still be prohibitive in the very small regularization scenario, since it strongly depends on the condition number $L/\la$. In \cref{sec:kernels} we suitably combine our optimization algorithm with projection techniques achieving optimal statistical learning error~\cite{marteau2019} in essentially $O(n \Deff_\la \log(\Rad^2/\la))$.

\paragraph{First-order algorithms for finite sums.}
In dimension $d$, accelerated algorithms for strongly-convex smooth (not necessarily self-concordant) finite sums, such as K-SVRG~\cite{Allen2017}, have a running time proportional $O( (n+ \sqrt{n L / \lambda})d)$. This can be improved with preconditioning to $O( (n+ \sqrt{d L / \lambda})d)$ for large $n$~\cite{Agarwal2017}. Quasi-Newton methods can also be used~\cite{gower2018accelerated},  but typically without the guarantees that we provide in this paper (which are logarithmic in the condition number in natural scenarios).

%
%

\section{Background: Newton methods and generalized self concordance}\label{sec:background}
In this section we start by recalling the definition of generalized self concordant functions and motivate it with examples. We then recall basic facts about Newton and approximate Newton methods, and present existing techniques to efficiently compute approximate Newton steps. We start by introducing the definition of generalized self-concordance, that here is an extension of the one in~\cite{bach2010self}.
\bd[generalized self-concordant (GSC) function]\label{df:gensc}
Let $\hh$ be a Hilbert space. We say that $f$ is a generalized self-concordant function on $\G \subset \hh$, when $\G$ is a bounded subset of $\hh$ and $f$ is a convex and three times differentiable mapping on $\hh$ such that
\[ \textstyle \forall x \in \hh,~ \forall h,k \in \hh,~ \nabla^{(3)}f(x)[h,k,k] \leq \sup_{g \in \G}|g \cdot h| ~ \nabla^2 f(x)[k,k].\]
\ed
We will usually denote by $\Rad$ the quantity $\sup_{g \in \G}\|g\| < \infty$ and often omit $\G$ when it is clear from the context (for simplicity think of $\G$ as the ball in $\hh$ centered in zero and with radius $R > 0$, then $\sup_{g \in \G}|g \cdot h| = \Rad\|h\|$).
The globally convergent second-order scheme we present in \cref{sec:global_cv_scheme} is specific to losses which satisfy this generalized self-concordance property. The following loss functions, which are widely used in machine learning, are generalized-self-concordant, and motivate this work.
\begin{example}[Application to finite-sum minimization]\label{ex:loss-text}
The following loss functions are generalized self-concordant functions, but not self-concordant: 
\\[.05cm]
(a) Logistic  regression: $f_i(x)= \log(1 + \exp(  - y_i w_i^\top x))$, where $x, w_i \in \R^d$ and $y_i \in \{-1,1\}$.
\\[.05cm]
(b) Softmax regression: $f_i(x) = \log \big( \sum_{j=1}^k \exp(x_j^\top w_i)\big) - x_{y_i}^\top w_i$, where now $x \in \R^{d \times k}$ and $y_i \in \{1,\dots,k\}$ and $x_j$ denotes the $j$-th column of $x$.
 \\[.05cm]
    (c) Generalized linear models with bounded features (see details in \cite[Sec.~2.1]{bach2014adaptivity}), which include conditional random fields~\cite{sutton2012introduction}.
\\[.05cm]
    (d) Robust regression: $f_i(x) = \varphi(y_i -w_i^\top x)$ with $\varphi(u) = \log(e^u+e^{-u})$.
\end{example}

Note that these losses are not \textit{self-concordant} in the sense of \cite{pilanci2017newton}. Moreover, even if the losses $f_i$ are self-concordant, the objective function $f$ is not necessarily self-concordant, making any attempt to prove the self-concordance of the objective function $f$ almost impossible.  

\paragraph{Newton method (NM).}
Given $x_0 \in \hh$, the Newton method consists in doing the following update:
\eqal{\label{eq:newton-intro}
x_{t+1} = x_t - \ns_\la(x_t),\qquad \ns_\la(x_t):= \Hess^{-1}_\la(x_t)\nabla f_\la(x_t).
}
The quantity $\ns_\la(x):= \Hess^{-1}_\la(x)\nabla f_\la(x)$ is called the Newton step at point $x$, and $x - \ns_\la(x)$ is the minimizer of the second order approximation of $f_\la$ around $x$. Newton methods enjoy the following key property: if $x_0$ is close enough to the optimum, the convergence to the optimum is quadratic and the number of iterations required to a given precision is independent of the condition number of the problem \cite{boyd2004convex}.  

However Newton methods have two main limitations:  (a) the region of quadratic convergence can be quite small and reaching the region can be computationally expensive, since it is usually done via first order methods \cite{Agarwal2017} that converge linearly depending on the condition number of the problem, (b) the cost of computing the Hessian can be really expensive when $n, d$ are large, and also (c) the cost of computing $\ns_\la(x_t)$ can be really prohibitive. In the rest of the section we recall some ways to deal with (b) and (c). Our main result of \cref{sec:global_cv_scheme} is to provide globalization scheme for the Newton method to tackle  problem (a), which is easily integrable with approximate techniques to deal with (b) ans (c), to make second-order technique competitive. 

\paragraph{Approximate Newton methods (ANM) and approximate solutions to linear systems.}
Computing exactly the Newton increment $\ns_\la(x_t)$, which corresponds essentially to the solution of a linear system, can be too expensive when $n, d$ are large. A natural idea is to approximate the Newton iteration, leading to {\em approximate Newton methods},
\eqal{\label{eq:approx-newton-intro}
x_{t+1} = x_t - \widetilde{\ns}_\la(x_t), \qquad \widetilde{\ns}_\la \approx \ns_\la(x_t).
}
In this paper, more generally we consider any technique to compute $\widetilde{\ns}_\la(x_t)$ that provides a {\em relative approximation} \cite{Deuflhard2011} of $\ns_\la(x_t)$ defined as follows.
\begin{definition}[relative approximation]
Let $\err < 1$, let $\Ab$ be an invertible positive definite Hermitian operator on $\hh$ and $b$ in $\hh$.
We denote by $\lso(\Ab,b,\err)$ the set of all  $\rho$-relative approximations of $z^* = \Ab^{-1}b$, i.e., $\lso(\Ab,b,\err) = \{z \in \hh ~|~ \|z - z^*\|_{\Ab} \leq \err \|z^*\|_{\Ab}\}$.
\end{definition}

\paragraph{Sketching and subsampling for approximate Newton methods.}
Many techniques for approximating linear systems have been used to compute $\widetilde{\ns}_\la$, in particular {\em sketching} of the Hessian matrix via fast transforms and {\em subsampling} (see \cite{pilanci2017newton,bollapragada2018exact,Agarwal2017} and references therein). Assuming for simplicity that $f_i = \ell_i(w_i^\top x)$, with $\ell_i: \R \to \R$ and $w_i \in \hh$, it holds:
\eqal{\label{eq:hessian-dec}
 \Hess(x) = \frac{1}{n} \sum_{i=1}^n \ell_i^{(2)}(w_i^\top x) w_i w_i^\top = V_x^\top V_x,
}
    with $V_x \in \R^{n\times d} = D_x W$, where $D_x \in \R^{n\times n}$ is a diagonal matrix defined as $(D_x)_{ii} = (\ell_i^{(2)}(w_i^\top x))^{1/2}$ and $W \in \R^{n \times d}$ defined as $W = (w_1, \dots, w_n)^\top$. 
    
    Both sketching and subsampling methods approximate $z^* = \Hess_\la(x)^{-1} \nabla f_\la(x)$ with $\tilde{z} = \widetilde{\Hess}_\la(x)^{-1} \nabla f_\la(x)$, in particular, in the case of subsampling $\widetilde{\Hess}(x) = \sum_{j=1}^Q p_{j} w_{i_j} w_{i_j}^\top$ where $Q \ll \min(n,d)$, $(p_j)_{j=1}^n$ are suitable weights and $(i_j)_{j=1}^Q$ are indices selected at random from $\{1,\dots,n\}$ with suitable probabilities. Sketching methods instead use $\widetilde{\Hess}(x) = \widetilde{V}_x^\top \widetilde{V}_x$, with $\widetilde{V}_x =  \Omega V_x$ with $\Omega \in \R^{Q \times n}$ a structured matrix such that computing $\widetilde{V}_x$ has a cost in the order of $O(n d \log n)$; to this end usually $\Omega$ is based on fast Fourier or Hadamard transforms~\cite{pilanci2017newton}.
Note that essentially all the techniques used in approximate Newton methods guarantee relative approximation. In particular the following results can be found in the literature (see \cref{lm:sub_unif,lm:sub_ny} in \cref{sec:lemmas-operators} and \cite{pilanci2017newton}, Lemma~2 for more details).

\blm\label{lm:appr-linear-system}
Let $x, b \in \hh$ and assume that $\ell_i^{(2)} \leq a$ for $a > 0$. With probability $1-\delta$ the following methods output an element in $\lso(\Hess_\la(x), b,\rho)$, in $O(Q^2 d + Q^3 + c)$ time, $O(Q^2+d)$ space:
\\[.05cm]
(a) Subsampling with uniform sampling (see \cite{RoostaKhorasaniMahoney19,Rudi15}), where $Q = O(\rho^{-2}a/\la\log \frac{1}{\la\delta})$ and $c = O(1)$. 
\\[.05cm]
(b) Subsampling with approximate leverage scores \cite{RoostaKhorasaniMahoney19,alaoui2015fast,Rudi15}), where $Q = O(\rho^{-2} \bar{\Deff_\la}\log 1/\la\delta), c = O(\min(n,a/\la)\bar{\Deff_\la}^2)$ and $\bar{\Deff_\la} = \Tr(W^\top W(W^\top W + \la/a I)^{-1})$ \cite{rudi2018fast}. Note that $\bar{\Deff_\la} \leq \min(n, d)$.\\
(c) Sketching with {\em fast Hadamard transform} \cite{pilanci2017newton}, where $Q = O(\rho^{-2} \bar{\Deff_\la}\log a/\la\delta), c = O(nd \log n)$.
\elm

\section{Globally convergent scheme for ANM algorithms on GSC functions\label{sec:global_cv_scheme}}
The algorithm is based on the observation that when $f_\la$ is generalized self concordant, there exists a region where $t$ steps of ANM converge as fast as $2^{-t}$. Our idea is to start from a very large regularization parameter $\la_0$, such that we are sure that $x_0$ is in the convergence region and perform some steps of ANM such that the solution enters in the convergence region of $f_{\la_1}$, with $\la_1 = q \la_0$ with $q < 1$, and to iterate this procedure until we enter the convergence region of $f_\la$. First we define the region of interest and characterize the behavior of NM and ANM in the region, then we analyze the globalization scheme.  

\paragraph{Preliminary results: the Dikin ellipsoid.}
We consider the following region that we prove to be contained in the region of quadratic convergence for the Newton method and that will be useful to build the globalization scheme. Let $c,R > 0$ and $f_\la$ be generalized self-concordant with coefficient~$R$, we call {\em Dikin ellipsoid} and denote by $\dik_\la(\cc)$ the region
\[\dik_\la(\cc):= \big\{x ~|~ \nd_\la(x) \leq \cc \sqrt{\la}/R  \big\}, \quad \textrm{with} \quad \nd_\la(x):= \|\nabla f_\la(x)\|_{\Hess^{-1}_\la(x)},\]
where $\nd_\la(x)$ is usually called the {\em Newton decrement} and $\|x\|_{\Ab}$ stands for $\|\Ab^{1/2} x\|$.
\blm\label{lm:dikin-is-cool}
Let $\la > 0, \cc \leq 1/7$, let $f_\la$ be generalized self-concordant and $x \in \dik_\la(\cc)$. Then it holds:
   $
    \frac{1}{4}\nd_\la(x)^2 \leq f_\la(x) - f_\la(\xla) \leq \nd_\la(x)^2$.
Moreover Newton method starting from $x_0$ has quadratic convergence, i.e., let $x_t$ be obtained via $t \in \N$ steps of Newton method in \cref{eq:newton-intro}, then
$ \nd_\la(x_t) \leq 2^{-(2^{t} - 1)} \nd_\la(x_0).$
Finally, approximate Newton methods starting from $x_0$ have a linear convergence rate, i.e., let $x_t$ given by \cref{eq:approx-newton-intro}, with $\nsa_t \in \lso(\Hess_\la(x_t),\nabla f_\la(x_t),\err)$ and $\err \leq 1/7$, then
$\nd_\la(x_t) \leq 2^{-t} \nd_\la(x_0).$
\elm

This result is proved in \cref{lm:dikin-is-cool_2} in \cref{app:newton_particular}. 
The crucial aspect of the result above is that when $x_0 \in \dik_\la(\cc)$, the convergence of the approximate Newton method is linear and does not depend on the condition number of the problem. However $\dik_\la(\cc)$ itself can be very small depending on $\sqrt{\la}/R$. In the next subsection we see how to enter in $\dik_\la(\cc)$ in an efficient way.

\paragraph{Entering the Dikin ellipsoid using a second-order scheme.}
The lemma above shows that $\dik_\la(\cc)$ is a good region where to use the approximate Newton algorithm on GSC functions. However the region itself is quite small, since it depends on $\sqrt{\la}/R$. Some other globalization schemes arrive to regions of interest by first-order methods or back-tracking schemes \cite{Agarwal2017,Monta15}. However such approaches require a number of steps that is  usually proportional to $\sqrt{L/\la}$ making them non-beneficial in machine learning contexts. Here instead we consider the following simple scheme where $\texttt{ANM}_\rho(f_\la,x,t)$ is the result of a $\rho$-relative approximate Newton method performing $t$ steps of optimization starting from $x$.

The main ingredient to guarantee the scheme to work is the following lemma (see \cref{lm:next_la} in \cref{app:dec_la_tech} for a proof).

\blm\label{lm:dec_la}
Let $\mu > 0$, $\cc < 1$ and $x \in \hh$.  Let $s = 1 + \Rad\|x\|/\cc$, then for $q \in [1 - 2/(3s), 1)$
\[\dik_\mu(\cc/3) \subseteq \dik_{q\mu}(\cc).\]
\elm 
Now we are ready to show that we can guarantee the loop invariant $x_k \in \dik_{\mu_k}(\cc)$.
Indeed assume that $x_{k-1} \in \dik_{\mu_{k-1}}(\cc)$. Then $\nd_{\mu_{k-1}}(x_{k-1}) \leq \cc\sqrt{\mu_{k-1}}/R$. By taking $t = 2, \rho = 1/7$, and performing $x_{k} = \texttt{ANM}_\rho(f_{\mu_{k-1}},x_{k-1},t)$, by \cref{lm:dikin-is-cool}, 
$\nd_{\mu_{k-1}}(x_{k}) \leq 1/4 \nd_{\mu_{k-1}}(x_{k-1}) \leq \cc/4~\sqrt{\mu_{k-1}}/R$,
i.e., $x_{k} \in \dik_{\mu_{k-1}}(\cc/4)$. If $q_{k}$ is large enough, this implies that $x_{k} \in \dik_{q_{k}\mu_{k-1}}(\cc) = \dik_{\mu_{k}}(\cc)$, by \cref{lm:dec_la}. Now we are ready to state our main theorem of this section.

\vspace*{.3cm}

\fbox{
\begin{minipage}[t]{0.9\textwidth}
\textbf{Proposed Globalization Scheme}
\begin{center}
\textit{Phase I: Getting in the Dikin ellispoid of $f_\la$}
\end{center}

\vspace*{-.1cm}

Start with $x_0 \in \hh, \mu_0 > 0$, $t, T \in \N$ and $(q_k)_{k \in \N} \in (0,1]$.\\
For $k \in \N$\\
${}\qquad x_{k+1} \leftarrow \texttt{ANM}_\rho(f_{\mu_k},x_{k},t)$\\
${}\qquad \mu_{k+1} \leftarrow q_{k+1}\mu_{k}$\\
Stop when $\mu_{k+1} < \la$ and set $x_{last}\leftarrow x_k$.
\begin{center}
\textit{Phase II: reach a certain precision starting from  inside the Dikin ellipsoid}
\end{center}

\vspace*{-.1cm}

Return $\widehat{x} \leftarrow \texttt{ANM}_\rho(f_{\la},x_{last},T)$
\end{minipage}
}

\paragraph{Fully adaptive method.} The scheme presented above converges with the following parameters. 

\bt \label{thm:easy_first_phase}\label{thm:main-thm}
Let $\epsilon >0$. Set  $\mu_0 = 7 \Rad \|\nabla f (0)\|$, $x_0 = 0$, and perform the globalization scheme above for $\err \leq 1/7, t = 2$, and
$ \textstyle q_k = \frac{1/3 + 7 \Rad \|x_{k}\|}{1 + 7 \Rad \|x_{k}\|}$, $ T = \lceil \log_2 \sqrt{1 ~\vee ~ (\la \epsilon^{-1} / R^2)}\rceil.$
Then denoting by $K$ the number of steps performed in the Phase I, it holds:
\[f_\la(\widehat{x}) - f_\la(x_\la^\star) \leq \epsilon, \qquad K \leq \left\lfloor \left(3 + 11 \Rad \|\xla\|\right)\log (7 \Rad \|\nabla f(0)\|/\la)\right\rfloor.\]
\et 
 Note that the theorem above (proven in \cref{app:proof_main_thm}) guarantees a solution with error $\epsilon$ with $K$ steps of ANM each performing 2 iterations of approximate linear system solving, plus a final step of ANM which performs $T$ iterations of approximate linear system solving. In case of $f_i(x) = \ell_i(w_i^\top x)$, with $\ell_i:\R\to\R$, $w_i \in \hh$ with $\ell_i^{(2)} \leq a$, for $a > 0$, the final runtime cost of the proposed scheme to achieve precision $\epsilon$, when combined with of the methods for approximate linear system solving from \cref{lm:appr-linear-system} (i.e. sketching), is $O(Q^2 + d)$ in memory and
\[O\Big((nd\log n + dQ^2 + Q^3)\Big(\Rad \|\xla\|\log\frac{R}{\la} + \log\frac{\la}{R\epsilon}\Big)\Big)~~\textrm{in time},\quad Q = O\Big(\bar{\Deff_\la}\log\frac{1}{\la\delta}\Big),\]
where $\bar{\Deff_\la}$, defined in \cref{lm:appr-linear-system}, measures the {\em effective dimension} of the correlation matrix $W^\top W$ with $W = (w_1,\dots, w_n)^\top \in \R^{n \times d}$, corresponding essentially to the number of eigenvalues of $W^\top W$ larger than $\la/a$. In particular note that $\bar{\Deff_\la} \leq \min(n, d, \text{rank}(W), a b^2/\la)$, with  $b:= \max_i\|w_i\| $, and usually way smaller than such quantities.
\br \label{rm:cond-dependence}
The proposed method does not depend on the condition number of the problem $L/\la$, but on the term $\Rad \|\xla\|$ which can be in the order of $\Rad/\sqrt{\la}$ in the worst case, but usually way smaller. For example, it is possible to prove that this term is bounded by an absolute constant not depending on $\la$, if at least one minimum for $f$ exists. In the appendix (see \cref{prp:adap_var}), we show a variant of this adaptive method which can leverage the regularity of the solution with respect to the Hessian, i.e., depending on the smaller quantity $\Rad \sqrt{\la}\|x^\star_{\la}\|_{\Hess^{-1}_{\la}(x^\star_{\la})}$ instead of $R\|\xla\|$.
\er
Finally note that it is possible to use $q_k = q$ fixed for all the iterations and way smaller than the one in \cref{thm:main-thm}, depending on some regularity properties of $\Hess$ (see \cref{prp:fixed_q} in \cref{app:maintheorems}).

\section{Application to the non-parametric setting: Kernel methods}\label{sec:kernels}
In supervised learning the goal is to predict well on future data, given the observed training dataset. Let $\X$ be the input space and $\Y \subseteq \R^p$ be the output space. We consider a probability distribution $P$ over $\X \times \Y$ generating the data and the goal is to estimate $g^*:\X \to \Y$ solving the problem 
\eqal{\label{eq:ideal-problem}
g^* = \argmin_{g:\X \to \Y} \mathcal{L}(g), \quad \mathcal{L}(g) = \mathbb{E}[\ell(g(x),y)],
}
for a given loss function $\ell:\Y\times \Y \to \R$. Note that $P$ is not known, and accessible only via the dataset $(x_i,y_i)_{i=1}^n$, with $n \in \N$, independently sampled from $P$. A prototypical estimator for $g^*$ is the regularized minimizer of the empirical risk $\widehat{\mathcal{L}}(g) = \frac{1}{n}\sum_{i=1}^n \ell(g(x_i), y_i)$ over a suitable space of functions $\cal G$. Given $\phi: \X \to \hh$ a common choice is to select ${\cal G}$ as the set of linear functions of $\phi(x)$, that is, ${\cal G} = \{ w^\top \phi(\cdot) ~|~ w \in \hh\}$. Then the regularized minimizer of $\widehat{\mathcal{L}}$, denoted by $\widehat{g}_\la$, corresponds to  
\eqal{\label{eq:ERM-text}
\widehat{g}_\la(x) = \widehat{w}_\la^\top \phi(x), \quad   \widehat{w}_\la = \argmin_{w \in \hh} \textstyle \frac{1}{n} \sum_{i=1}^n f_i(w) + \la\|w\|^2, \quad f_i(w) = \ell(w^\top \phi(x_i), y_i).
}
Learning theory guarantees how fast $\widehat{g}_\la$ converges to the best possible estimator $g^*$ with respect to    the number of observed examples, in terms of the so called {\em excess risk} $\mathcal{L}(\widehat{g}_\la) - \mathcal{L}(g^*)$. The following theorem recovers the minimax optimal learning rates for squared loss and extend them to any generalized self-concordant loss function.

{\em Note on $\Deff_\la$.} In this section, we always denote with $\Deff_\la$ the effective dimension of the problem in \cref{eq:ideal-problem}. 
When the loss belongs to the family of generalized linear models (see \cref{ex:loss-text}) and if the model is well-specified, then $\Deff_\la$ is defined exactly as in \cref{eq:eff-dim-text} otherwise we need a more refined definition (see \cite{marteau2019} or  \cref{eq:deff-statistics} in \cref{app:kernels}).

\bt[from \cite{marteau2019}, Thm. 4]\label{thm:opt-rates-ERM-text}
Let $\la > 0, \delta \in (0,1]$. Let $\ell$ be generalized self-concordant with parameter $R > 0$ and $\sup_{x \in X} \|\phi(x)\| \leq C < \infty$. Assume that there exists $g^*$ minimizing $\mathcal{L}$. Then there exists $c_0$ not depending on $n, \lambda, \delta, \Deff_\la, C, g^*$, such that if $\sqrt{\Deff_\la/n},\bias_\la \leq \la^{1/2}/\Rad$, and $n \geq C/\la \log(\delta^{-1}C/\la)$ the following holds with probability $1-\delta$:
\begin{equation}
\label{eq:opt-rates-ERM-text}
\mathcal{L}(\widehat{g}_\la) - \mathcal{L}(g^*) \leq c_0\Big(\frac{\Deff_\la}{n} + \bias^2_\la\Big)\log (1/\delta), \qquad \bias_\la:=  \la \|g^*\|_{\Hess_\la^{-1}(g^*)}.
\end{equation}
\et
Under standard regularity assumptions of the learning problems \cite{marteau2019}, i.e., (a) the {\em capacity condition} $\sigma_j(\Hess(g^*)) \leq C j^{-\alpha}$, for $\alpha \geq 1,  C > 0$ (i.e., a decay of eigenvalues $\sigma_j(\Hess(g^*))$ of the Hessian at the optimum), and (b) the {\em source condition}  $g^* = \Hess(g^*)^r v$, with $v \in \hh$ and $r > 0$ (i.e., the control of the optimal $g^\ast$ for a specific Hessian-dependent norm),   $\Deff_\la \leq C' \la^{-1/\alpha}$ and $\bias^2_\la \leq C^{\prime \prime} \la^{1+2r}$, leading to the following optimal learning rate,
\eqal{\label{eq:rates-text}
\mathcal{L}(\widehat{g}_\la) - \mathcal{L}(g^*)  \leq c_1 n^{-\frac{1+2r \alpha}{1 + \alpha + 2r\alpha}} \log (1/\delta), \quad \textrm{when} \quad \la = n^{-\frac{\alpha}{1+\alpha+2r\alpha}}.
}
Now we propose an algorithmic scheme to compute efficiently an approximation of $\widehat{g}_\la$   that achieves the same optimal learning rates. First we need to introduce the technique we are going to use. 

\paragraph{Nystr\"om projection.} It consists in suitably selecting $\{\bar{x}_1, \dots, \bar{x}_M\} \subset \{x_1,\dots, x_n\}$, with $M \ll n$ and computing $\bar{g}_{M,\la},$ i.e., the solution of \cref{eq:ERM-text} over $\hh_M = \text{span}\{\phi(\bar{x}_1),\dots, \phi(\bar{x}_M)\}$ instead of $\hh$. In this case the problem can be reformulated as a problem in $\R^M$ as
\eqal{\label{eq:prob-ny-char-text}
\bar{g}_{M,\la} = \bar{\alpha}_{M,\la}^\top {\bf T}^{-1} v(x), \qquad \bar{\alpha}_{M,\la} = \argmin_{\alpha \in \R^M} \bar{f}_\la(\alpha), \qquad \bar{f}(\alpha) = \frac{1}{n} \sum_{i=1}^n \bar{f}_i(\alpha) + \la\|\alpha\|^2,
}
where $\bar{f}_i(\alpha) = \ell( v(x_i)^\top {\bf T}^{-1}\alpha,~y_i)$ and $v(x) \in \R^M$, $v(x) = (k(x,\bar{x}_1), \dots, k(x,\bar{x}_M))$ with $k(x,x') = \phi(x)^\top \phi(x')$ the associated positive-definite kernel~\cite{shawe2004kernel}, while ${\bf T}$ is the upper triangular matrix such that ${\bf K} = {\bf T}^\top {\bf T} $, with ${\bf K} \in \R^{M \times M}$ with ${\bf K}_{ij} = k(\bar{x}_i,\bar{x}_j)$. In the next theorem we characterize the sufficient $M$ to achieve minimax optimal rates, for two standard techniques of choosing the Nystr\"om points $\{\bar{x}_1, \dots, \bar{x}_M\}$.

\bt[Optimal rates for learning with Nystr\"om]\label{thm:opt-rates-ny-text}
Let $\la > 0, \delta \in (0,1]$. Assume the conditions of \cref{thm:opt-rates-ERM-text}. Then the excess risk of $\bar{g}_{M,\la}$ is bounded with prob. $1-2\delta$ as in \cref{eq:opt-rates-ERM-text} (with $c_1' \propto c_1$), when\\
${}\qquad(1)~$ Uniform Nystr\"om method \cite{Rudi15,Rudi17} is used and $~M ~\geq~ C_1/\la~\log(C_2/\la\delta)$.\\
${}\qquad(2)~$ Approximate leverage score method \cite{alaoui2015fast,Rudi15,Rudi17} is used and $~M ~\geq~ C_3~\Deff_\la~\log(C_4/\la\delta)$.\\
Here $C, C_1, C_2, C_4$ do not depend on $\la, n, M, \Deff_\la, \delta$.
\et
\cref{thm:opt-rates-ny-text} generalizes results for learning with Nystr\"om and squared loss \cite{Rudi15}, to GSC losses. It is proved in \cref{thm:nystrom_estimator_bounds}, in \cref{app:selecting_nys}. As in \cite{Rudi15}, \cref{thm:opt-rates-ny-text} shows that Nystr\"om is a valid technique for dimensionality reduction. Indeed it is essentially possible to project the learning problem on a subspace $\hh_M$ of dimension $M = O(c/\la)$ or even as small as $M = O(\Deff_\la)$ and still achieve the optimal rates of \cref{thm:opt-rates-ERM-text}. Now we are ready to introduce our algorithm.

\paragraph{Proposed algorithm.}
The algorithm conceptually consists in (a) performing a projection step with Nystr\"om, and  (b) solving the resulting optimization problem with the globalization scheme proposed in \cref{sec:global_cv_scheme} based on ANM in \cref{eq:approx-newton-intro}. 
In particular, we want to avoid to apply explicitly ${\bf T}^{-1}$ to each $v(x_i)$ in \cref{eq:prob-ny-char-text}, which would require $O(n M^2)$ time. Then we will use the following approximation technique based only on matrix vector products, so we can just apply ${\bf T}^{-1}$ to $\alpha$ at each iteration, with a total cost proportional only to $O(n M  + M^2)$ per iteration.
Given $\alpha, \nabla \bar{f}_\la(\alpha)$, we approximate $z^* =  \bar{\Hess}_\la(\alpha)^{-1} \nabla \bar{f}_\la(\alpha)$, where $\bar{\Hess}_\la$ is the Hessian of $\bar{f}_\la(\alpha)$, with $\tilde{z}$ defined as
\[ \tilde{z} = \texttt{prec-conj-grad}_t(\bar{\Hess}_\la(\alpha), \nabla \bar{f}_\la(\alpha)),\]
where $\texttt{prec-conj-grad}_t$ corresponds to performing $t$ steps of preconditioned conjugate gradient~\cite{golub2012matrix} with preconditioner computed using a subsampling approach for the Hessian among the ones presented in \cref{sec:background}, in the paragraph starting with \cref{eq:hessian-dec}. The pseudocode for the whole procedure is presented in Alg.~\ref{alg:algo}, \cref{sec:algorithm}.
This technique of approximate linear system solving has been studied in \cite{Rudi17} in the context of empirical risk minimization for squared loss.

\blm[\cite{Rudi17}]\label{lm:pcg-appr-tech}
Let $\la > 0, \alpha, b \in \R^M$. The previous method, applied with $t = O(\log 1/\rho)$, outputs an element of $\lso(\bar{\Hess}_{\la}(\alpha), b,\rho)$, with probability $1-\delta$ with complexity $O((nM + M^2 Q + M^3 + c)t)$ in time and $O(M^2 + n)$ in space, with $Q = O(C_1/\la \log (C_1/\la\delta)), c = O(1)$ if uniform sub-sampling is used or $Q = O(C_2\Deff_\la \log (C_1/\la\delta)), c = O(\Deff_\la^2\min(n,\frac{1}{\la}))$ if sub-sampling with leverage scores is used \cite{rudi2018fast}.
\elm

A more complete version of this lemma is shown in \cref{prp:approx_ns_kernels} in \cref{app:ans_kernels}. 
We conclude this section with a result proving the learning properties of the proposed algorithm.

\bt[Optimal rates for the proposed algorithms]\label{thm:final-alg-kernel-text}
Let $\la > 0$ and $\epsilon < \la/\Rad^2$. Under the hypotheses of \cref{thm:opt-rates-ny-text}, if we set $M$ as in \cref{thm:opt-rates-ny-text}, $Q$ as in \cref{lm:pcg-appr-tech} and setting the globalization scheme as in \cref{thm:main-thm}, then the proposed algorithm (Alg.~\ref{alg:algo}, \cref{sec:algorithm}) finishes in a finite number of newton steps $N_{ns} = O(\Rad \|g^{*}\|\log(C/\la) +\log(C/\epsilon))$ and returns a predictor $g_{Q,M,\la}$ of the form $g_{Q,M,\la} = \alpha^\top {\bf T}^{-1}v(x)$. With probability at least $1-\delta$, this predictor satisfies: 
\begin{equation}
\label{eq:opt-alg}
\mathcal{L}(g_{Q,M,\la}) - \mathcal{L}(g^*) \leq c_0\Big(\frac{\Deff_\la}{n} + \bias^2_\la + \epsilon\Big)\log (1/\delta), \qquad \bias_\la:=  \la \|g^*\|_{\Hess_\la^{-1}(g^*)}.
\end{equation}
\et
The theorem above (see \cref{prp:algo_stat_gen}, \cref{app:main_kernels} for exacts quantifications) shows that the proposed algorithm is able to achieve the same learning rates of plain empirical risk minimization as in \cref{thm:opt-rates-ERM-text}. The total complexity of the procedure, including the cost of computing the preconditioner,  the selection of the Nystr\"om points via approximate leverage scores and also the computation of the leverage scores \cite{rudi2018fast} is then
$$O\left(\Rad \|g^*\| \log (\Rad^2/\la)\left(n~\Deff_\la \log(C\la^{-1}\delta^{-1})~c_X +  ~+~ \Deff_\la^3 \log^3(C\la^{-1} \delta^{-1}) + \min(n,C/\la)~\Deff_\la^2 \right)\right)$$ in time and $O(\Deff_\la^2\log^2(C\la^{-1}\delta^{-1}))$ in space, where $c_X$ is the cost of computing the inner product $k(x,x')$ (in the kernel setting assumed when the input space $X$ is $X=\R^p$ it is $c=O(p)$). As noted in \cite{rudi2018fast}, under the standard regularity assumptions on the learning problem seen above, $\Deff_\la^2 \leq \Deff_\la/\la \leq n$ when the optimal $\la$ is chosen. So the total computational complexity  is
\[ O\left(\Rad  \log (\Rad^2/\la)~ \log^3(C\la^{-1} \delta^{-1})~\|g^*\|\cdot  n \cdot\Deff_\la \cdot c_X \right) ~~\textrm{in time}, \quad O(\Deff_\la^2\cdot \log^2(C \la^{-1}\delta^{-1})) ~~\textrm{in space}.\]
First note, the fact that due to the statistical properties of the problem the complexity does not depend even implicitly on $\sqrt{C/\la}$, but only on $\log(C/\la)$, so the algorithm runs in essentially $O(n \Deff_\la)$, compared to $O(\Deff_\la \sqrt{n C/\la})$ of the accelerated first-order methods we develop in \cref{sec:exp-appendix} and the $O(n \Deff_\la \sqrt{ C/\la})$ of other Newton schemes (see \cref{sec:comparison}). To our knowledge,
this is the first algorithm to achieve optimal statistical learning rates for generalized self-concordant losses and with complexity only $\widetilde{O}(n \Deff_\la)$. This generalizes similar results for squared loss \cite{Rudi17,rudi2018fast}.

\section{Experiments}
\label{sec:exp}

The code necessary to reproduce the following experiments is available on GitHub at \url{https://github.com/umarteau/Newton-Method-for-GSC-losses-}.

We compared the performances of our algorithm for kernel logistic regression on two large scale classification datasets ($n\approx 10^7$), Higgs and Susy, pre-processed as in \cite{Rudi17}. We implemented the algorithm in pytorch and performed the computations on $1$ Tesla P100-PCIE-16GB GPU. For Susy ($n=5 \times 10^6, p=18$): we used Gaussian kernel with $k(x,x') = e^{-\|x-x'\|^2/(2\sigma^2)}$, with $\sigma=5$, which we obtained through a grid search (in \cite{Rudi17}, $\sigma = 4$ is taken for the ridge regression); $M=10^4$ Nystr\"om centers and a subsampling $Q = M$ for the preconditioner, both obtained with uniform sampling. Analogously for Higgs ($n=1.1 \times 10^7, p=28$): , we used a Gaussian kernel with $\sigma=5$ and $M=2.5\times10^4$ and $Q=M$, using again uniform sampling.
\begin{figure}
    \centering
    \includegraphics[width=0.44\textwidth]{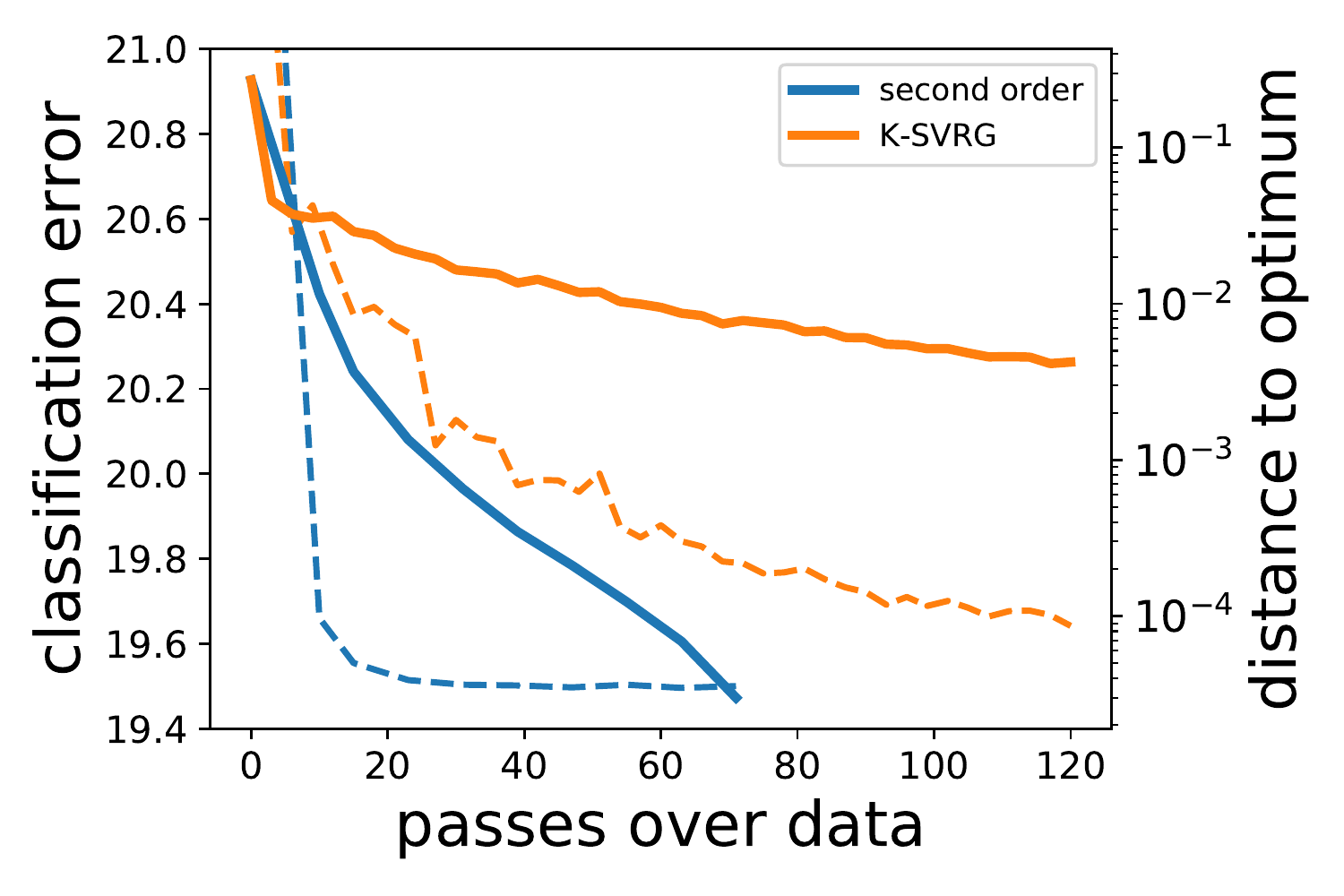}~~
    \includegraphics[width=0.44\textwidth]{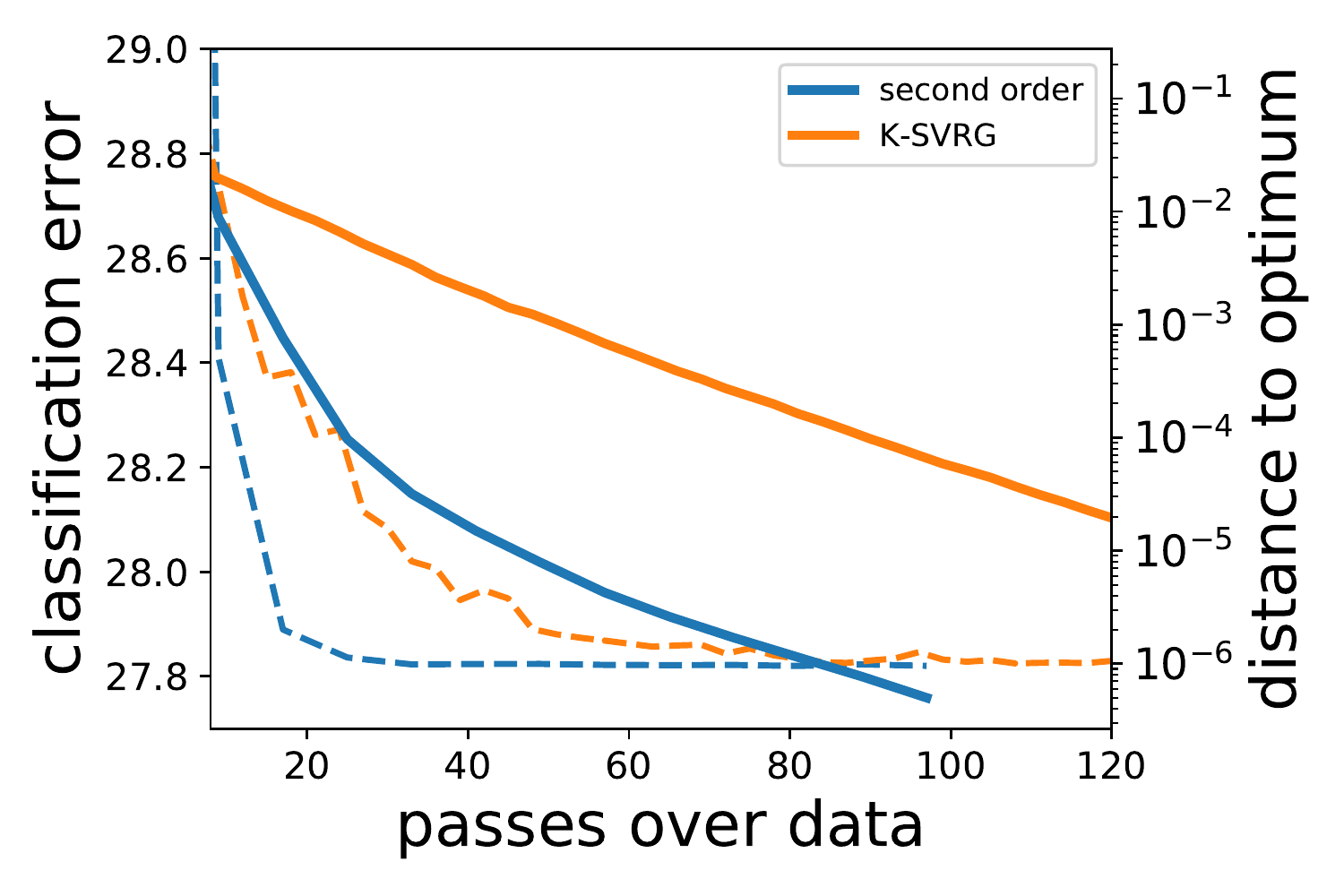}
    \caption{Training loss and test error as as function of the number of passes on the data for our algorithm vs. K-SVRG. on the \textbf{(left)} Susy and \textbf{(right)} Higgs data sets.}
    \label{fig:time-comparison}
\end{figure}
To find reasonable $\la$ for supervised learning applications, we cross-validated $\la$ finding the minimum test error at $\la = 10^{-10}$ for Susy and $\la = 10^{-9}$ for Higgs (see \cref{fig:test_susy,fig:test_higgs} in \cref{sec:exp-appendix}) for such values our algorithm and the competitor achieve an error of 19.5\% on the test set for Susy, comparable to the state of the art (19.6\% \cite{Rudi17}) and analogously for Higgs (see \cref{sec:exp-appendix}). We then used such $\la$'s as regularization parameters and compared our algorithm with a well known accelerated stochastic gradient technique {\em Katyusha SVRG} (K-SVRG)~\cite{Allen2017}, tailored to our problem using mini batches. In \cref{fig:time-comparison} we show the convergence of the training loss and classification error with respect to the number of passes on the data, of our algorithm compared to K-SVRG. It is possible to note our algorithm is order of magnitude faster in achieving convergence, validating empirically the fact that the proposed algorithm scales as $O(n \Deff_\la)$ in learning settings, while accelerated first order methods go as $O( (n+ \sqrt{n L / \lambda})\Deff_\la)$. Moreover, as mentioned in the introduction, this highlights the fact that precise optimization is necessary to achieve a good performance in terms of test error. Finally, note that since a pass on the data is much more expensive for K-SVRG than for our second order method (see \cref{sec:exp-appendix} for details), the difference in computing time between the second order scheme and K-SVRG is even more in favour of our second order method (see \cref{fig:susy_compare_small,fig:higgs_compare} in \cref{sec:exp-appendix}).

\subsection*{Acknowledgments}

We acknowledge support from the European Research Council (grant SEQUOIA 724063).

\bibliography{biblio_arxiv.bib}

\newpage

\appendix

{\Huge{Organization of the Appendix}}

\begin{itemize}
    \item [\large{\textbf{\ref{app:sc}.}}]\hyperref[app:sc]{{\textbf{\large{Main results on generalized self-concordant functions}}}} \vspace{0.2cm} \\
Notations, definitions and basic results concerning generalized self-concordant functions.
    \vspace{0.1cm} 
    \item [\large{\textbf{\ref{app:res-APN}.}}] \hyperref[app:res-APN]{{\textbf{\large{Results on approximate Newton methods}}}} \vspace{0.2cm} \\
    In this section, the interaction between the notion of Dikin ellipsoid, approximate Newton methods and generalized self-concordant functions is studied. The results needed in the main paper are all concentrated in \cref{app:newton_particular}. In particular the results in \cref{lm:dikin-is-cool} are proven in a more general form in \cref{lm:dikin-is-cool_2}. 
    \vspace{0.1cm} 
    \item [\large{\textbf{\ref{app:dec-lambda}.}}] \hyperref[app:dec-lambda]{\textbf{\large{Proof of bounds for the globalization scheme}}}
    \vspace{0.2cm} \\
    In this section, we leverage the results of the previous two sections to analyze the globalization scheme.
    \begin{itemize}
        \item[{\textbf{\ref{app:dec_la_tech}.}}] \hyperref[app:dec_la_tech]{\textbf{{Main technical lemmas}}} \\ We start by proving the result on the inclusion of Dikin ellipsoids (\cref{lm:dec_la}).
        \item[{\textbf{\ref{app:maintheorems}.}}] \hyperref[app:maintheorems]{\textbf{{Proof of main theorems}}} \\
        In particular, a general version of \cref{thm:main-thm} is proven. Moreover \cref{rm:cond-dependence} is proven in \cref{prp:adap_var}, while the fixed scheme to choose $(q_k)_{k \in \N}$ is proven in \cref{prp:fixed_q}.
        \item[{\textbf{\ref{app:proof_main_thm}.}}] \hyperref[app:proof_main_thm]{\textbf{{Proof of Thm. 1}}}\\ 
        Finally, we prove the properties of the globalization schemes presented in \cref{thm:main-thm}.
    \end{itemize}
    \vspace{0.1cm} 
    \item [\large{\textbf{\ref{app:kernels}.}}] \hyperref[app:kernels]{\textbf{\large{Non-parametric learning with generalized self-concordant functions}}} \vspace{0.2cm} \\
     In this section, some basic results about non-parametric learning with generalized self-concordant functions are recalled and the main results of \cref{sec:kernels} are proven.
    \begin{itemize}
        \item [{\textbf{\ref{app:kern-setting}.}}] \hyperref[app:kern-setting]{\textbf{{General setting and assumptions, statistical result for regularized ERM.}}}\\
        More details about the generalization properties of empirical risk minimization as well as the optimal rates in \cref{thm:opt-rates-ERM-text} are recalled.
        \item [{\textbf{\ref{app:kern-rand-proj}.}}] \hyperref[app:kern-rand-proj]{\textbf{{Reducing the dimension: projecting on a subspace using Nystr\"om sub-sampling.}}}
        \item [{\textbf{\ref{app:kern-ny-examples}.}}] \hyperref[app:kern-ny-examples]{\textbf{{Sub-sampling techniques}}}. \\
        The basics of uniform sub-sampling and sub-sampling with approximate leverage scores are recalled.
        \item
        [{\textbf{\ref{app:selecting_nys}.}}] \hyperref[app:selecting_nys]{\textbf{{Selecting the $M$ Nystr\"om points}}}\\
        \cref{thm:opt-rates-ny-text} is proven in a more general version in \cref{thm:nystrom_estimator_bounds}.
        \item[\textbf{\ref{app:approximate_b}}]\hyperref[app:approximate_b]{\textbf{Performing the globalization scheme to approximate $\beta_{M,\la}$}}\\
        A general scheme is proposed to solve the projected problem approximately using the globalization scheme.
        \begin{itemize}             \item[{\textbf{\ref{app:ans_kernels}.}}] \hyperref[app:ans_kernels]{\textbf{{Performing approximate Newton steps}}} \\
        We start by describing the way of computing approximate Newton steps.
        A generalized version of \cref{lm:pcg-appr-tech} is proven in \cref{prp:approx_ns_kernels}.
            \item[{\textbf{\ref{sssec:gen_scheme_optim}.}}] \hyperref[sssec:gen_scheme_optim]{\textbf{{Applying the globalization scheme to control $\hnd_{M,\la}(\beta)$}}}\\
            We then completely analyse the approximating of $\beta_{M,\la}$ from an optimization point of view (see \cref{prp:gen_scheme_optim}).
        \end{itemize}
        \item [{\textbf{\ref{app:main_kernels}.}}] \hyperref[app:main_kernels]{\textbf{{Final algorithm and results}}}\\
        Finally, the proof of \cref{thm:final-alg-kernel-text} is provided, using the results of the previous subsections.
    \end{itemize}
    \vspace{0.1cm}
    \item [\large{\textbf{\ref{sec:algorithm}.}}] \hyperref[sec:algorithm]{\textbf{\large{Algorithm}}}
    \vspace{0.2cm} \\
    In this section, the pseudocode for the algorithm presented in \cref{sec:kernels} and analyzed in \cref{thm:main_kernels} is provided.
    \vspace{0.1cm} 
    \item [\large{\textbf{\ref{sec:exp-appendix}.}}] \hyperref[sec:exp-appendix]{\textbf{\large{Experiments}}}
    \vspace{0.2cm} \\
    In this section, more details about the experiments are provided.
    \vspace{0.1cm} 
    \item [\large{\textbf{\ref{app:proj}.}}] \hyperref[app:proj]{\textbf{\large{Solving a projected problem to reduce dimension}}}
    \vspace{0.2cm} \\
    In this section, more details about the problem of randomized projections are provided.
    \begin{itemize}
        \item [{\textbf{\ref{app:proj-to-ideal}.}}] \hyperref[app:proj-to-ideal]{\textbf{{Relating the projected to the original problem}}}\\
        In particular, results to relate the ERM with the projected ERM in terms of excess risk are provided for generalized self-concordant functions.
    \end{itemize}
    \vspace{0.1cm}
    \item [\large{\textbf{\ref{app:proj_stat}.}}] \hyperref[app:proj_stat]{\textbf{\large{Relations between statistical problems and empirical problem.}}}
    \vspace{0.2cm} \\
   In this section, we provide results to relate excess expected risk with excess empirical risk for generalized self-concordant functions.
    \vspace{0.1cm}
    \item [\large{\textbf{\ref{sec:lemmas-operators}.}}] \hyperref[sec:lemmas-operators]{\textbf{\large{Multiplicative approximations for Hermitian operators}}}
    \vspace{0.2cm} \\
    In this section, some general analytic results on multiplicative approximations for Hermitian operators are derived. Moreover they are used to provide a simplified proof for the results in \cref{lm:appr-linear-system}. See in particular \cref{lm:sub_unif,lm:sub_ny} and \cite{pilanci2017newton}, Lemma~2.
    \vspace{0.1cm}
\end{itemize}
\vspace{1cm}

\section{\label{app:sc}Main results on generalized self-concordant functions}

In this section, we start by introducing a few notations. We define the key notion of generalized self-concordance in \cref{app:sc_2}, and present the main results concerning generalized self-concordant functions. In \cref{app:sc_3}, we describe how generalized self-concordance behaves with respect to an expectation or to certain relaxations.

\paragraph{Notations}

Let $\la \geq 0$ and $\mathbf{A}$ be a bounded positive semidefinite Hermitian operator on $\hh$. We denote with $\mathbf{I}$ the identity operator, and 
\eqal{\|x\|_{\mathbf{A}} &:= \|\mathbf{A}^{1/2} x\|,\\
\mathbf{A}_\la &:= \mathbf{A} + \la \mathbf{I}.}

Let $f$ be a twice differentiable convex function on a Hilbert space $\hh$. We adopt the following notation for the Hessian of $f$: 

\[\forall x \in \hh,~\Hess_f(x):= \nabla^2 f(x) \in \mathcal{L}(\hh).\]

For any $\la > 0$, we define the $\la$-regularization of $f$:
\[f_\la:= f + \frac{\la}{2}\|\cdot\|^2.\]
 $f_\la$ is $\la$-strongly convex and has a unique minimizer which we denote with $x^{f,\la}_{\star}$. Moreover, define
 \[\forall x \in \hh,~\Hess_{f,\la}(x):= \nabla^2 f_\la(x) = \Hess_f(x) + \la \Id,\qquad \nd_{f,\la}(x):= \|\nabla f_\la(x)\|_{\Hess_{f,\la}^{-1}(x)}.\]
The quantity $\nd_{f,\la}(x)$ is called the \textbf{Newton decrement} at point $x$ and will play a significant role.

When the function $f$ is clear from the context, we will omit the subscripts with $f$ and use $\Hess,\Hess_\la,\nd_\la ...$.

\subsection{Definitions and results on generalized self-concordant functions \label{app:sc_2}}

In this section, we introduce the main definitions and results for self-concordant functions. These results are mainly the same as in appendix B of \cite{marteau2019}.

\bd[generalized self-concordant function]\label{df:genscc}
Let $\hh$ be a Hilbert space. Formally, a generalized self-concordant function on $\hh$ is a couple $(f,\G)$  where:
\begin{enumerate}[label = \roman*]
\item $\G$ is a bounded subset of $\hh$; we will usually denote $\|\G\|$ or $\Rad$ the quantity $\sup_{g \in \G}\|g\| < \infty$;
\item $f$ is a convex and three times differentiable mapping on $\hh$ such that
\[\forall x \in \hh,~ \forall h,k \in \hh,~ \nabla^{(3)}f(x)[h,k,k] \leq \sup_{g \in \G}|g \cdot h| ~ \nabla^2 f(x)[k,k].\]

\end{enumerate}
\ed

To make notations lighter, we will often omit $\G$ from the notations and simply say that $f$ stands both for the mapping and the couple $(f,\G)$.

\bd[Definitions]\label{df:first_definitions_sc}
Let $f$ be a generalized self-concordant function. We define the following quantities. 
\begin{itemize}
    \item $\forall h \in \hh,~\tn_f(h):= \sup_{g \in \G}|h \cdot g|$;
    \item $\forall x \in \hh,~ \forall \la > 0,~ \radd_{f,\la}(x):= \frac{1}{\sup_{g \in \G}\|g\|_{\Hess^{-1}_{f,\la}(x)}}$;
    \item $\forall \cc \geq 0,~\forall \la > 0,~ \dik_{f,\la}(\cc):= \left\{x~:~\nd_{f,\la}(x) \leq \cc \radd_{f,\la}(x)\right\}$.
\end{itemize}
We also define the following functions: 
\eqal{
\label{eq:def_fun}
\psi(t) = \frac{e^t - t  -1 }{t^2},\quad \dikin(t) = \frac{1 - e^{-t}}{t},\quad \dikins(t) = \frac{e^t -1}{t}.
}
\ed 

Note that $\psi$, $\dikins$ are increasing functions and that $\dikin$ is a decreasing function. Moreover, $\frac{\dikins(t)}{\dikin(t)} = e^t$. Once again, if $f$ is clear, we will often omit the reference to $f$ in the quantities above, keeping only $\tn,\radd_\la,\dik_\la...$

We condense results obtained in \cite{marteau2019} under a slightly different form. The proofs, however, are exactly the same. 

While in \cite{marteau2019}, only the regularized case is dealt with, the proof techniques are exactly the same to obtain \cref{prp:thmcool1}. \cref{prp:thmcool2} is proved explicitly in Proposition 4 of \cite{marteau2019} and \cref{lm:loc} is proved in Proposition 5. \\\
Omitting the subscript $f$, we get the following results.

\bp[Bounds for the non-regularized function $f$]\label{prp:thmcool1} Let $f$ be a generalized self-concordant function. Then the following bounds hold (we omit $f$ in the subscripts):
\eqal{
\label{eq:h}
\forall x \in \hh,~ \forall h \in \hh,~ e^{-\tn(h)}\Hess(x) \preceq \Hess(x + h) \preceq e^{\tn(h)}\Hess(x),
}
\eqal{
\label{eq:g}
\forall x,h \in \hh,~\forall \la > 0,~ \|\nabla f(x + h) - \nabla f(x)\|_{\Hess^{-1}_\la(x)} \leq \dikins(\tn(h)) \| h \|_{\Hess_\la(x)},
}
\eqal{
\label{eq:fv}
\forall x,h \in \hh,~ \psi(-\tn(h)) \|h\|^2_{\Hess(x)} \leq f(x + h) - f(x) - \nabla f(x).h \leq \psi(\tn(h)) \|h\|^2_{\Hess(x)}.
}

\ep

We get the analoguous bounds in the regularized case.

\bp[Bounds for the regularized function $f_\la$] \label{prp:thmcool2}
Let $f$ be a generalized self-concordant function and $\la >0$ be a regularizer. Then the following bounds hold:
\eqal{
\label{eq:hl}
\forall x,h \in \hh,~ e^{-\tn(h)}\Hess_\la(x) \preceq \Hess_{\la}(x + h) \preceq e^{\tn(h)}\Hess_\la(x),
}
\eqal{
\label{eq:gl}
\forall x,h \in \hh,~ \dikin(\tn(h)) \| h \|_{\Hess_\la(x)} \leq  \|\nabla f_\la(x + h) - \nabla f_\la(x)\|_{\Hess^{-1}_\la(x)} \leq \dikins(\tn(h)) \| h \|_{\Hess_\la(x)},
}
\eqal{
\label{eq:fvl}
\forall x,h \in \hh,~ \psi(-\tn(h)) \|h\|^2_{\Hess_\la(x)} \leq f_{\la}(x + h) - f_\la(x) - \nabla f_\la(x).h \leq \psi(\tn(h)) \|h\|^2_{\Hess_\la(x)}.
}

\ep 

\bcor
Let $f$ be a $\G$ generalized self-concordant function and $\la >0$ be a regularizer, and $\xla$ the unique minimizer of $f_\la$. Then the following bounds hold for any $x \in \hh$:
\eqal{
\label{eq:go}
\dikin(\tn(x - \xla)) \| x - \xla \|_{\Hess_\la(x)} \leq  \underbrace{\|\nabla f_\la(x)\|_{\Hess^{-1}_\la(x)}}_{\nd_\la(x)} \leq \dikins(\tn(x - \xla)) \| x - \xla \|_{\Hess_\la(x)}
,}
\eqal{
\label{eq:fvo}
 \psi(-\tn(x - \xla)) \|x - \xla\|^2_{\Hess_\la(\xla)} \leq f_{\la}(x) - f_\la(\xla) \leq \psi(\tn(x - \xla)) \|x - \xla\|^2_{\Hess_\la(\xla)}
.}
\ecor

Moreover, the following localization lemma holds.

\blm[localization]\label{lm:loc}
Let $\la > 0$ be fixed. If $\frac{\nd_\la(x)}{\radd_\la(x)} <1$, then 
\eqal{
\label{eq:loc}
\tn(x - \xla) \leq -\log\left(1-\frac{\nd_\la(x)}{\radd_\la(x)} \right).
}
In particular, this shows:
\[\forall \cc <1,~\forall \la >0,~ x \in \dik_\la(\cc) \implies \tn(x - \xla) \leq - \log(1 - \cc). \]
\elm

We now state a Lemma which shows that the difference to the optimum in function values is equivalent to the squared newton decrement in a small Dikin ellipsoid. We will use this result in the main paper.

\blm[Equivalence of norms]\label{lm:equiv_norms}
Let $\la >0$ and $x \in \dik_\la(\frac{1}{7})$. Then the following holds:
\[\frac{1}{4}\nd_\la(x)^2 \leq f_\la(x) - f_\la(\xla)\leq \nd_\la(x)^2.\]
\elm 

\bpr 
Apply \cref{lm:loc} knowing $x \in \dik_\la(\frac{1}{7})$ to get $\nm{x - \xla} \leq \log(7/6)$. Then apply \cref{eq:fvl} and \cref{eq:gl} to get:
\eqals{
 f_\la(x) - f_\la(\xla) &\leq \psi(\nm{x - \xla})\|x-\xla\|_{\Hess_\la(\xla)}^2\\
&\leq e^{\nm{x - \xla}}\psi(\nm{x - \xla})\|x-\xla\|_{\Hess_\la(x)}^2 \\
&\leq  \frac{e^{\nm{x - \xla}}\psi(\nm{x - \xla})}{\dikin(\nm{x-\xla})^2}\nd_\la(x)^2.  }
Replacing with the bound above, we get
\[\forall \la >0,~ \forall x \in \dik_\la(\frac{1}{7}),~  f_\la(x) - f_\la(\xla) \leq \nd_\la(x)^2. \]
For the lower bound, proceed in exactly the same way.
\epr 

\subsection{Comparison between generalized self-concordant functions\label{app:sc_3}}

The following result is straightforward. 

\blm[Comparison between generalized self-concordant functions]
Let $\G_1 \subset \G_2 \subset \hh$ be two bounded subsets. If $(f,\G_1)$ is  generalized self-concordant, then $(f,\G_2)$ is also generalized self-concordant. Moreover, 
\[\forall x \in \hh,~ \forall \la > 0,~ \radd_{(f,\G_1),\la}(x) \geq \radd_{(f,\G_2),\la}(x).\]
\elm

In particular, we will often use the following fact. If $(f,\G)$ is generalized self-concordant, and $\G$ is bounded by $\Rad$, then $(f,B_{\hh}(\Rad))$ is also generalized self-concordant. Moreover,  
\[\radd_{(f,B_{\hh}(\Rad)),\la}(x) = \frac{\sqrt{\la + \la_{\min}(\Hess_f(x))}}{\Rad} \geq \frac{\sqrt{\la}}{\Rad}.\]

We now state a result which shows that, given a family of generalized self-concordant functions, the expectancy of that family is also generalized self-concordant. This can be seen as a reformulation of Proposition 2 of \cite{marteau2019}. 

\bp[Expectation]\label{prp:expectancy}
Let $\Z$ be a polish space equipped with its Borel sigma-algebra, and $\hh$ be a Hilbert space. Let $((f_z,\G_z))_{z \in \Z}$ be a family of generalized self-concordant functions such that the mapping $(z,x) \mapsto f_z(x)$ is measurable.

Assume we are given a random variable $Z$ on $\Z$, whose support we denote with $\supp(Z)$, such that 

\begin{itemize}
    \item the random variables $\|f_Z(0)\|, \|\nabla f_Z(0)\|, \Tr(\nabla^2 f_Z(0))$ are are bounded;
    \item $\G:= \bigcup_{z \in \supp(Z)}\G_z$ is a bounded subset of $\hh$.
\end{itemize}

Then the mapping $f: x \in \hh \mapsto \Exp{f_Z(x)}$ is well defined, $(f,\G)$ is  generalized self-concordant, and we can differentiate under the expectation.

\ep

\bcor
Let $n \in \N$ and $(f_i,\G_i)_{1\leq i\leq n}$ be a family of generalized self-concordant functions. Define
\[f(x) = \frac{1}{n}\sum_{i=1}^n{f_i(x)},~\G = \bigcup_{i=1}^n{\G_i}.\]
Then $(f,\G)$ is generalized self-concordant.

\ecor 

\newpage

\section{\label{app:res-APN}Results on approximate Newton methods}

In this section, we assume we are given a generalized self-concordant function $f$ in the sense of \cref{app:sc}. As $f$ will be fixed throughout this part, we will omit it from the notations. Recall the definitions from \cref{df:first_definitions_sc}:
\[\nd_\la(x):= \|\nabla f_\la(x)\|_{\Hess_\la^{-1}(x)},~\frac{1}{\rla(x)}:= \sup_{g\in \G}{\|g\|_{\Hess_\la^{-1}(x)}},~\dik_\la(\cc):= \left\{x~:~\frac{\nd_\la(x)}{\rla(x)}\leq \cc \right\}. \]

Define the following quantities: 

\begin{itemize}
    \item the true Newton step at point $x$ for the $\la$-regularized problem:
    $$\ns_\la(x):= \Hess_\la^{-1}(x)\nabla f_\la(x).$$
    \item the renormalized Newton decrement $\ndt_\la(x)$:
    \[\ndt_\la(x):= \frac{\nd_\la(x)}{\rla(x)}.\]
\end{itemize}

Moreover, note that a direct application of \cref{eq:hl} yields the following equation which relates the radii at different points:

\eqal{
\label{eq:radii relation} 
\forall \la >0,~ \forall x \in \hh,~ \forall h \in \hh,~ e^{-\nm{h}} \rla(x) \leq \rla(x+h) \leq e^{\nm{h}}\rla(x).
}

In this appendix, we develop a complete analysis of so-called approximate Newton methods in the case of generalized self-concordant losses. By "approximate Newton method", we mean that instead of performing the classical update $x_{t+1} = x_t - \ns_\la(x_t)$, we perform an update of the form $x_{t+1 } = x_t - \nsa_t$ where $\nsa_t$ is an approximation of the real Newton step. We will characterize this approximation by measuring its distance to the real Newton step using two parameters $\err$ and $\era$:
\[\|\nsa_t - \ns_\la(x_t)\| \leq \err \nd_\la(x_t) + \era.\]
We start by presenting a few technical results in \cref{app:newton_tech}. We continue by proving that an approximate Newton method has linear convergence guarantees in the right Dikin ellipsoid in \cref{app:newton_gen}. In \cref{app:newton_particular}, we adapt these results to a certain way of computing approximate Newton steps, which will be the one we use in the core of the paper. In \cref{app:newton_better}, we mention ways to reduce the computational burden of these methods by showing that since all Hessians are equivalent in Dikin ellipsoids, one can actually sketch the Hessian at one given point in that ellipsoid instead of re-sketching it at each Newton step. For the sake of simplicity, this is not mentioned in the core paper, but works very well in practice. 

\subsection{Main technical results \label{app:newton_tech}}

We start with a technical decomposition of the Newton decrement at point $x - \nsa$ for a given $\nsa \in \hh$. 

\blm[Technical decomposition] Let $\la >0$, $x \in \hh$ be fixed.
Assume we perform a step of the form $x - \nsa$ for a certain $\nsa \in \hh$. 
Define

$$\delta := \|\nsa - \ns_\la(x)\|_{\Hess_\la(x)},\qquad \widetilde{\delta}:= \frac{\delta}{\rla(x)}.$$ 
The following holds:
\eqal{\ndt_\la(x - \nsa) &\leq e^{\ndt_\la(x) + \widetilde{\delta}}\left[\psi(\ndt_\la(x) + \widetilde{\delta})(\ndt_\la(x) + \widetilde{\delta})^2 + \widetilde{\delta}\right]; \\
\nd_\la(x - \nsa_\la(x)) &\leq e^{\ndt_\la(x) + \widetilde{\delta}}\left[\psi(\ndt_\la(x) + \widetilde{\delta})(\ndt_\la(x) + \widetilde{\delta})(\nd_\la(x) + \delta) + \delta\right].}
\elm

\bpr

Note that by definition, $\nabla f_\la(x) = \Hess_\la(x) \ns_\la(x)$. Hence

\eqals{\|\nabla f^\la(x - \nsa)\|_{\Hess_\la^{-1}(x)} & = \|\nabla f^\la(x - \nsa) - \nabla f^\la(x) + \Hess_\la(x)\ns_\la(x)\|_{\Hess_\la^{-1}(x)}\\
&\leq \|\nabla f^\la(x - \nsa) - \nabla f^\la(x) + \Hess_\la(x)\nsa\|_{\Hess_\la^{-1}(x)} \\
&+ \|\Hess_\la(x)(\ns_\la(x)-\nsa)\|_{\Hess_\la^{-1}(x)} \\
&= \|\int_{0}^1{[\Hess_\la(x - s\nsa) -\Hess_\la(x)]\nsa ds} \|_{\Hess_\la^{-1}(x)} + \delta \\
& \leq \int_0^1{\|\Hess_\la^{-1/2}(x)\Hess_\la(x - s\nsa)\Hess_\la^{-1/2}(x) - \Id\|ds}~\|\nsa\|_{\Hess_\la(x)} + \delta.
}

Now using \cref{eq:hl}, one has $\|\Hess_\la^{-1/2}(x)\Hess_\la(x - s\nsa)\Hess_\la^{-1/2}(x) - \Id\| \leq e^{s\nm{\nsa}}-1$, whose integral on $s$ is $\psi(\nm{\nsa})\nm{\nsa}$ where $\psi$ is defined in \cref{df:first_definitions_sc}. Morever, bounding 
\[\|\nsa\|_{\Hess_\la(x)} \leq \|\nsa - \ns_\la(x)\|_{\Hess_\la(x)} + \|\ns_\la(x)\|_{\Hess_\la(x)} = \delta + \nd_\la(x),\]
it holds
\[\|\nabla f^\la(x - \nsa)\|_{\Hess_\la^{-1}(x)} \leq  \psi(\nm{\nsa}) \nm{\nsa}~(\nd_\la(x) + \delta) + \delta .\]

\paragraph{1. }
Now note that using \cref{eq:hl}, it holds:
$\nd_\la(x-\nsa) \leq e^{\nm{\nsa}/2}\|\nabla f^\la(x - \nsa)\|_{\Hess_\la^{-1}(x)}$ and hence: 
\eqal{\label{eq:lm11}\nd_\la(x-\nsa) \leq e^{\nm{\nsa}/2}\left(\psi(\nm{\nsa}) \nm{\nsa}~(\nd_\la(x) + \delta) + \delta \right) .}

\paragraph{2. } Moreover, using \cref{eq:radii relation},

\eqal{\label{eq:lm12}
\ndt_\la(x-\nsa) \leq e^{\nm{\nsa}}\left(\psi(\nm{\nsa}) \nm{\nsa}~(\ndt_\la(x) + \widetilde{\delta}) + \widetilde{\delta} \right)
.}

Noting that 
$$\nm{\nsa} \leq \frac{\|\nsa\|_{\Hess_\la(x)}}{\rla(x)} \leq  \ndt_\la(x) + \widetilde{\delta},$$
and bounding \cref{eq:lm11} simply by taking $e^{\nm{\nsa}/2} \leq e^{\nm{\nsa}}$, we get the two bounds in the lemma.

\epr

We now place ourselves in the case where we are given an approximation of the Newton step of the following form. Assume $\la$ and $x$ are fixed, and that we approximate $\ns_\la(x)$ with $\nsa$ such that there exists $\err \geq 0$ and $\era \geq 0$ such that it holds:
\[\|\nsa - \ns_\la(x)\|_{\Hess_\la(x)} \leq \err \nd_\la(x) + \era.\]
We define/prove the three different following regimes. 

\blm[3 regimes]\label{lm:3regimes}
Let $x \in \dik_\la\left(\frac{1}{7}\right)$ and $\la >0$ be fixed. Let
\[ 0 \leq \err \leq \frac{1}{7},~\era \geq 0 \text{ s.t. } \tilde{\varepsilon}_0:= \frac{\era}{\rla(x)} \leq \frac{1}{21}.\]
Let $\nsa$ be an approximation of the Newton steps satisfying $\|\nsa - \ns_\la(x)\|_{\Hess_\la(x)} \leq \err \nd_\la(x) + \era$. The three following regimes appear.
\begin{itemize}
    \item If $\ndt_\la(x) \geq \err$ and $\ndt_\la(x)^2 \geq \tilde{\varepsilon}_0$, then we are in the \textbf{quadratic regime}, i.e.
 \[\frac{10\ndt_\la(x - \nsa_\la(x))}{3} \leq \left(\frac{10\ndt_\la(x)}{3}\right) ^2 ,~ \nd_\la(x- \nsa_\la(x)) \leq \frac{10}{3} \ndt_\la(x) \nd_\la(x).\]
 \item If $ \err \geq \ndt_\la(x) $ and $\err \ndt_\la(x) \geq \erat$, then we are in the \textbf{linear regime}, i.e.  
  \[\frac{10}{3}\ndt_\la(x - \nsa_\la(x)) \leq \left(\frac{10\err}{3}\right)\left(\frac{10}{3} \ndt_\la(x)\right) ,~ \nd_\la(x- \nsa_\la(x)) \leq \frac{10}{3} \ndt_\la(x) \nd_\la(x).\]
  \item If $\erat \geq \ndt_\la(x)^2, \err ~ \ndt_\la(x) $, then the \textbf{maximal precision} of the approximation is reached, and it holds: 
    \[\ndt_\la(x - \nsa_\la(x)) \leq 3 \erat \leq \frac{1}{7} ,~ \nd_\la(x- \nsa_\la(x)) \leq 3 \era .\]
\end{itemize}
\elm

\bpr
Using the previous lemma, 
\eqals{\ndt_\la(x - \nsa_\la(x)) &\leq e^{(1 + \err)\ndt_\la(x) + \erat}\left[\psi((1+\err)\ndt_\la(x) + \erat)((1 + \err)\ndt_\la(x) + \erat)^2 + \err \ndt_\la(x) + \erat\right]\\
& \leq \square_1(\ndt_\la(x),\err,\erat)~\ndt_\la(x)^2 + \square_2(\ndt_\la(x),\err,\erat)~\err \ndt_\la(x) + \square_3(\ndt_\la(x),\err,\erat)~\erat
,}
and 
\[\nd_\la(x - \nsa_\la(x)) \leq \square_1(\ndt_\la(x),\err,\erat)~\ndt_\la(x) \nd_\la(x) + \square_2(\ndt_\la(x),\err,\erat)~\err \nd_\la(x) + \square_3(\ndt_\la(x),\err,\erat)~\era,\]
where the following defintions are used: 
\eqals{\square_1(\ndt,\err,\erat) &:= e^{(1 + \err)\ndt + \erat}\psi((1 + \err)\ndt + \erat)(1 + \err)^2,\\
\square_2(\ndt,\err,\erat) &:= e^{(1 + \err)\ndt + \erat},\\
\square_3(\ndt,\err,\erat) &:= e^{(1 + \err)\ndt + \erat}\left[2 \psi((1 + \err)\ndt + \erat)(1 + \err)\ndt + 1\right].}
Now assume $\erat \leq \frac{1}{21}$, $\ndt_\la(x), \err \leq \frac{1}{7}$.
Replacing these values in the functions above bounds $\square_1,\square_2$ and $\square_3$, and using the case distinction, we get the result.
\epr

\subsection{General analysis of an approximate Newton method \label{app:newton_gen}}

The following proposition describes the behavior of an approximate newton method where $\err$ and $\era$ are fixed a priori. 

\bp[General approximate Newton scheme results]\label{prp:gen_anm}
Let $\cc \leq \frac{1}{7}$ be fixed and $x_0 \in \dik_\la(\cc)$ be a given starting point. \\
Let $\err \leq \frac{1}{7}$ and $\era$ such that $\era \leq \frac{\cc}{4}~\rla(x_0)$. \\
Define the following approximate Newton scheme: 

\[\forall t \geq 0,~ x_{t+1} = x_t - \nsa_t,~\qquad\|\nsa_t - \ns_\la(x_t)\|_{\Hess_\la(x_t)} \leq \err \nd_\la(x_t) + \era.\]
The following guarantees hold.
\begin{itemize}
    \item $\forall t \geq 0,~ x_t \in \dik_\la(\cc)$.
    \item Let $t_c = \left\lfloor \log_2 \log_2 \frac{3}{10 \err}\right\rfloor + 1$.
    \[\forall t \leq t_c,~ \frac{10 \ndt_\la(x_t)}{3} \leq \max\left(\frac{12 \era}{\rla(x_0)}, 2^{-2^t}\right),\]
    \[\forall t \geq t_c,~ \frac{10 \ndt_\la(x_t)}{3} \leq \max\left(\frac{12 \era}{\rla(x_0)}, \left(\frac{10\err}{3}\right)^{t-t_c  +1 }\right).\]
    \item We can bound the relative decrease for both the Newton decrement and the renormalized Newton decrement:
    \begin{align*}
        &\forall t \leq t_c,~  &&\nd_\la(x_t) \leq \max\left(3 \era, \left(\frac{1}{2}\right)^{2^t-1}\nd_\la(x_0)\right),\\
        & &&\ndt_\la(x_t) \leq \max\left(\frac{18 \era}{5\rla(x_0)}, \left(\frac{1}{2}\right)^{2^t-1}\ndt_\la(x_0)\right).\\
    &\forall t \geq t_c,~ &&\nd_\la(x_t) \leq \max\left(3 \era, \left(\frac{10\err}{3}\right)^{t-t_c  +1 }\nd_\la(x_0)\right),\\
    &&&\ndt_\la(x_t) \leq \max\left(\frac{18 \era}{5 \rla(x_0)}, \left(\frac{10\err}{3}\right)^{t-t_c  +1 }\ndt_\la(x_0)\right).
        \end{align*}
\end{itemize}
\ep

\bpr
Start by noting, using \cref{eq:radii relation},

\eqal{\label{eq:cond_eps}\forall x \in \dik_\la\left(\frac{1}{7}\right),~ \varepsilon \leq \frac{\rla(x)}{21},~ \frac{6}{7}\rla(x_0)\leq  \rla(x) \leq \frac{7}{6}\rla(x_0).}
In particular, this holds for any $x \in \dik_\la(\cc),~ \cc \leq \frac{1}{7}$.
Thus,
\[\forall \cc \leq \frac{1}{7},~ \forall x_0 \in \dik_\la(\cc),~\frac{\era}{\rla(x_0)} \leq \frac{\cc}{4} \implies \forall x \in \dik_\la(\cc),~\frac{\era}{\rla(x)} \leq \frac{\cc}{3}.\]
\paragraph{1.} Proving the first point is simple by induction.
Indeed, assume $\ndt_\la(x_t) \leq \cc$. We can apply \cref{lm:3regimes} since the conditions on $\varepsilon$ and $\err$ guarantee that the conditions of this lemma are satisfied. 

If we are in either the linear or quadratic regime, the fact that $\frac{10\err}{3},\frac{10 \ndt_\la(x_t)}{3} \leq \frac{10}{21}$ show that $\ndt_\la(x_{t+1}) \leq \frac{10}{21} \ndt_\la(x_t) \leq \cc$.

If we are in the last case, $\ndt_\la(x_{t+1}) \leq \frac{3 \era}{\rla(x_t)} \leq \cc$.

\paragraph{2.} Let us prove the second bullet point by induction. 
Start by assuming the property holds at $t$. By the previous point, the hypothesis of \cref{lm:3regimes} are satisfied at $x_t$ with $\err$ and $\varepsilon$. Assume we are in the limiting case; we easily show that in this case, 
\[\frac{10 \ndt_\la(x_{t+1})}{3} \leq \frac{10}{3} ~3 \frac{\era}{\rla(x_t)} \leq \frac{35 \era}{3\rla(x_0)}.\]
Here, the last inequality comes from \cref{eq:cond_eps}. 
If we are not in the limiting case, let us distinguish between the two following cases. \\\

If  $t \leq t_c - 1$,
\begin{align*}\frac{10 \ndt_\la(x_{t+1})}{3} &\leq \frac{10 \ndt_\la(x_{t})}{3} \max\left(\frac{10 \ndt_\la(x_{t})}{3},\frac{10 \err}{3}\right) \\
&\leq \max \left(\frac{35 \era}{3\rla(x_0)},\frac{10 \ndt_\la(x_{t})}{3} \max\left(\left(\frac{1}{2}\right)^{2^t},\frac{10 \err}{3}\right)\right),\end{align*}
where the last inequality comes from using the induction hypothesis and the fact that $\frac{10 \ndt_\la(x_{t})}{3} \leq 1$. Using once again the induction hypotheses and the fact that $t \leq \left\lfloor \log_2 \log_{2} \frac{3}{10 \err}\right\rfloor$ which implies $\frac{10 \err }{3} \leq \left(\frac{1}{2}\right)^{2^t}$, we finally get 
\[\frac{10 \ndt_\la(x_{t+1})}{3} \leq \max \left(\frac{35 \era}{3 \rla(x_0)},\left(\frac{1}{2}\right)^{2^{t+1}}\right).\]

The fact that the second property holds for $t = t_c$ is trivial
Now consider the case where $t \geq t_c$. Using the same technique as before but noting that in this case
\[\frac{10 \ndt_\la(x_t)}{3} \leq \max\left(\frac{35 \era}{3 \rla(x_0)},\left(\frac{10 \err}{3}\right)^{t-t_c+1}\right) \leq \max\left(\frac{35 \era}{3 \rla(x_0)},\frac{10 \err}{3}\right) ,\]
We easily use \cref{lm:3regimes} to reach the desired conclusion. 

\paragraph{3.} Let $t < t_c$. Then using \cref{lm:3regimes}:
\[\forall s \leq t,~ \nd_\la(x_{s+1}) \leq \max\left(3 \era,\max(\frac{10\err}{3},\frac{10\ndt_\la(x_s)}{3}) \nd_\la(x_s)\right).\]
Using the fact that for any $s \leq t$, $\frac{10\ndt_\la(x_s)}{3} \leq \max(\frac{35 \era}{3\rla(x_0)},\left(\frac{1}{2}\right)^{2^s})$:
\[\forall s \leq t,~ \nd_\la(x_{s+1}) \leq \max\left(3 \era,\frac{35 \era}{3}\frac{\nd_\la(x_s)}{\rla(x_0)},\max(\frac{10\err}{3},\left(\frac{1}{2}\right)^{2^s}) \nd_\la(x_s)\right).\]
Now using the fact that for any $s \leq t,~ \ndt_\la(x_s) \leq \frac{1}{7}$, we see that $\frac{\nd_\la(x_s)}{\rla(x_0)} \leq \frac{7}{6}\ndt_\la(x_s) \leq \frac{1}{6}$ and hence $\frac{35 \era}{3}\frac{\nd_\la(x_s)}{\rla(x_0)} \leq 3 \era$. Moreover, since $s \leq t < t_c$, $\max(\frac{10\err}{3},\left(\frac{1}{2}\right)^{2^s}) = \left(\frac{1}{2}\right)^{2^s}$. Thus:
\[\forall s \leq t,~ \nd_\la(x_{s+1}) \leq \max\left(3 \era,\left(\frac{1}{2}\right)^{2^s} \nd_\la(x_s)\right).\]
Combining these results yields: 
\[\nd_\la(x_{t+1}) \leq \max\left(3 \era, \left(\frac{1}{2}\right)^{2^{t+1}-1} \nd_\la(x_0)\right).\]
This shows the first equation, that is: 

\[\forall t \leq t_c,~ \nd_\la(x_{t}) \leq \max\left(3 \era, \left(\frac{1}{2}\right)^{2^{t}-1} \nd_\la(x_0)\right).\]
The case for $t \geq t_c$ is completely analogous. We can also reproduce the same proof to get the same bounds for $\ndt$, since the bounds in \cref{lm:3regimes} are the same for both.

\epr

\subsection{Main results in the paper \label{app:newton_particular}}

In the main paper, we mention two types of Newton method. First, we present a result of convergence on the full Newton method:

\blm[Quadratic convergence of the full Newton method]\label{lm:full_newton}
Let $\cc \leq \frac{1}{7}$ and $x_0 \in \dik_\la(\cc)$. Define 
\[x_{t+1} =  x_t - \ns_\la(x_t).\]
Then this scheme converges quadratically, i.e.:
\[\forall t \in \N,~ \frac{\nd_\la(x_t)}{\nd_\la(x_0)},\frac{\ndt_\la(x_t)}{\ndt_\la(x_0)}\leq 2^{-(2^{t} - 1)}.\]
Thus :
\begin{itemize}
    \item $\forall t \in \N,~ x_t \in \dik_\la(\cc)$.
    \item For any $\tilde{\cc} \leq \cc$ then $\forall t \geq  \left\lceil \log_2 \left(1  + \log_2 \frac{\cc}{\tilde{\cc}}\right) \right\rceil,~ x_t \in \dik_\la(\tilde{\cc})$.
    \item For any $\varepsilon > 0$, $\forall t \geq  \left\lceil \log_2 \left(1  + \log_2 \frac{\nd_\la(x_0)}{\sqrt{\varepsilon}}\right) \right\rceil,~ \nd_\la(x_t) \leq \sqrt{\varepsilon},~ f_\la(x)-f_\la(\xla) \leq \varepsilon$.
    \item If we perform the Newton method and return  the first $x_t$ such that $\nd_\la(x_t) \leq \sqrt{\varepsilon}$, then the number of Newton steps computations is at most 
    $1 + \left\lceil \log_2 \left(1  + \log_2 \frac{\nd_\la(x_0)}{\sqrt{\varepsilon}}\right) \right\rceil$. 
\end{itemize}
\elm 

\bpr 

A full Newton method is an approximate Newton method where $\err,\era = 0$. Thus apply \cref{prp:gen_anm}; note that in this case $t_c = +\infty$. The last point shows that if $\cc \leq \frac{1}{7}$, and if we perform the Newton method with a full Newton step, then 
\[\forall t \geq 0,~ \ndt_\la(x_t) \leq 2^{-(2^t -1)} \nd_\la(x_0),~\ndt_\la(x_t) \leq 2^{-(2^t -1)} \nd_\la(x_0).\]
This shows the quadratic convergence, and the first two points directly follow. For the third point, the result for $\nd_\la(x_t)$ directly follows from the previous equation, and the one on function values is a direct consequence of \cref{lm:equiv_norms} and the fact that $x_t \in \dik_\la(1/7)$.\\\

For the last point, note that $\nd_t(x_t) = \nabla f_\la(x_t)\cdot \ns_\la(x_t)$ is accessible. Moreover, the bound on $t$ is given in the point before, and since one has to compute $\ns_\la(x_s)$ for $0 \leq s \leq t$, there are at most $t+1$ computations.
\epr 

In the main paper, we compute approximate Newton steps by considering methods which naturally yield only a relative error $\err$ and no absolute error $\era$. Indeed, we take the following notation.

\paragraph{Approximate solutions to linear problems.}
Let $\Ab$ be a positive definite Hermitian operator on $\hh$, $b$ in $\hh$, and a wanted relative precision $\err$.

We say that $x$ is a $\err$-relative approximation to the linear problem $\Ab x = b$ and write $x \in \lso(\Ab,b,\err)$ if the following holds:
\[\|\Ab^{-1}b - x\|_{\Ab} \leq \err \|b\|_{\Ab^{-1}} = \err \|\Ab^{-1} b\|_{\Ab}.\]
Note that if $x \in \lso(\Ab,b,\err)$ for $\err < 1$, then $$ (1-\err) \|b\|_{\Ab^{-1}} \leq x \cdot b \leq (1+\err) \|b\|_{\Ab^{-1}} . $$

The following lemma shows that if, instead of computing the exact Newton step, we compute a relative approximation of the Newton step belonging to $\lso(\Hess_\la(x),\nabla f_\la(x),\err)$ for a given $\err < 1$, then one has linear convergence. Moreover, we show that we can still perform a method which automatically stops.

\bp[relative approximate Newton method] \label{prp:approximate_newton_scheme}
Let $\la > 0$, $\err \leq \frac{1}{7}$, $\cc \leq \frac{1}{7}$ and a starting point $x_0 \in \dik_\la(\cc)$. 
Assume we perform the following Newton scheme:
\[\forall t \geq 0,~ x_{t+1} = x_t - \nsa_t,~\qquad \nsa_t \in \lso(\Hess_\la(x_t),\nabla f_\la(x_t),\err) .\]
Then the scheme converges linearly, i.e. 
\[\forall t \in \N,~ \frac{\nd_\la(x_t)}{\nd_\la(x_0)},\frac{\ndt_\la(x_t)}{\ndt_\la(x_0)} \leq 2^{-t}.\]
Thus,
\begin{itemize}
   \item $\forall t \in \N,~ x_t \in \dik_\la(\cc)$.
    \item For any $\tilde{\cc} \leq \cc$ then $\forall t \geq  \left\lceil \log_2  \frac{\cc}{\tilde{\cc}} \right\rceil,~ x_t \in \dik_\la(\tilde{\cc})$.
    \item For any $\varepsilon > 0$, $\forall t \geq  \left\lceil \log_2  \frac{\nd_\la(x_0)}{\sqrt{\varepsilon}} \right\rceil,~ \nd_\la(x_t) \leq \sqrt{\varepsilon},~ f_\la(x)-f_\la(\xla) \leq \varepsilon$
    \item If the method is performed and returns  the first $x_t$ such that $x_t \cdot \nsa_t \leq \frac{6}{7}\varepsilon$, then at most 
    $2 + \left\lfloor \log_2\left( \sqrt{\frac{4}{3}}\frac{\nd_\la(x_0)}{\sqrt{ \varepsilon}}\right) \right\rfloor$ approximate Newton steps computations have been performed, and $\nd_\la(x_t) \leq \sqrt{\varepsilon},~ f_\la(x)-f_\la(\xla) \leq \varepsilon$.
\end{itemize}
\ep

\bpr 
Apply \cref{prp:gen_anm} with $\era = 0$ and $\err = \frac{1}{7}$, since if $\err \leq \frac{1}{7}$, then a fortiori the approximation satisfies the condition for $\err = \frac{1}{7}$. The last point clearly states that
\[\forall t \in \N,~ \frac{\nd_\la(x_t)}{\nd_\la(x_0)},\frac{\ndt_\la(x_t)}{\ndt_\la(x_0)} \leq \left(\frac{10}{21}\right)^{t} \leq 2^{-t}.\]
From this, using \cref{lm:equiv_norms} for the third point, the first three points are easily proven. \\\
For the last point, note that since $\nsa_t \in \lso(\Hess_\la(x_t),\nabla f_\la(x_t),\err)$, the following holds: $\nabla f_\la(x_t) \cdot \nsa_t = \nd_\la(x_t)^2 + \nabla f_\la(x_t) \cdot \left(\nsa_t - \Hess_\la^{-1}(x_t) \nabla f_\la(x_t)  \right)$. 
Now bound 
\[|\nabla f_\la(x_t) \cdot\left( \nsa_t - \Hess_\la^{-1}(x_t) \nabla f_\la(x_t) \right)| \leq \nd_\la(x_t) ~\|\nsa_t - \Hess_\la^{-1}(x_t) \nabla f_\la(x_t)\|_{\Hess_\la(x_t)} \leq \err \nd_\la(x_t)^2.\]
Thus: 
\[(1-\err)\nd_\la(x_t)^2 \leq \nabla f_\la(x_t) \cdot \nsa_t \leq (1+\err)\nd_\la(x_t)^2.\]
Since $\rho = \frac{1}{7}$, we see that if $\nabla f_\la(x_t) \cdot \nsa_t \leq \frac{6}{7}\varepsilon$, then $\nd_\la(x_t)^2 \leq \varepsilon$. 
Moreover, since we stop at the first $t$ where $\nabla f_\la(x_t) \cdot \nsa_t \leq \frac{6}{7}\varepsilon$, then if $t$ denotes the time at which we stop,
\[ \frac{6}{7}\varepsilon < \nabla f_\la(x_{t-1}) \cdot \nsa_{t-1}  \leq \frac{8}{7} \nd_\la(x_{t-1})^2 \implies \nd_\la(x_{t-1}) ^2\geq \frac{3}{4}\varepsilon.\]
Since $\nd_\la(x_{t-1})^2 \leq 2^{-2(t-1)}\nd_\la(x_0)^2$, this implies in turn that $t-1 \leq \log_2\left( \sqrt{\frac{4}{3}}\frac{\nd_\la(x_0)}{\sqrt{ \varepsilon}}\right)$.
Thus, necessarily, $t \leq 1 + \left\lfloor \log_2\left( \sqrt{\frac{4}{3}}\frac{\nd_\la(x_0)}{\sqrt{ \varepsilon}}\right) \right\rfloor$, and since we compute approximate Newton steps for $s = 0,...,t$, we finally have that the number of approximate Newton steps is bounded by 
\[2 + \left\lfloor \log_2\left( \sqrt{\frac{4}{3}}\frac{\nd_\la(x_0)}{\sqrt{ \varepsilon}}\right) \right\rfloor.\]
\epr 

Last but not least, we summarize all these theorem in the following simple result.

\blm\label{lm:dikin-is-cool_2}
Let $\la > 0, \cc \leq 1/7$, let $f_\la$ be generalized self-concordant and $x \in \dik_\la(\cc)$. It holds:
   $
    \frac{1}{4}\nd_\la(x)^2 \leq f_\la(x) - f_\la(\xla) \leq \nd_\la(x)^2$.
Moreover, the full Newton method starting from $x_0$ has quadratic convergence, i.e. if $x_t$ is obtained via $t \in \N$ steps of the Newton method \cref{eq:newton-intro}, then
$ \nd_\la(x_t) \leq 2^{-(2^{t} - 1)} \nd_\la(x_0).$
Finally, the approximate Newton method starting from $x_0$ has linear convergence, i.e. if $x_t$ is obtained via $t \in \N$ steps of  \cref{eq:approx-newton-intro}, with $\nsa_t \in \lso(\Hess_\la(x_t),\nabla f_\la(x_t),\err)$ and $\err \leq 1/7$, then
$\nd_\la(x_t) \leq 2^{-t} \nd_\la(x_0).$
\elm

\bpr 
 The three points are obtained in the following lemmas, assuming $x \in \dik_\la(1/7)$.
 \begin{itemize}
 \item  For $\frac{1}{4}\nd_\la(x)^2 \leq f_\la(x) - f_\la(\xla) \leq \nd_\la(x)^2$, see \cref{lm:equiv_norms} in \cref{app:sc_2}.
 \item The convergence rate of the full Newton method starting in $\dik_\la(1/7)$ is obtained in \cref{lm:full_newton}.
 \item The convergence rate of the approximate Newton method starting in $\dik_\la(1/7)$ is obtained in \cref{prp:approximate_newton_scheme}. 
 \end{itemize}
\epr

\subsection{Sketching the Hessian only once in each Dikin ellispoid \label{app:newton_better}}

In this section, we provide a lemma which shows in essence that if we are in a small Dikin ellipsoid, then we can keep the Hessian of the starting point and compute approximations of $\Hess_\la^{-1}(x_0)\nabla f_\la(x_t)$; they will be good approximations to $\Hess_\la^{-1}(x_t)\nabla f_\la(x_t)$ as well.

\blm 
Let $\cc < 1$ and $x_0 \in \dik_\la(\cc)$ be fixed.

Let $\tHess$ be an approximation of the Hessian at $x_0$, approximation wich we quantify with 
\[t:= \|\Hess^{-1/2}_\la(x_0)\left(\Hess_\la(x_0) - \tHess\right)\Hess^{-1/2}_\la(x_0)\|.\]
Assume 
\[1+t < 2(1 - \cc)^2.\]
Let $b \in \hh$. If $\nsa \in \lso(\tHess_\la,b,\terr)$, then
$$\forall x \in \dik_\la(\cc),~ \nsa \in \lso(\Hess_\la(x),b,\err),~ \err = \frac{(\terr -1)(1-\cc)^2 + (1+t)}{2(1 - \cc)^2 - (1+t)}.$$
In particular, if $\cc \leq \frac{1}{30}$, $x_0 \in \dik_\la(\cc)$, 
\[\forall x \in \dik_\la(\cc),~\forall b \in \hh,~ \nsa \in \lso(\Hess_\la(x_0),b,\frac{1}{20}) \implies \nsa \in \lso(\Hess_\la(x),b,\frac{1}{7}).\]

\elm 

\bpr 
First, start with a general theoretical result. 
\paragraph{1. } Let $\Ab$ and $\Bb$ be two positive semi-definite hermitian operators. Let $\la > 0$, $b \in \hh$ and $\nsa \in \lso(\Bb_\la,b,\terr)$. 
Decompose 
\eqals{
\|\Ab_\la^{-1}b - \nsa\|_{\Ab_\la} &\leq \|\Ab_\la^{-1}b - \Bb_\la^{-1}b\|_{\Ab_\la}  + \|\Bb_\la^{-1}b - \nsa\|_{\Ab_\la} \\
&\leq \|\Ab_\la^{1/2}(\Ab_\la^{-1} - \Bb_\la^{-1})\Ab_\la^{1/2}\|~ \|b\|_{\Ab_\la^{-1}} + \|\Ab_\la^{1/2}\Bb_\la^{-1/2}\| ~\|\Bb_\la^{-1}b - \nsa\|_{\Bb_\la}.}

Now using the fact that $\Ab_\la^{-1} - \Bb_\la^{-1} = \Bb_\la^{-1}(\Bb - \Ab)\Ab_\la^{-1}$, 

\eqals{
\|\Ab_\la^{1/2}(\Ab_\la^{-1} - \Bb_\la^{-1})\Ab_\la^{1/2}\| &\leq \|\Ab_\la^{-1/2}(\Bb-\Ab)\Ab_\la^{-1/2}\|~\|\Ab_\la^{1/2}\Bb_\la^{-1} \Ab_\la^{1/2}\| \\
&=  \|\Ab_\la^{-1/2}(\Bb-\Ab)\Ab_\la^{-1/2}\|~\|\Ab_\la^{1/2}\Bb_\la^{-1/2} \|^2.
}
Moreover, 
\[\|\Bb_\la^{-1}b - \nsa\|_{\Bb_\la} \leq \terr \|b\|_{\Bb_\la^{-1}} \leq \|\Ab^{1/2}\Bb^{-1/2}\|~\|b\|_{\Ab_\la^{-1}}.\]
Putting things together, and noting that from \cref{lm:equiv_operators}, $\|\Ab^{1/2}\Bb^{-1/2}\|^2 \leq \frac{1}{1-\|\Ab_\la^{-1/2}(\Bb - \Ab)\Ab_\la^{-1/2}\|}$ as soon as $\|\Ab_\la^{-1/2}(\Bb - \Ab)\Ab_\la^{-1/2}\| < 1$, it holds: 
$$\nsa \in \lso(\Ab_\la,b,\err),~\err = \frac{\terr + \|\Ab_\la^{-1/2}(\Bb - \Ab)\Ab_\la^{-1/2}\|}{1 - \|\Ab_\la^{-1/2}(\Bb - \Ab)\Ab_\la^{-1/2}\|} .$$

The aim is now to apply this lemma to $\Ab = \Hess(x)$ and $\Bb = \tHess$.

\paragraph{2. } Let $x,x_0 \in \dik_\la(\cc)$. Using \cref{lm:comb_approx}, we see that 
\[1 + \|\Hess_\la^{-1/2}(x)(\tHess - \Hess(x))\Hess_\la^{-1/2}(x)\| \leq (1 + t)(1 + \|\Hess_\la^{-1/2}(x)(\Hess(x_0) - \Hess(x))\Hess_\la^{-1/2}(x)\| ).\]
Using \cref{eq:hl}, it holds:
\[ (e^{-\nm{x-x_0}}- 1)\Id \preceq \Hess_\la^{-1/2}(x)(\Hess(x_0) - \Hess(x))\Hess_\la^{-1/2}(x) \preceq (e^{\nm{x_0 - x}} - 1)\Id .\]
Thus, 
\[\|\Hess_\la^{-1/2}(x)(\Hess(x_0) - \Hess(x))\Hess_\la^{-1/2}(x)\| \leq \max(1-e^{-\nm{x-x_0}},e^{\nm{x-x_0}} - 1) = e^{\nm{x-x_0}} - 1 .\]
Finally, using the fact that $x_0,x \in \dik_\la(\cc)$ for $\cc < 1$ yields $\nm{x-x_0} \leq 2 \log \frac{1}{1-\cc}$. Hence
\[1 + \|\Hess_\la^{-1/2}(x)(\Hess(x_0) - \Hess(x))\Hess_\la^{-1/2}(x)\| \leq \frac{1}{(1-\cc)^2}.\]
Thus,
\[\|\Hess_\la^{-1/2}(x)(\tHess - \Hess(x))\Hess_\la^{-1/2}(x)\| \leq \frac{1 + t}{(1 - \cc)^2} - 1.\]
The result then follows.
\epr

\newpage

\section{Proof of bounds for the globalization scheme}\label{app:dec-lambda}

In this section, we prove that the scheme of decreasing $\mu$ towards $\la$ converges. 

\subsection{Main technical lemmas \label{app:dec_la_tech}}

\blm[Next $\mu$]\label{lm:next_la} Let $\mu >0$, $\cc <1$.
$$\nd_\mu(x) \leq \frac{\cc}{3}~\frac{\sqrt{\mu}}{\Rad} \implies \nd_{\tmu}(x) \leq \cc ~ \frac{\sqrt{\tmu}}{\Rad}, \qquad \tmu:=  q~\mu,\qquad q \geq \frac{\frac{1}{3} + \frac{\Rad \sqrt{\mu} \|x\|_{\Hess_\mu^{-1}(x)}}{\cc}}{1 + \frac{\Rad \sqrt{\mu} \|x\|_{\Hess_\mu^{-1}(x)}}{\cc}} .$$

$$ x\in \dik_\mu\left( \frac{\cc}{3}\right)  \implies x\in \dik_{\tmu}\left( \cc\right), \qquad \tmu:=  q~\mu,\qquad q \geq \frac{\frac{1}{3} + \frac{\mu \|x\|_{\Hess_\mu^{-1}(x)}}{\cc~ \radd_\mu(x)}}{1 + \frac{\mu \|x\|_{\Hess_\mu^{-1}(x)}}{\cc~ \radd_\mu(x)}} .$$

\elm

\bpr
For any $\tmu < \mu$, note that 
$$ \forall x \in \hh,~ \|\Hess_{\tmu}^{-1/2}(x)\Hess_{\mu}^{1/2}(x)\| = \sqrt{\frac{\la_{\min}(\Hess(x)) + \mu}{\la_{\min}(\Hess(x)) + \tmu}}\leq \sqrt{\mu/\tmu}.$$
 This shows that $\|\cdot\|_{\Hess_{\tmu}^{-1}(x)} \leq  \sqrt{\frac{\mu}{\tmu}}~\|\cdot\|_{\Hess_{\mu}^{-1}(x)}$, and in particular that $\frac{1}{\radd_{\tmu}(x)} \leq \sqrt{\mu/\tmu}\frac{1}{\radd_{\mu}(x)}$.
 
 Using this fact, it holds: 
\eqals{\ndt_{\tmu}(x) &= \frac{\|\nabla f_{\tmu}(x)\|_{\Hess^{-1}_{\tmu}(x)}}{\radd_{\tmu}(x)} \\
&= \frac{\|\nabla f_\mu(x) -(\mu - \tmu)x\|_{\Hess^{-1}_{\tmu}(x)}}{\radd_{\tmu}(x)} \\
& \leq \frac{\mu}{\tmu}~\frac{\|\nabla f_\mu(x)\|_{\Hess_{\mu}^{-1}(x)}}{\radd_\mu(x)}
 + \left(\frac{\mu}{\tmu}-1\right) \frac{\|\mu x\|_{\Hess^{-1}_\mu(x)}}{\radd_\mu(x)}.}
 
 Hence, if $\ndt_\mu(x) \leq \frac{\cc}{3}$, a condition to obtain $\ndt_{\tmu}(x) \leq \cc$ is the following:
 
 \[\frac{\mu}{\tmu} \left(\frac{\cc}{3} + t\right) \leq \cc + t \Leftrightarrow \tmu \geq \mu \frac{\cc/3 +  t}{ \cc +  t}\qquad t =  \frac{\|\mu x\|_{\Hess^{-1}_\mu(x)}}{\radd_\mu(x)}.\]
 This yields the second point of the lemma. The analysis is completely analoguous for the first.
\epr

\blm[Useful bounds for $q$]\label{lm:bounds_q}
Let $\mu >0$. Then the following hold:
\[ \forall x \in \hh,~ \frac{\mu \|x\|_{\Hess_\mu^{-1}(x)}}{\radd_\mu(x)} \leq \Rad \sqrt{\mu} \|x\|_{\Hess_\mu^{-1}(x)} \leq \Rad \|x\|
.\]
Moreover, we can bound all of these quantities using $\xmu$: 
\begin{itemize}
\item For any $\cc < 1$, $x \in \hh$, if $ x \in \dik_\mu(\cc/3)$, then the following holds:
\[\frac{\mu \|x\|_{\Hess^{-1}_\mu(x)}}{\cc ~ \radd_\mu(x)} \leq \frac{1}{3}\left(1 + \frac{1}{1-\cc/3}\right) + \frac{1}{1-\cc/3} \frac{\|\mu \xmu\|_{\Hess_\mu^{-1}(\xmu)}}{ \cc ~ \radd_\mu(\xmu)}. \]
\item For any $\cc < 1$, $x \in \hh$, if $\frac{\Rad \nd_\mu(x)}{\sqrt{\mu}} \leq \frac{\cc}{3}$, then the following holds:
\[\frac{\Rad \sqrt{\mu} \|x\|_{\Hess_\mu^{-1}(x)}}{\cc} \leq \left(1 + \frac{1}{1-\cc/3}\right)\frac{1}{3} + \sqrt{\frac{1}{1-\cc/3}} \frac{\Rad \sqrt{\mu} \|\xmu\|_{\Hess_\mu^{-1}(\xmu)}}{\cc}.\]
Likewise, it can be shown that under the same conditions: 
\[\frac{\Rad \|x\|}{\cc} \leq \frac{\Rad \|\xmu\|}{\cc} + \frac{1}{3}\dikins(-\log(1 - \cc/3)).\]
\end{itemize}
\elm

\bpr 
The first bound is obvious. Moreover, the fact that $\ndt_\mu(x) \leq \frac{\cc}{3}$ implies that $\nm{x - \xmu}  \leq \log\frac{1}{1-\cc/3}$. Thus, we get the classical bounds on the Hessian using \cref{eq:h}: 
\[e^{-\nm{x - \xmu}} \Hess(x) \preceq \Hess(\xmu) \preceq e^{\nm{x - \xmu}} \Hess(x).\]

\paragraph{1. Bound on $\mu \|x\|_{\Hess_\mu^{-1}(x)}$.}
 Using \cref{eq:hl,eq:gl},

\eqals{
\mu \|x\|_{\Hess^{-1}_\mu(x)} & = \|\nabla f_\mu(x) - \nabla f(x) + \nabla f(\xmu) - \nabla f(\xmu)\|_{\Hess^{-1}_\mu(x)}  \\
&\leq \nd_\mu(x) + \int_{0}^1{\|\Hess_\mu(x)^{-1/2}\Hess(x_t)(x - \xmu)\|~dt} + \|\nabla f(\xmu)\|_{\Hess_\mu(x)},~x_t = t x + (1-t)\xmu.}

Now bound $\|\Hess_\mu(x)^{-1/2}\Hess(x_t)(x - \xmu)\| \leq \|\Hess_\mu(x)^{-1/2}~\Hess_\mu(x_t)^{1/2}\|~\|x - \xmu\|_{\Hess(x_t)}$ and use \cref{eq:hl} and \cref{eq:h} to get:
\[\|\Hess_\mu(x)^{-1/2}\Hess(x_t)(x - \xmu)\| \leq e^{t~\nm{x-\xmu}} \|x - \xmu\|_{\Hess(x)}.\]
Integrating this yields:
\[\int_{0}^1{\|\Hess_\mu(x)^{-1/2}\Hess(x_t)(x - \xmu)\|~dt} \leq \dikins(\nm{x -\xmu})~ \|x - \xmu\|_{\Hess(x)} \leq e^{\nm{x-\xmu}}~\nd_\mu(x).\]
Where the last inequality is obtained using the bounds between gradient and hessian distance \cref{eq:gl}. Finally, using the bound on $\nm{x-\xmu}$,
\[\mu \|x\|_{\Hess_\mu^{-1}(x)} \leq \left(1 + \frac{1}{1-\cc/3}\right)\nd_\mu(x) + \sqrt{\frac{1}{1-\cc/3}}\|\nabla f(\xmu)\|_{\Hess_\mu^{-1}(\xmu)}.\]

\paragraph{2. Bound on $\Rad  \|x\|$.}
Start by decomposing 
\[\Rad \|x\| \leq \Rad \|\xmu\| + \Rad  \|x - \xmu\|.\]
Now bound 
\[\Rad  \|x - \xmu\| \leq \frac{\Rad}{\sqrt{\mu}} \|x - \xmu\|_{\Hess_\mu(x)}.\]
Using \cref{eq:hl}, $ \|x - \xmu\|_{\Hess_\mu(x)} \leq \dikins(-\log(1-\cc/3))\nd_\mu(x).$ Hence:
\[\Rad \|x\| \leq \Rad \|\xmu\| + \dikins(-\log(1-\cc/3)) \frac{\Rad \nd_\mu(x)}{\sqrt{\mu}} .\]

\paragraph{3. Now assume $x \in \dik_\mu(\cc/3)$.} Using the bound on $\mu \|x\|_{\Hess^{-1}_\mu(x)}$, and noting that 
\[\frac{1}{\radd_\mu(x)} \leq e^{\nm{x-\xmu}/2} \frac{1}{\radd_\mu(\xmu)},\]
it holds:
\[\frac{\mu \|x\|_{\Hess^{-1}_\mu(x)}}{\cc ~ \radd_\mu(x)} \leq \frac{1}{3}\left(1 + \frac{1}{1-\cc/3}\right) + \frac{1}{1-\cc/3} \frac{\|\mu \xmu\|_{\Hess_\mu^{-1}(\xmu)}}{ \cc ~ \radd_\mu(\xmu)}.\]
\paragraph{4. Now assume $\frac{\Rad \nd_\mu(x)}{\sqrt{\mu}} \leq \frac{\cc}{3}$.}. We know that in particular, $x \in \dik_\mu(\cc/3)$ and hence:
\eqals{\Rad \sqrt{\mu} \|x\|_{\Hess_\mu^{-1}(x)} &\leq \left(1 + \frac{1}{1-\cc/3}\right)\frac{\Rad \nd_\mu(x)}{\sqrt{\mu}} + \sqrt{\frac{1}{1-\cc/3}}\frac{\Rad \mu \|\xmu\|_{\Hess_\mu^{-1}(\xmu)}}{\sqrt{\mu}}\\
&\leq \left(1 + \frac{1}{1-\cc/3}\right)\frac{\cc}{3} + \sqrt{\frac{1}{1-\cc/3}} \Rad \sqrt{\mu} \|\xmu\|_{\Hess_\mu^{-1}(\xmu)}.}
Hence 
\[\frac{\Rad \sqrt{\mu} \|x\|_{\Hess_\mu^{-1}(x)}}{\cc} \leq \left(1 + \frac{1}{1-\cc/3}\right)\frac{1}{3} + \sqrt{\frac{1}{1-\cc/3}} \frac{\Rad \sqrt{\mu} \|\xmu\|_{\Hess_\mu^{-1}(\xmu)}}{\cc}.\]
Likewise:  
\[\frac{\Rad \|x\|}{\cc} \leq \frac{\Rad \|\xmu\|}{\cc} + \frac{1}{3}\dikins(-\log(1 - \cc/3)).\]
\epr

We can get the following simpler bounds.

\bcor[Application to $\cc = \frac{1}{7}$]\label{cor:bounds_q} Applying \cref{lm:bounds_q} to $\cc = \frac{1}{7}$, we get the following bounds.
Let $\mu >0$.
\begin{itemize}
\item For any $x \in \hh$, if $ x \in \dik_\mu(\cc/3)$, then the following holds:
\[\frac{7 \mu \|x\|_{\Hess^{-1}_\mu(x)}}{ \radd_\mu(x)} \leq 1 + \frac{8\|\mu \xmu\|_{\Hess_\mu^{-1}(\xmu)}}{  \radd_\mu(\xmu)}. \]
\item For any $\cc < 1$, $x \in \hh$, if $\frac{\Rad \nd_\mu(x)}{\sqrt{\mu}} \leq \frac{\cc}{3}$, then the following hold:
\[ 7\Rad \sqrt{\mu} \|x\|_{\Hess_\mu^{-1}(x)} \leq 1 +  8 \Rad \sqrt{\mu} \|\xmu\|_{\Hess_\mu^{-1}(\xmu)}.\]
\[7\Rad \|x\| \leq 7\Rad \|\xmu\| + 1.\]
\end{itemize}
\ecor

\subsection{Proof of main theorems \label{app:maintheorems}}

In this section, we bound the number of iterations of our scheme in different cases.

Recall the proposed globalization scheme in the paper, where $\texttt{ANM}_\rho(f,x,t)$ is a method performing $t$ successive $\rho$-relative approximate Newton steps of $f$ starting at $x$.

\fbox{
\begin{minipage}[t]{0.9 \textwidth}
\textbf{Proposed Globalization Scheme}
\begin{center}
\textit{Phase I: Getting in the Dikin ellispoid of $f_\la$}
\end{center}
Start with $x_0 \in \hh, \mu_0 > 0$, $t, T \in \N$ and $(q_k)_{k \in \N} \in (0,1]$.\\
For $k \in \N$\\
${}\qquad x_{k+1} \leftarrow \texttt{ANM}_\rho(f_{\mu_k},x_{k},t)$\\
${}\qquad \mu_{k+1} \leftarrow q_{k+1}\mu_{k}$\\
Stop when $\mu_{k+1} < \la$ and set $x_{last}\leftarrow x_k$.
$K \leftarrow k$
\begin{center}
\textit{Phase II: reach a certain precision starting from  inside the Dikin ellipsoid}
\end{center}
Return $\widehat{x} \leftarrow \texttt{ANM}_\rho(f_{\la},x_{last},T)$
\end{minipage}
}

\vspace*{0.3cm}

Throughout this section, we will denote with $K$ the value of $k$ when the scheme stops, i.e. the first value of $k$ such that $\mu_{k+1} < \la$. 

\paragraph{Adaptive methods } We start by presenting an adaptive way to select $\mu_{k+1}$ from $\mu_k$, with theoretical guarantees. The main result is the following.

\bp[Adaptive, simple version]\label{prp:adapt}
Assume that we perform phase I starting at $x_0$ such that 
\[\frac{\Rad \nd_{\mu_0}(x_0)}{\sqrt{\mu_0}} \leq \frac{1}{7}.\]
Assume that at each step $k$, we compute $x_{k+1}$ using $t = 2$ iterations of the $\rho$-relative approximate Newton method. 
Then if at each iteration, we set: 
\[\mu_{k+1} = q_{k+1} ~\mu_{k},\qquad q_{k+1}:= \frac{\frac{1}{3} + 7\Rad \|x_{k+1}\|}{1 + 7\Rad \|x_{k+1}\|} .\]
Then the following hold:
\paragraph{1.} $\forall k \leq K+1,~\frac{\Rad \nd_{\mu_k}(x_k)}{\sqrt{\mu_k}} \leq \frac{1}{7}$ .
\paragraph{2.} The decreasing parameter $q_{k+1}$ is bounded above before reaching $K$:
$$\forall k \leq K,~ q_{k+1} \leq  \frac{\frac{4}{3} + 7\Rad \|x^\star_{\mu_k}\|}{2 + 7\Rad \|x^\star_{\mu_k}\|} \leq  \frac{\frac{4}{3} + 7\Rad \|x^\star_{\la}\|}{2 + 7\Rad \|x^\star_{\la}\|}.$$
\paragraph{3. } $K$ is finite, 
$$ K \leq \left\lfloor \frac{\log \frac{\mu_0}{\la}}{\log \frac{2 + 7\Rad \|x^\star_{\la}\|}{\frac{4}{3} + 7\Rad \|x^\star_{\la}\|}}\right\rfloor \leq \left\lfloor \left(3 + 11\Rad \|\xla\|\right)\log \frac{\mu_0}{\la}\right\rfloor,$$ 
and $\frac{\Rad \nd_\la(x_{K+1})}{\sqrt{\la}} \leq \frac{1}{7} $.
\ep 

\bpr 

Let us prove the three points one by one.
\paragraph{1.} This is easily proved by induction, the keys to the induction hypothesis being:
\begin{itemize}
    \item  Using the induction hypothesis, $x_k \in \dik_{\mu_k}(\cc)$ and hence, using \cref{prp:approximate_newton_scheme} shows that after two iterations of the approximate Newton scheme, $\frac{\nd_{\mu_k}(x_{k+1})}{\nd_{\mu_k}(x_k)}\leq \frac{1}{3}$ which implies $\frac{\Rad \nd_{\mu_k}(x_{k+1})}{\sqrt{\mu_k}} \leq \frac{\cc}{3}$.
    \item Now using \cref{lm:next_la}, we see that that since 
    $$7 \Rad \|x_{k+1}\| = \frac{\Rad \|x_{k+1}\|}{\cc} \geq  \frac{\Rad \sqrt{\mu_k} \|x_{k+1}\|_{\Hess^{-1}_{\mu_k}(x_{k+1})}}{\cc },$$
    the hypotheses to guarantee the bound for $q_{k+1}$ hold, hence 
    \[\frac{\Rad \nd_{\mu_{k+1}}(x_{k+1})}{\sqrt{\mu_{k+1}}} \leq \cc.\]
\end{itemize}
\paragraph{2.} Using the second bullet point of \cref{cor:bounds_q}, we see that the previous point implies
\[\forall k \leq K,~ 7\Rad \|x_{k+1}\| \leq 7 \Rad \|x^\star_{\mu_k}\| + 1 \implies q_{k+1} \leq \frac{4/3 + 7 \Rad \|x^\star_{\mu_k}\|}{2 + 7 \Rad \|x^\star_{\mu_k}\|}.\]
Now using the fact that for any $k \leq K$, $\mu_k > \la$, we can use the simple fact that $\|x^\star_\la\| \geq \|x^\star_{\mu_k}\|$ to get the desired bound for $q_{k+1}$. 
\paragraph{3.} Using the previous point clearly shows the following bound: 
\[\forall k \leq K+1,~\mu_{k} \leq  \left(\frac{\frac{4}{3} + 7\Rad \|x^\star_{\la}\|}{2 + 7\Rad \|x^\star_{\la}\|}\right)^{k}\mu_0.\]
As this clearly converges to $0$ when $k$ goes to infinity, $K$ is necessarily finite. Applying this for $k = K$, we see that:
\[\la \leq \mu_K \leq \left(\frac{\frac{4}{3} + 7\Rad \|x^\star_{\la}\|}{2 + 7\Rad \|x^\star_{\la}\|}\right)^{K}\mu_0.\]
This shows that $K \leq \frac{\log \frac{\mu_0}{\la}}{\log\frac{2 + 7\Rad \|x^\star_{\la}\|}{\frac{4}{3} + 7\Rad \|x^\star_{\la}\|} }$.

The final bound is obtained noting that 
\[\frac{2 + 7\Rad \|x^\star_{\la}\|}{\frac{4}{3} + 7\Rad \|x^\star_{\la}\|} = 1 + \frac{1}{t},\qquad t = 2 + \frac{21}{2}\Rad \|x^\star_{\la}\|,\]
and using the classical bound:
\[\frac{1}{\log(1 + \frac{1}{t})} \leq t+1.\]
Finally, the fact that $\frac{\Rad\nd_\la(x_{K+1})}{\sqrt{\la}} \leq \cc$ is just a consequence of the fact that $\mu_{K+1} \leq \la \leq \mu_K$ and thus that $\la = q \mu_K$ with $q \geq q_{K+1}$, which is shown to satisfy the condition in \cref{lm:next_la}. Hence, the lemma holds not only for $\mu_{K+1}$ but also for $\la$.
\epr

\br[$\mu_0$] \label{rmk:mu0}
In the previous proposition, we assume start at $x_0,\mu_0$ such that 
\[\frac{\Rad \nd_{\mu_0}(x_0)}{\sqrt{\mu_0}} \leq \frac{1}{7}.\]
A simple way to have such a pair is simply to select:
\[x_0 = 0,~ \mu_0 = 7\Rad \| \nabla f(0) \|,\]
since $\frac{\Rad \nd_{\mu_0}(x_0)}{\sqrt{\mu_0}}  =\frac{\Rad \|\nabla f(0)\|_{\Hess_{\mu_0}^{-1}(0)}}{\sqrt{\mu_0}} \leq \frac{\Rad \|\nabla f(0)\|}{\mu_0} $.
\er

Alternatively, if one can approximately compute $\|x\|_{\Hess^{-1}_\mu(x)}$, one can propose the following variant, whose proof is completely analogous.

\bp[Adaptive, small variant version]\label{prp:adap_var}
Assume that we perform phase I starting at $x_0$ such that 
\[\frac{\Rad \nd_{\mu_0}(x)}{\sqrt{\mu_0}} \leq \frac{1}{7}.\]
Then if at each iteration, we set: 
\[t_{k+1} = 7\sqrt{\frac{7}{6}}\Rad \sqrt{\mu_k }\sqrt{x_{k+1}\cdot s_{k+1}},s_{k+1} \in \lso(\Hess_{\mu_k}(x_{k+1}),x_{k+1},\frac{1}{7}),\]
and 
\[\mu_{k+1} = q_{k+1} ~\mu_{k},\qquad q_{k+1}:= \frac{\frac{1}{3} + t_{k+1} }{1 + t_{k+1}} .\]
Then the following hold:
\paragraph{1.} $\forall k \leq K,~\frac{\Rad \nd_{\mu_k}(x_k)}{\sqrt{\mu_k}} \leq \frac{1}{7}$.
\paragraph{2.} The decreasing parameter $q_{k+1}$ is bounded above before reaching $K$:
\[\forall k \leq K,~q_{k+1} \leq \sup_{\mu_0 \geq \mu \geq \la}\frac{\frac{7}{3} + 10 \Rad \sqrt{\mu}\|x^\star_{\mu}\|_{\Hess^{-1}_{\mu}(x^\star_{\mu})}}{3 + 10\Rad \sqrt{\mu}\|x^\star_{\mu}\|_{\Hess^{-1}_{\mu}(x^\star_{\mu})}} \leq \frac{\frac{7}{3} + 10 \Rad \|\xla\|}{3 + 10\Rad \|\xla\|}.\]
\paragraph{3. } $K$ is finite, 
$$ K   \leq \left(\frac{9}{2} + 15\sup_{\la \leq \mu \leq \mu_0}\Rad \sqrt{\mu}\|x^\star_{\mu}\|_{\Hess^{-1}_{\mu}(x^\star_{\mu})}\right)\log \frac{\mu_0}{\la},$$ 
and $\frac{\Rad \nd_\la(x_{K+1})}{\sqrt{\la}} \leq \frac{1}{7} $.
\ep 

\bpr 
The main thing to note is that because of the properties of $\frac{1}{7}$-approximations, if $s_{k+1} \in \lso(\Hess_{\mu_k}(x_{k+1}),x_{k+1},\frac{1}{7})$,
\[(1-\frac{1}{7})\|x_{k+1}\|_{\Hess^{-1}_{\mu_k}(x_{k+1})}^2 \leq x_{k+1}\cdot s_{k+1} \leq (1+\frac{1}{7})\|x_{k+1}\|_{\Hess^{-1}_{\mu_k}(x_{k+1})}^2 .\]
Hence,
\[\|x_{k+1}\|_{\Hess^{-1}_{\mu_k}(x_{k+1})} \leq \sqrt{\frac{7}{6}}\sqrt{x_{k+1}\cdot s_{k+1}} \leq \sqrt{\frac{4}{3}}\|x_{k+1}\|_{\Hess^{-1}_{\mu_k}(x_{k+1})}.\]
Hence, $t_{k+1} \geq 7 \Rad \sqrt{\mu_k}\|x_{k+1}\|_{\Hess^{-1}_{\mu_k}(x_{k+1})}$, and we can apply \cref{lm:next_la} to get the first point.\\\
To get the second point, we bound $t_{k+1}$ above: 
\[t_{k+1} \leq 7 \sqrt{\frac{4}{3}} \Rad \sqrt{\mu_k}\|x_{k+1}\|_{\Hess^{-1}_{\mu_k}(x_{k+1})}.\]
Now use \cref{cor:bounds_q} to find:
\[t_{k+1} \leq  \sqrt{\frac{4}{3}}\left( 1 + 8\Rad \sqrt{\mu_k}\|x^\star_{\mu_k}\|_{\Hess^{-1}_{\mu_k}(x^\star_{\mu_k})}\right) \leq 2 + 10 \Rad \sqrt{\mu_k}\|x^\star_{\mu_k}\|_{\Hess^{-1}_{\mu_k}(x^\star_{\mu_k})}.\]
Thus, 
\[q_{k+1} \leq \frac{\frac{7}{3} + 10 \Rad \sqrt{\mu_k}\|x^\star_{\mu_k}\|_{\Hess^{-1}_{\mu_k}(x^\star_{\mu_k})}}{3 + 10\Rad \sqrt{\mu_k}\|x^\star_{\mu_k}\|_{\Hess^{-1}_{\mu_k}(x^\star_{\mu_k})}}.\]
Note that as long as $k \geq K$, 
\[q_{k+1} \leq \sup_{\mu \geq \la}\frac{\frac{7}{3} + 10 \Rad \sqrt{\mu}\|x^\star_{\mu}\|_{\Hess^{-1}_{\mu}(x^\star_{\mu})}}{3 + 10\Rad \sqrt{\mu}\|x^\star_{\mu}\|_{\Hess^{-1}_{\mu}(x^\star_{\mu})}} \leq \frac{\frac{7}{3} + 10 \Rad \|\xla\|}{3 + 10\Rad \|\xla\|}.\]
This guarantees convergence. \\\

For the last point, the proof is exactly the same as in the previous proposition. 

\epr

\paragraph{General non-adaptive result.}

As mentioned in the core of the article, in practice, we do not select $q_{k+1}$ at each iteration using a safe adaptative value, but rather decrease $\mu_{k+1} = q \mu_k$ with a constant $q$, which we see as a parameter to tune. The following result shows that for $q$ large enough, this is justified, and that the lower bound we get for $q$ depends on the radius of the Dikin ellipsoid $\radd_\mu(x)$, instead of $\frac{\sqrt{\mu}}{\Rad}$ in the previous bounds, which is somewhat finer and shows that if the data is structured such that this radius is very big, then $q$ might actually be very small.

\bp[Fixed $q$]\label{prp:fixed_q}
Assume that we perform phase I starting at $x_0$ such that 
\[x_0 \in \dik_{\mu_0}(\frac{1}{7}).\]
Assume we perform the method with a fixed $q_{k+1} = q$, satisfying
\[q \geq \sup_{\la \leq \mu \leq \mu_0}\frac{\frac{4}{3} + 8\frac{\mu \|x^\star_{\mu}\|_{\Hess_{\mu}^{-1}(x^\star_{\mu})}}{\radd_{\mu}(x^\star_{\mu})}}{2 + 8\frac{\mu \|x^\star_{\mu}\|_{\Hess_{\mu}^{-1}(x^\star_{\mu})}}{\radd_{\mu}(x^\star_{\mu})}}.\]
Then the following hold:
\paragraph{1.} $\forall k \leq K+1,~ x_k \in \dik_{\mu_k}(\frac{1}{7})$.
\paragraph{2. } $K$ is finite, 
$$ K   \leq \frac{1}{1-q}\log \frac{\mu_0}{\la}, $$ 
and $x_{K+1} \in \dik_{\la}(\frac{1}{7}) $.
\ep 

\bpr 
Let us prove the two points.

\paragraph{1.} Let us prove the result by induction. The initialization is trivial. Now assume $x_k \in \dik_{\mu_k}(\frac{1}{7})$. Performing two iterations of the approximate Newton method guarantees that 
\[x_{k+1} \in \dik_{\mu_k}(\frac{1}{21}),\]
as show in \cref{prp:approximate_newton_scheme}.
Now using \cref{lm:next_la}, we see that $x_{k+1} \in \dik_{q \mu_k}(\frac{1}{7})$, provided that 

\[q \geq \frac{\frac{1}{3} + \frac{7 \mu_k \|x_{k+1}\|_{\Hess_{\mu_k}^{-1}(x_{k+1})}}{\radd_{\mu_k}(x_{k+1})}}{1 + \frac{7 \mu_k \|x_{k+1}\|_{\Hess_{\mu_k}^{-1}(x_{k+1})}}{\radd_{\mu_k}(x_{k+1})}} .\]

Now using \cref{cor:bounds_q}, we get that 
\[\frac{7 \mu_k \|x_{k+1}\|_{\Hess_{\mu_k}^{-1}(x_{k+1})}}{\radd_{\mu_k}(x_{k+1})} \leq 1 + \frac{8\mu_k \|x^\star_{\mu_k}\|_{\Hess_{\mu_k}^{-1}(x^\star_{\mu_k})}}{\radd_{\mu_k}(x^\star_{\mu_k})} \leq 1 + 8 \sup_{\la \leq \mu \leq \mu_0} \frac{\mu \|x^\star_{\mu}\|_{\Hess_{\mu}^{-1}(x^\star_{\mu})}}{\radd_{\mu}(x^\star_{\mu})}.\]
Hence the result. 

\paragraph{2.} This point just follows, using the bound $\frac{1}{\log \frac{1}{q}} \leq \frac{1}{1-q}$.

\epr 

\subsection{Proof of \cref{thm:easy_first_phase} \label{app:proof_main_thm}}

Using \cref{rmk:mu0}, the fact that $x_0 = 0$ and  $\mu_0 = 7 \Rad \|\nabla f(0)\|$, as well as the hypotheses of the theorem, we can apply \cref{prp:adapt}, and show that the number of steps $K$ performed in the first phase is bounded: 

\[K \leq \left\lfloor (3 + 11 \Rad \|\xla\|)\log(7\Rad \|\nabla f(0)\|/\la)\right\rfloor.\]
Moreover, this proposition also shows that $\Rad \nd_\la(x_{last})/\sqrt{\la} \leq \frac{1}{7}$. Hence, we can use \cref{prp:approximate_newton_scheme}: if 
\[t\geq T = \left\lceil \log_2\sqrt{\frac{\la \varepsilon^{-1}}{\Rad^2}} \right\rceil  \geq \left\lceil \log_2\frac{\nd_\la(x_{last})}{\sqrt{\varepsilon}} \right\rceil, \] 
then it holds $\nd_\la(\hat{x}) \leq \sqrt{\varepsilon}$ and $f_\la(\hat{x}) - f_\la(\xla) \leq \varepsilon$. 
\qed{}

\newpage

\section{Non-parametric learning with generalized self-concordant functions\label{app:kernels}}
In this section, the aim is to provide a fast algorithm in the case of Kernel methods which achieves the optimal statistical guarantees. 

\subsection{General setting and assumptions, statistical result for regularized ERM.\label{app:kern-setting}}

In this section, we consider the supervised learning problem of learning a predictor $f:\X \rightarrow \Y$ from training samples $(x_i,y_i)_{1\leq i \leq n}$ which we assume to be realisations from a certain random variable $ Z = (X,Y) \in \Z = \X\times \Y $ whose distribution is $\rho$.  In what follows, for simplification purposes, we assume $\Y = \R$; however, this analysis can easily be adapted (although with heavier notations) to the setting where $\Y = \R^p$. Our aim is to compute the predictor of minimal generalization error 
\eqal{\label{pb:learning}\inf_{f \in \hh}{L(f):= \Expb{z \sim \rho}{\ell_z(f(x))}},}
where $\hh$ is a space of candidate solutions and $\ell_z:  \R \rightarrow \R$ is a loss function comparing the prediction $f(x)$ to the objective $y$.  

\paragraph{Kernel methods.}

Kernel methods consider a space of functions $\hh_K$ implicitly constructed from a symmetric positive semi-definite Kernel  $K: \X \times \X \rightarrow$ and whose basic functions are the $K_x: t \in \X \mapsto K(x,t)$ and the linear combinations of such functions 
$f = \sum_{j=1}^m{\alpha_j K_{x_j}}$. \\

It is endowed with a scalar product such that: $\forall x_1,x_2 \in \X,~ K_{x_1}\cdot K_{x_2} = K(x_1,x_2)$, and as a consequence, $\hh_K$ satisfies the self-reprocucing property:
 $$ \forall x \in \X,~ \forall f \in \hh,~f(x) = \langle f,K_x\rangle_{\hh}.$$  In order to find a good predictor for \cref{pb:learning}, the following estimator, called the regularized ERM estimator, is often computed:
$$\hf_\la:= \argmin_{f \in \hh} \hL_\la(f):= \frac{1}{n}\sum_{i=1}^n{\ell_{z_i}(f(x_i))} + \frac{\la}{2}\|f\|_{\hh}^2.$$

The properties of this estimator have been studied in \cite{devito} for the square loss and \cite{marteau2019} for generalized self-concordant functions. In \cref{app:proj_stat}, we recall the full setting of \cite{marteau2019}, and extend it to include the statistical properties of the projected problem.

\paragraph{Assumptions}

In this section, we will make the following assumptions, which are reformulations of the assumptions of \cite{marteau2019}, which we recall in \cref{app:proj_stat}, in order to have the statistical properties of the regularized ERM. First, we assume that the $(x_i,y_i)$ are i.i.d. samples.

\begin{restatable}[i.i.d. data]{ass}{asmiid1}
\label{asm:iid1}
The samples $(z_i)_{1 \leq i \leq n}=(x_i,y_i)_{1\leq i \leq n}\in\Z^n$ are independently and identically distributed according to $\rho$. 
\end{restatable}
In the case where $\Y = \R$, we make the following assumptions on the loss, which leads to the self concordance of the mappings $f \mapsto \ell_z(f(x))$ and that of $L$, $\hL$...

\begin{restatable}[Technical assumptions]{ass}{asmtech1}
\label{asm:asmtech1}
The mapping $(z,t) \in \Z \times \R \mapsto \ell_z(t)$ is measurable. Moreover, 
\begin{itemize}
    \item there exists $\Rell < \infty$ such that for all $z \in \supp(Z)$, 
    \[\forall t \in \R,~ |\ell_z^{(3)}(t)|\leq \Rell \ell_z^{\prime \prime}(t),\]
    \item the random variables $|\ell_Z(0)|, |\ell^{\prime}_Z(0)|, | \ell^{\prime \prime}_Z(0)|$ are are bounded;
    \
    \item The kernel is bounded, i.e. $\forall x \in \supp(X),~K(x,x) \leq \ka^2$ for a certain $\ka$. 
\end{itemize}
\end{restatable}

Using these assumptions, we see that the following properties are satisfied. Define $L_z(f): = \ell_z(f(x))$. Then the $L_z$ satisfy the following properties:
\begin{itemize}
    \item For any $z \in \Z$, $(L_z,\left\{\Rell K_{x}\right\})$ is a generalized self-concordant function in the sense of \cref{df:first_definitions_sc}.
    \item The mapping $(z,f) \in \Z \times \hh \mapsto L_z(f)$ is measurable;
    \item the random variables $\|L_Z(0)\|, \|\nabla L_Z(0)\|, \Tr(\nabla^2 L_Z(0))$ are bounded by $|\ell_Z(0)|$, $\ka|\ell^{\prime}_Z(0)|$, $\ka^2 |\ell_Z^{\prime \prime}(0)|$;
    \item $\G:= \left\{\Rell K_{x}~:~z \in \supp(Z)\right\}$ is a bounded subset of $\hh$, bounded by $\Rad = \Rell\ka$.
\end{itemize}

This shows that \cref{asm:gen_sc} and \cref{asm:asmtech} are satisfied by the $L_z$ and hence, using \cref{prp:scProblems} in the next appendix, $L$ is well-defined, generalized self-concordant with $\G$. Moreover, the empirical loss 
\[\hL = \frac{1}{n}\sum_{i=1}^n{L_{z_i}},\]
is also generalized self-concordant with $\hG:= \left\{\Rell K_{x_i}~:~1 \leq i \leq n\right\}$.\\\

Finally, as in \cref{app:proj_stat}, we make an assumption on the regularity of the problem; namely, we assume that a solution to the learning problem exists in $\hh$.

\begin{restatable}[Existence of a minimizer]{ass}{asmminimizer1}
\label{asm:minimizer1}
There exists $f^\star \in \hh$ such that $ L(f^\star) = \inf_{f \in \hh} L(f).$
\end{restatable}

We adopt all the notations from \cref{app:proj_stat}, doing the distinction between expected an empirical problems by adding a $\widehat{\cdot}$ over the quantities related to the empirical problem. We continue using the standard notations for $L$: for any $f \in \hh$ and $\la >0$,
\[L_\la(f) = L(f) + \frac{\la}{2}\|f\|^2,\qquad \hL_\la(f) = \hL(f) + \frac{\la}{2}\|f\|^2  \]
\[\Hess(f) = \nabla^2 L(f),\qquad \Hess_\la(f) = \nabla^2 L_\la(f) = \Hess(f) + \la \Id \]
\[\hHess(f) = \nabla^2 \hL(f),\qquad \hHess_\la(f) = \nabla^2 \hL_\la(f) = \hHess(f) + \la \Id \]
Recall that $\hf_\la$ is defined as the minimizer of $\hL_\la$. 

Define the following bounds on the second order derivatives: 
\[\forall f \in \hh,~ \btns(f) = \sup_{z \in \supp(Z)}\ell^{\prime \prime}_z(f(x)).\]
\paragraph{Statistical properties of the estimator}

The statistical properties of the estimator $\hf_\la$ have been studied in \cite{marteau2019} in the case of generalized self concordance, an are reported in the main lines in \cref{app:proj_stat}. The statistical rates of this estimator and the optimal choice of $\la$ is determined by two parameters, defined in \cref{df:main_quant_stat} and which we adapt to the Kernel problem here.

\begin{itemize}
    \item the \emph{bias} $\bias_\la = \|\Hess_\la(\fstar)^{-1/2} \nabla L_\la(\fstar)\| =\la \|\fstar\|_{\Hess^{-1}_\la(\fstar)}$, which characterizes the regularity of the optimum. The faster $\bias_\la$ decreases to zero, the more regular $\fstar$ is.
    \item the \emph{effective dimension} 
    \eqal{\label{eq:deff-statistics}
    \Deff_\la = \Exp{\|\Hess_\la(\fstar)^{-1/2}\nabla L_Z(\fstar)\|^2}.
    }
    This quantity characterizes the size of the space $\hh$ with respect to the problem; the slower it explodes as $\la$ goes to zero, the smaller the size of $\hh$. 
\end{itemize}

For more complete explanations on the meaning of these quantities, we refer to \cite{marteau2019}. \\\

Moreover, as mentioned in \cref{df:main_quant_stat}, one can define
\eqal{\label{eq:bb}\Bos:= \sup_{z \in \supp(Z)}{\|\nabla L_z(\fstar)\|},~ \Bts:= \sup_{z \in \supp(Z)}{\Tr(\nabla^2 L_z(\fstar))},~ \Qs = \frac{\Bos}{\sqrt{\Bts}}, ~ \bts = \btns(\fstar).}

We assume the following regularity condition on the minimizer $\fstar$, in order to get statistical bounds.  

\begin{restatable}[Source condition]{ass}{asmsource1}
\label{asm:source}
There exists $r >0 $ and $g \in \hh$ such that $\fstar = \Hess^r(\fstar)g$. This implies the following decrease rate of the bias:
\[\bias_\la \leq \Lc \la^{1/2+r},\qquad \Lc = \|g\|_{\hh}.\]
\end{restatable}

This is a stronger assumption than the existence of the minimizer as $r > 0$ is crucial for our analysis. 

We also quantify the effective dimension $\Deff_\la$: (however, since it always holds for $\alpha  =1$, this is not, strictly speaking, an additional assumption). 

\begin{restatable}[Effective dimension]{ass}{asmeffdim}
\label{asm:deff}
The effective dimension decreases as $\Deff_\la \leq \Qc \la^{-1/\alpha}$. 
\end{restatable}

If these two assumptions hold, define: 

\[\beta = \frac{\alpha}{1 + \alpha(1+2r)},\qquad \gamma = \frac{(1+2r)\alpha}{1 + \alpha(1+2r)}.\]

Under these assumptions, one can obtain the following statistical rates (which can be found in \cite{marteau2019} or in \cref{cor:bounds_source_capacity}).

 \bp\label{prp:bounds_source_capacity_K}
Let $\delta \in (0,1/2]$. Under \cref{asm:iid1,asm:deff,asm:asmtech1,asm:minimizer1,asm:source}, when $n \geq N$ and $\la = (C_0/n)^{\beta}$, then with probability at least $1-2\delta$,
\[L(\hf_\la) - L(\fstar) \leq C_1 n^{-\gamma}\log \frac{2}{\delta},\]
with $C_0 = 256 (\Qc/\Lc)^2,~C_1 = 8 (256)^{\gamma}~(\Qc^{\gamma} ~ 
    \Lc^{1-\gamma})^2$ and $N$ defined in \cite{marteau2019}, and satisfying
$N = O(\poly(\Bos,\Bts,\Lc,\Qc,\Rad,\log(1/\delta)))$. 
\ep

\subsection{Reducing the dimension: projecting on a subspace using Nystr\"om sub-sampling.\label{app:kern-rand-proj}}

\paragraph{Computations}

Using a representer theorem, one of the key properties of Kernel spaces is that, owing to the reproducing property, 
$$\hf_\la \in \hh_n:= \left\{ \sum_{i=1}^n{\alpha_i K_{x_i}}~:~(\alpha_i) \in \R^n\right\}.$$
This means that solving the regularized empirical problem can be turned into a finite dimensional problem in $\alpha$. Indeed
 $\hf_\la = \sum_{i=1}^n{\alpha_i K_{x_i}}$ where $\alpha =  (\alpha_i)_{1\leq i\leq n}$ is the solution to the following problem:
\[\alpha = \argmin_{\alpha \in \R^n}{\frac{1}{n}\sum_{i=1}^n{\ell_{z_i}(\alpha^\top {\bf K}_{nn}e_i)} + \frac{\la}{2}\alpha^\top {\bf K}_{nn} \alpha},\qquad {\bf K}_{nn} = (K(x_i,x_j))_{1 \leq i,j\leq n}\in\R^{n \times n}.\]

The previous problem is usually too costly to solve directly for large values of $n$, both in time and memory, because of the operations involving ${\bf K}_{nn}$.  A solution consists in looking for a solution in a smaller dimensional sub-space $\hh_M$ constructed from sub-samples of the data $\left\{ \tilde{x}_1,...,\tilde{x}_M\right\} \subset \left\{x_1,...,x_n\right\}$:
\[\hh_M:= \left\{\sum_{j=1}^M{\tilde{\alpha}_j K_{\tilde{x}_j}}~:~ \tilde{\alpha} \in \R^M\right\}.\]

In this case, the minimizer $\hf_{M,\la}  = \argmin_{f \in \hh_M}{\hL_\la(f)}$ can be written $\hf_{M,\la} = \sum_{j=1}^M{\tilde{\alpha}_j K_{\tilde{x}_j}}$, where $\tilde{\alpha}$ is the solution to the following problem:
\[\tilde{\alpha} = \argmin_{\alpha \in \R^M}{\frac{1}{n}\sum_{i=1}^n{\ell_{z_i}(\alpha^\top {\bf K}_{Mn}e_i)} + \frac{\la}{2}\alpha^\top {\bf K}_{MM} \alpha}, \]
where 
\[{\bf K}_{nM} = (K(x_i,\tilde{x}_j))_{\substack{1 \leq i \leq n \\ 1 \leq j \leq M}},~{\bf K}_{Mn} = {\bf K}_{nM}^\top,~ {\bf K}_{MM}:= (K(\tilde{x}_i,\tilde{x}_j))_{1 \leq i,j \leq M}.\]

Let $\bf T$ be an upper triangular matrix such that ${\bf T}^\top {\bf T} = {\bf K}_{MM}$. One can re-parametrize the previous problem in the following way. For any $\beta \in \R^M$, define $f_\beta = \sum_{j=1}^M{[{\bf T}^{\dagger}\beta]_j ~K_{\tilde{x}_j}}$. This implies in particular that $\|f_\beta\|_{\hh} = \|\beta\|_{\R^M}$. Then $\hf_{M,\la} = f_{\beta_{M,\la}}$, where 

\[\beta_{M,\la} = \argmin_{\beta \in \R^M}{\hL_{M,\la}(\beta):= \frac{1}{n}\sum_{i=1}^n{\ell_{z_i}(e_i^\top {\bf K}_{nM}{\bf T}^{\dagger}\beta)} + \frac{\la}{2}\|\beta\|^2}.\]

Using the properties the $\ell_z$, one easily shows that 
$ \beta \mapsto \ell_{z_i}(e_i^\top {\bf K}_{nM}{\bf T}^{\dagger}\beta)$ is $\left\{\Rell {\bf T}^{-\top}{\bf K}_{Mn}e_i\right\}$ generalized self-concordant, and $\|\Rell {\bf T}^{-\top}{\bf K}_{Mn}e_i\| \leq \Rell \sqrt{K(x_i,x_i)}$. Thus, $\hL_M$ is also generalized self-concordant, and the associated $\hG_M$ is bounded by $\Rad= \Rell \ka$. It will therefore be possible to apply the second order scheme presented in this paper to approximately compute $\beta_{M,\la}$.

\paragraph{Statistics}

Let $\hnd_{\la,M}(\beta)$ denote the Newton decrement of $\hL_{\la,M}$ at point $\beta$ and $\Pj_M$ denote the orthogonal projection on $\hh_M$. Then the following statistical result shows that provided $\beta$ is a good enough approximation of the optimum, and provided $\hh_M$ is large enough, then $f_\beta$ has the same generalization error as the empirical risk minimizer $\hf_\la$. 

Recall the following result proved in \cref{prp:proj_pb_stat_gen_part} in \cref{app:bounds_proj}. 
\bp[Behavior of an approximation to the projected problem] \label{prp:proj_pb_stat_gen_kernels} Suppose that 
\cref{asm:iid1,asm:asmtech1,asm:minimizer1} are satisfied.
Let $n \in \N$, $\delta \in (0,1/2]$, $0 < \la \leq \Bts$. Whenever
\[n \geq \triangle_1 \frac{\Bts}{\la}\log \frac{8 \square_1^2\Bts}{\la \delta},\qquad \Cone \sqrt{ \frac{ \Deff_\lambda \vee (\Qs)^2}{n} ~\log \frac{2}{\delta}}\leq \frac{\la^{1/2}}{\Rad}, \qquad \Cone \bias_\la \leq \frac{\la^{1/2}}{\Rad},\]
if
\[\|\Hess^{1/2}(\fstar)(\Id - \Pj_M)\|^2 \leq \la \frac{\sqrt{2}}{480},~ 126 \hnd_{M,\la}(\beta) \leq \frac{\la^{1/2}}{\Rad}, \]

the following holds, with probability at least $1 - 2 \delta$.

\[L(f_{\beta}) - L(\fstar) \leq \Cfone~ \bias_\la^2 + \Cftwo~ \frac{\Deff_\la \vee (\Qs)^2}{n}~\log \frac{2}{\delta} + \Cfthree~ \hnd^2_{M,\la}(\beta), \qquad \Rad \|f_{\beta} - \fstar\|_{\hh} \leq 10, \]

where $\Cfone \leq 6.0\mathrm{e}{4}$, $\Cftwo \leq 6.0\mathrm{e}{6} $ and $\Cfthree \leq 810$, $\Cone$ is defined in \cref{lm:technical_before_proj}, and the other constants are defined in \cref{thm:general-result}. 
\ep 

In particular, if we apply the previous result for a fixed $\la$, the following theorem holds (for a proof, see \cref{app:quantibounds}).
\bt[Quantitative result with source $r > 0$]\label[theorem]{thm:quantisource_kernels} Suppose that 
\cref{asm:iid1,asm:deff,asm:asmtech1,asm:minimizer1,asm:source} are satisfied.
Let $n \geq N$ and $\delta \in (0,\frac{1}{2}]$. If $\lambda = \left( \left(\frac{\Qc}{\Lc}\right)^{2}~\frac{1}{n}\right)^{\frac{\alpha}{\alpha(1+2r) + 1}}$, and if 
\[\|\Hess^{1/2}(\fstar)(\Id - \Pj_M)\|^2 \leq \la \frac{\sqrt{2}}{480},~\hnd_{M,\la}(\beta) \leq \Qc^{\gamma} ~ 
    \Lc^{1-\gamma} n^{-\gamma/2},\]
    then with probability at least $1-2 \delta$,
    \[L(f_\beta) - L(\fstar) \leq \Kfinal \left(\Qc^{\gamma} ~ 
    \Lc^{1-\gamma}\right)^2 ~~\frac{1}{n^{\gamma}} \log \frac{2}{\delta}, \qquad \Rad\|f_\beta - \fstar\| \leq 10,\]
    where $N$ is defined in \cref{eq:Nvalue} and $\Kfinal \leq 7.0e6$.
\et

The proof of the previous result is quite technical and can be found in \cref{app:proj_stat}, in \cref{thm:quantisource}.

\subsection{A note on sub-sampling techniques\label{app:kern-ny-examples}}
Let $Z$ be a random variable on a Polish space $\Z$ and $(v_z)_{z \in \Z}$ be a family of vectors in $\hh$ such that $||v||_{L^\infty(Z)}:= \sup_{z \in \supp(Z)}\|v_z\| < \infty $ is bounded. Assume that $z_1,...,z_n$ are i.i.d. samples from $Z$.

Define the following trace class Hermitian operators: 

$$\Ab = \Exp{v_Z \otimes v_Z},~\hAb = \frac{1}{n}\sum_{i=1}^n{v_{z_i} \otimes v_{z_i}}.$$ 

Define 
\eqal{\label{eq:deg_freed}\Ny^{\Ab}(\la):= \Tr(\Ab_\la^{-1}\Ab),\qquad \Nyi^{\Ab}(\la):= \sup_{z \in \supp(Z)}\|\Ab_\la^{-1/2}v_z\|^2 .}
We typically have:
\[\Ny^{\Ab}(\la) \leq \Nyi^{\Ab}(\la) \leq \frac{\|v\|_{L^\infty(Z)}^2}{\la}.\]

We define the leverage scores associated to the points $z_i$ and $\Ab$: 
\eqal{\label{eq:lev_scores} \forall 1 \leq i \leq n,~\forall t > 0,~ l^{\Ab}_i(t) = \|\hAb^{-1/2}_{t}v_{z_i}\|^2 = n\left((\Gb_{nn} + t n \Id)^{-1}\Gb_{nn}\right)_{ii},}
where $\Gb_{nn} = (v_{z_i}\cdot v_{z_j})_{1 \leq i,j \leq n}$ denotes the Gram matrix associated to the family $v_{z_i}$. 

As in \cite{Rudi15}, definition 1, we give the following definition for leverage scores.
\begin{definition}[$q$-approximate leverage scores]
given $t_0$, a family $(\tilde{l}^{\Ab}_i(t))_{1\leq i \leq n}$ is said to be a family of $q$-approximate leverage scores with respect to $\Ab$ if  
\[\forall 1 \leq i \leq n,~ \forall t \geq t_0,~ \frac{1}{q}~ l^{\Ab}_i(t) \leq \tilde{l}^{\Ab}_i(t) \leq q~l^{\Ab}_i(t).\]
\end{definition}

We say that a subset of $m$ points $\left\{\tilde{z}_1,...,\tilde{z}_m\right\} \subset \left\{z_i ~:~1 \leq i\leq n\right\}$ is:
\begin{itemize}
    \item \textbf{Sampled using $q$-approximate leverage scores for $t$} if the $\tilde{z}_j  = z_{i_j}$ where the $i_j$ are $m$ i.i.d. samples from $\left\{1,...,n\right\}$ using the probability vector $p_i = \frac{\tilde{l}^{\Ab}_i(t)}{\sum_{\tilde{i}=1}^n{\tilde{l}^{\Ab}_{\tilde{i}}(t)}}$.
    In that case, we define $\hAb_{m}:= \frac{1}{m}\sum_{j=1}^m{\frac{1}{np_{i_j}} v_{\tilde{z}_j} \otimes v_{\tilde{z}_j}}$.
    \item \textbf{Sampled uniformly} if the $\left\{i_j~:~1 \leq j \leq m\right\}$ is a uniformly chosen subset of  $\left\{1,...,n\right\}$ of size $m$. 
    In this case, we define 
    $\hAb_{m}:= \frac{1}{m}\sum_{j=1}^m{ v_{\tilde{z}_j} \otimes v_{\tilde{z}_j}}$.
\end{itemize}

In \cref{app:sub_nystrom}, we present technical lemmas which allow us to show that if $m$ is large enough, the following hold:
\begin{itemize}
    \item $\|\Ab_\eta(\Id - \Pj_m)\|^2 \leq 3 \eta$, where $\Pj_m$ is the orthogonal projection on the subspace induced by the $v_{\tilde{z}_j}$;
    \item $\hAb_{m,\la}$ is equivalent to $\hAb_\la$.
\end{itemize}

\br[cost of computing $q$-approximate leverage scores]\label{rmk:cost_lev_scores}
In \cite{rudi2018fast}, one can show that the complexity of computing $q$-approximate leverage scores can be achieved in: 
$\csla = O(q^2 \Ny^{\Ab}(\la)^2 \min(n,1/\la))$ time (where a unit of time is a scalar product evaluation) and $O(\Ny^{\Ab}(\la)^2 + n)$ in memory.
\er

\subsection{Selecting the $M$ Nystr\"om points\label{app:selecting_nys}}

In order for \cref{thm:quantisource_kernels} to hold, we must subsample the $M$ points such as to guarantee $\|\Hess^{1/2}(\fstar)(\Id - \Pj_M)\|^2  \leq \frac{\sqrt{2}\la}{480}$.

Since we must sub-sample the $M$ points a priori, i.e. before performing the method, it is necessary to have sub-sampling schemes which do not depend heavily on the point. Define the covariance operator:
$$\Cov = \Exp{K_X \otimes K_X}.$$
Since $\Hess(\fstar) = \Exp{\ell^{\prime \prime}_Z(f(X)) ~K_X \otimes K_X }$, it is easy to see that $\Hess(\fstar)  \preceq \bts \Cov$. Note that for $\Cov$, since $\widehat{\Cov} = \frac{1}{n}\sum_{i=1}^n{K_{x_i} \otimes K_{x_i}}$, the leverage scores have the following form:
\[\forall 1 \leq i \leq n,~ l_i^{\Cov}(t) = n\left(({\bf K}_{nn} + \la n \Id)^{-1} {\bf K}_{nn}\right)_{ii}.\]

\bp[Selecting Nystr\"om points]\label{prp:nystrom_select}
Let $\delta >0$. Let $\eta = \min(\|\Cov\|,\frac{\la \sqrt{2}}{1440 (\bts\vee 1)})$. Assume the samples  $\left\{\tilde{x}_1,...,\tilde{x}_M\right\}$ are obtained with one of the following.\\
\noindent{\bf 1.} $n \geq  M \geq \left(10 + 160 \Nyi^\Cov(\eta)\right)\log \frac{8 \ka^2}{\eta \delta}$ using uniform sampling;\\
\noindent{\bf 2.} $  M \geq \left(6 + 486 q^2 \Ny^\Cov(\eta)\right)\log \frac{8 \ka^2}{\eta \delta}$ using $q$-approximate leverage scores with respect to $\Cov$ for $t = \eta$, $t_0 \vee \frac{19 \ka^2}{n}\log \frac{n}{2 \delta} < \eta$, $n \geq 405 \ka^2 \vee 67 \ka^2 \log \frac{12 \ka^2}{\delta}$.\\
Then it holds, with probability at least $1-\delta$:
\[\|\Cov_{\eta}^{1/2}(\Id - \Pj_M)\| \leq 3 \eta \implies \|\Hess^{1/2}(\fstar)(\Id - \Pj_M)\|^2 \leq \la \frac{\sqrt{2}}{480}.\]

\ep 

\bpr 
The proof is a direct application of the lemmas in \cref{app:sub_nystrom}. Indeed, note that since $\Cov = \Exp{K_X \otimes K_X}$, then the results can be applied with $Z \leftarrow X$ and $v_z \leftarrow K_x$. Indeed, from \cref{asm:asmtech1}, it holds:
\[\sup_{x\in \supp(X)}\|K_x\|^2 \leq \ka^2.\]
\epr 

We can now combine \cref{prp:nystrom_select} and \cref{prp:proj_pb_stat_gen_kernels} to obtain the following statistical bounds for the optimizer of the projected Nystr\"om problem $\beta_{M,\la}$.

\bt \label{thm:nystrom_estimator_bounds}
Suppose that 
\cref{asm:iid1,asm:asmtech1,asm:minimizer1} are satisfied.
Let $n \in \N$, $\delta \in (0,1/2]$, $0 < \la \leq \Bts \wedge 720\sqrt{2} (\bts\vee 1)\|\Cov\| $. Assume
\[n \geq \triangle_1 \frac{\Bts}{\la}\log \frac{8 \square_1^2\Bts}{\la \delta},\qquad \Cone \sqrt{  \frac{\Deff_\lambda \vee (\Qs)^2}{n} ~\log \frac{2}{\delta}} \leq \frac{\la^{1/2}}{\Rad},\qquad  \Cone \bias_\la \leq \frac{\la^{1/2}}{\Rad},\]
Let $\eta = \frac{\la \sqrt{2}}{1440 (\bts\vee 1)}$. Assume the samples  $\left\{\tilde{x}_1,...,\tilde{x}_M\right\}$ are obtained with one of the following.\\
\noindent{\bf 1.} $n \geq  M \geq \left(10 + 160 \Nyi^\Cov(\eta)\right)\log \frac{8 \ka^2}{\eta \delta}$ using uniform sampling;\\
\noindent{\bf 2.} $  M \geq \left(6 + 486 q^2 \Ny^\Cov(\eta)\right)\log \frac{8 \ka^2}{\eta \delta}$ using $q$-approximate leverage scores with respect to $\Cov$ for $t = \eta$, $t_0 \vee \frac{19 \ka^2}{n}\log \frac{n}{2 \delta} < \eta$, $n \geq 405 \ka^2 \vee 67 \ka^2 \log \frac{12 \ka^2}{\delta}$.\\
The following holds, with probability at least $1 - 3 \delta$.
\[L(f_{\beta_{M,\la}}) - L(\fstar) \leq \Cfone~ \bias_\la^2 + \Cftwo~ \frac{\Deff_\la \vee (\Qs)^2}{n}~\log \frac{2}{\delta} ,\qquad \Rad \|\beta_{M,\la}\| \leq \Rad \|\fstar\|  + 10,\]
where $\Cfone \leq 6.0\mathrm{e}{4}$, $\Cftwo \leq 6.0\mathrm{e}{6} $ and $\Cfthree \leq 810$, $\Cone$ is defined in \cref{lm:technical_before_proj}, and the other constants are defined in \cref{thm:general-result}. 

\et

\bpr 
This is simply a reformulation of \cref{prp:proj_pb_stat_gen_kernels}, noting that $\hnd_{M,\la}(\beta_{M,\la}) = 0$ and that \cref{prp:nystrom_select} implies the condition on the Hessian at the optimum.
\epr

Provided source condition holds with $r>0$, the conditions of this theorem are not void.

\subsection{\label{app:approximate_b}Performing the globalization scheme to approximate $\beta_{M,\la}$}

In order to apply \cref{prp:proj_pb_stat_gen_kernels}, one needs to control $\hnd_{M,\la}(\beta)$.

We will apply our general scheme to $\hL_{M,\la}$ in order to obtain such a control. 

\subsubsection{Performing approximate Newton steps \label{app:ans_kernels}}

The key element in the globalization scheme is to be able to compute $\frac{1}{7}$-approximate Newton steps.

 Note that at a given point $\beta$ and for a given $\mu > 0$ the Hessian is of the form: 

\[\hHess_{M,\mu}(\beta) = \frac{1}{n} {\bf T}^{-\top}{\bf K}_{Mn} {\bf D}_n(\beta) {\bf K}_{nM}{\bf T}^{-1} + \mu \Id_M ,\]
where ${\bf D}_n(\beta) =  \diag((d_i(\beta))_{1 \leq i\leq n})$ is a diagonal matrix whose elements are given by $d_i(\beta) = \ell_{z_i}^{\prime\prime}(e_i^\top {\bf K}_{nM} {\bf T}^{-1} \beta)$.

Note that we can always write
\[\hHess_{M,\mu}(\beta) = \frac{1}{n} \sum_{i=1}^n{u_i(\beta) u_i(\beta)^\top}+ \mu \Id,\qquad u_i(\beta) = \sqrt{d_i(\beta)}{\bf T}^{-\top}{\bf K}_{Mn} e_i\]

The gradient can be put in the following form:

\[\nabla \hL_{M,\mu}(\beta) = \frac{1}{n}{\bf T}^{-\top}{\bf K}_{Mn}  v + \mu \beta,\qquad v = (\ell_{z_i}^{\prime}(e_i^\top {\bf K}_{nM}{\bf T}^{-1}\beta))_{1 \leq i \leq n}.\]

Computing the gradient at one point therefore costs $O(nM + M^2)$, this being the cost of computing ${\bf K}_{nM}$ times a vector costs $O(nM)$ and computing ${\bf T}^{-1}$ times a vector takes $O(M^2)$ since ${\bf T}$ is triangular. Moreover, the cost in memory is $O(M^2 + n)$, $M^2$ being needed for the saving of ${\bf T}$ and $n$ for the saving of the gradient; ${\bf K}_{nM}$ times a vector can also be done in $O(n)$ memory, provided we compute it by blocks. \\\

On the other hand, computing the full Hessian matrix would cost $nM^2$ operations, which is un-tractable. However, computing a Hessian vector product can be done in $O(nM + M^2)$ time, as for the gradient, which suggest using an iterative solver with preconditioning.

\paragraph{Computing $x \in \lso(\Ab,b,\err)$ through pre-conditioned conjugate gradient descent.}

Assume we wish to solve the problem 
$\Ab x = b$ where $\Ab \in \R^{M \times M}$ is a positive definite matrix and $b$ is a vector of $\R^M$. If one uses the conjugate gradient method starting from zero, then if $x_k$ denotes the $k$-the iterate of the conjugate gradient algorithm, Theorem 6.6 in \cite{Saad03} shows that  
\[x_k \in \lso(\Ab,b,\err),\qquad \err = 2\left(\frac{\sqrt{\cond(\Ab)} -1 }{\sqrt{\cond(\Ab)} + 1}\right)^k.\]

where $\cond(\Ab)$ is the condition number of the matrix $\Ab$, namely the ratio $\frac{\la_{\max}(\Ab)}{\la_{\min}(\Ab)}$. If $\cond(\Ab)$ is large, this convergence can be very slow. The idea of preconditioning is to compute an approximation matrix $\tAb$ such that 
\eqal{\label{eq:tab}\frac{1}{2}\tAb \preceq \Ab \preceq \frac{3}{2} \tAb.}

We then compute ${\bf B}$ a triangular matrix such that ${\bf B}^\top {\bf B} = \tAb$ using a cholesky decomposition, which can be done in $O(M^3)$, and note that ${\bf B}^{-\top} \Ab {\bf B}^{-1}$ is very well conditioned; indeed, its condition number is bounded by $3$. 

Perform a conjugate gradient method to solve the pre-conditioned problem ${\bf B}^{-\top} \Ab {\bf B}^{-1} z = {\bf B}^{-\top} b$, and denote with $z_\tau$ the $\tau$-th iteration of this method. Then 
using the bound on the condition number, we find 

\[z_\tau \in \lso({\bf B}^{-\top}\Ab {\bf B}^{-1},{\bf B}^{-\top}b,\rho),\qquad \rho = 2 \left(\frac{\sqrt{3}-1}{\sqrt{3}+1}\right)^{\tau},\]
which in turn implies that by setting $x_\tau:= {\bf B}^{-1}z_\tau$, 
\[x_\tau \in \lso(\Ab,b,\rho),~ \rho = 2 \left(\frac{\sqrt{3}-1}{\sqrt{3}+1}\right)^{\tau}.\]

This shows that after at most $\tau = 3$  iterations, provided $\tAb$ satisfies \cref{eq:tab}, $x_\tau \in \lso(\Ab,b,\frac{1}{7})$. The cost of this method is therefore $ O(M^3 + nM) $ in time, and $O(n + M^2)$ due to the computing of the preconditioner and computing matrix vector products by block. This does not include the cost of finding a suitable $\tAb$.

\paragraph{Computing a suitable approximation of $\hHess_{M,\mu}(\beta)$}

To compute a good pre-conditioner, we will subsample $Q$ points $i_1,...,i_Q$ points from $\left\{1,...,n\right\}$,  and sketch the Hessian using these $Q$ points.

\bp[Computing approximate newton steps]\label{prp:approx_ns_kernels} Let $\delta >0$. Let $\beta \in \R^M$ and $\mu \geq \la $, and assume $\frac{19 \btns(f_\beta)\ka^2}{n}\log \frac{n}{2 \delta} < \la$ and $n \geq 405 \btns(f_\beta)\ka^2 \vee 67 \btns(f_\beta)\ka^2 \log \frac{12 \btns(f_\beta)\ka^2}{\delta}$. Let $\tilde{\mu} = \min(\mu,\|\Hess(f_\beta)\|)$. Assume one of the following properties is satisfied

{\noindent \bf 1.} $Q \geq \left(10 + 160 \Nyi^{\Hess(f_\beta)}(\tilde{\mu})\right)\log \frac{8 \btns(f_\beta)\ka^2}{\tilde{\mu} \delta}$ with uniform sampling of the  $\left\{i_1,...,i_Q\right\}$. We set ${\bf D}_Q = \diag(\ell^{\prime \prime}_{z_{i_j}}(f_\beta(x_{i_j})))_{1 \leq j \leq Q}$\\
{\noindent \bf 2.} $Q \geq \left(6 + 486 q^2 \Ny^{\Hess(f_\beta)}(\tilde{\mu})\right)\log \frac{8 \btns(f_\beta)\ka^2}{\tilde{\mu} \delta}$ using $q$-approximate leverage scores associated to $\Hess(f_{\beta})$ for $t = \tilde{\mu}$.
We set ${\bf D}_Q = \diag\left(\frac{\ell^{\prime \prime}_{z_{i_j}}(f_\beta(x_{i_j}))}{p_{i_j}}\right),$
where the $p_{i_j}$ are the probabilities computed from the leverage scores.\\

Assume we use a pre-conditioner ${\bf B}$ such that 
$${\bf B}^\top {\bf B} =  \frac{1}{Q}{\bf T}^{-\top} {\bf K}_{MQ} {\bf D}_Q {\bf K}_{QM} {\bf T}^{-1} + \mu \Id_M,\qquad {\bf K}_{QM} = (K(x_{i_j},\tilde{x}_k))_{\substack{1 \leq j \leq Q \\ 1 \leq k \leq M}}.$$

If we perform $\tau = \log(\rho/2) / \log((\sqrt{3} +1)/{\sqrt{3} -1})$ iterations of the conjugate gradient descent on the pre-conditioned Newton system using ${\bf B}$ as a preconditioner, then with probability at least $1-\delta$, this procedure is returns $\nsa \in \lso(\hHess_{M,\la}(\beta),\nabla \hL_{M,\la}(\beta),\rho)$, and the computational time is of order $O(\tau(Mn + M^2 Q + M^3 + \csla))$, and the memory requirements can be reduced to $O(M^2 + n)$. Here $\csla$ stands for the complexity of computing Nystrom leverage scores, and using \cref{rmk:cost_lev_scores} or \cite{rudi2018fast}, $\csla = O(1)$ if uniform sampling is used, and $\csla = O(\Ny^{\Hess(f_\beta)}(\tilde{\mu})^2/\la)$ if Nystrom sub-sampling is used. Note that for $\tau = 3$, $\rho = \frac{1}{7}$.

\ep

\bpr  Start by defining the following operators:

\begin{itemize}
    \item  $K_n: f\in \hh \rightarrow  (f(x_i))_{1\leq i\leq n} \in \R^n$;
    \item  $K_M: f \in \hh \rightarrow  (f(\tilde{x}_j))_{1\leq j\leq M} \in \R^M$;
    \item $V = K_M^* {\bf T}^{-1}$, where ${\bf T}$ is an upper triangular matrix such that ${\bf T}^\top {\bf T} = {\bf K}_{MM} = K_M K_M^*$. 
\end{itemize}

Note that $K_n V = {\bf K}_{nM} {\bf T}^{-1}$. 

Now note that 
\[\forall f \in \hh,~ \Hess(f) = \Exp{v_z \otimes v_z},\qquad \hHess(f) = \frac{1}{n} \sum_{i=1}^n{v_{z_i} \otimes v_{z_i}} ,\qquad v_z = \sqrt{\ell^{\prime \prime}_z(f(x))}K_x. \]

Since for any $f \in \hh$, $\hHess(f) = \frac{1}{n} K_n^* {\bf D}_n(f) K_n $, where ${\bf D}_n(f) = \diag(\ell^{\prime \prime}_{z_i}(f(x_i)))$, we see that 

\[\hHess_{M,\mu}(\beta) = V^* \hHess(f_{\beta}) V + \mu \Id_M. \]

Thus, the last lemma of \cref{app:sub_nystrom} can be applied, using the fact that $\|v_z\|^2 \leq \btns(f)\ka^2$, to get that in both cases of the proposition, under the corresponding assumptions: 

\[ \frac{1}{2}\left(\frac{1}{Q}{\bf T}^{-\top} {\bf K}_{MQ} {\bf D}_Q {\bf K}_{QM} {\bf T}^{-1} + \mu \Id_M\right) \preceq \hHess_{M,\mu}(\beta) \preceq \frac{3}{2}\left(\frac{1}{Q}{\bf T}^{-\top} {\bf K}_{MQ} {\bf D}_Q {\bf K}_{QM} {\bf T}^{-1} + \mu \Id_M\right).\]

The rest of the proposition follows from the previous discussion.

\epr

\subsubsection{Applying the globalization scheme to control $\hnd_{M,\la}(\beta)$ \label{sssec:gen_scheme_optim} }

In order to apply 
\cref{prp:approx_ns_kernels} to each point $\beta$ in our method, we need to have a globalized version of the condition of this proposition.

First, we start by localizing the different values of $\beta$ we will visit throughout the algorithm.

\bd[path of regularized solutions]
Let $\la > 0$, $\varepsilon >0$. Define the path of regularized solutions
\eqal{
\label{df:path_M}
\pthm_{\la}:= \left\{\beta_{M,\mu}~:~ \mu \geq \la\right\}
.}
And the $\varepsilon$ approximation of this path: 
\eqal{
\label{df:path_M_eps}
\pthm_{\la,\varepsilon}:= \left\{\beta \in \R^M~:~ d(\beta,\pthm_{\la}) \leq \varepsilon \right\}.
}

\ed 

Note that we always have $\pthm_\la \subset \mathcal{B}_{\R^M}(\|\beta_{M,\la}\|)$. 
We now state a lemma proving that all the values visited during the algorithm will lie in an approximation of this path.

\blm Define
Let $\beta \in \R^M$ such that $\hnd_{M,\mu}(\beta) \leq \frac{\mu^{1/2}}{7\Rad}$ for some $\mu \geq \la$. Then the following holds:
\[\beta \in \pthm_{\la,\frac{1}{6\Rad}}.\]
\elm

\bpr
Bound 
\[\Rad \|\beta - \beta_{M,\mu}\| \leq \frac{\Rad}{\mu^{1/2}} \|\beta - \beta_{M,\mu}\|_{\hHess_{M,\mu}(\beta)} \leq \frac{1}{\dikin(\tn_{M}(\beta - \beta_{M,\mu}))} \frac{\Rad \hnd_{M,\mu}(\beta)}{\mu^{1/2}}.\]
Just apply \cref{eq:gl} to obtain $\Rad \|\beta - \beta_{M,\mu}\| \leq \frac{1}{6}$. 
\epr

We now introduce the following quantities which will allow to control the number of sub-samples throughout the whole algorithm.

\bd Define 
\begin{itemize}
  \item   
$\btb:= \sup_{\beta \in \pthm_{\la,1/6\Rad}}{\btns(f_\beta)}$.
\item $\Nyb^{\Hess}(\la) = \sup_{\beta \in \pthm_{\la,1/6\Rad}}\Ny^{\Hess(f_\beta)}(\la)$.
\item $\Nyib^{\Hess}(\la) = \sup_{\beta \in \pthm_{\la,1/6\Rad}}\Nyi^{\Hess(f_\beta)}(\la)$.
\item $\overline{\|\Hess\|} = \min_{\beta \in \pthm_{\la,1/6\Rad}}\|\Hess(f_\beta)\|$.
\end{itemize}

\ed

\bp[Performance of the globalization scheme]\label{prp:gen_scheme_optim}
Let $\varepsilon > 0$, $\delta >0$, $\tilde{\la} = \min(\la,\overline{\|\Hess\|})$. Assume $\frac{19 \btb \ka^2}{n}\log \frac{n}{2 \delta} < \tilde{\la}$ and $n \geq 405 \btb \ka^2 \vee 67 \btb \ka^2 \log \frac{12 \btb \ka^2}{\delta}$.

Assume we perform the globalization scheme with the parameters in \cref{thm:easy_first_phase},  where in order to compute any $\rho$ approximation of a regularized Newton step, we use a conjugate gradient descent on the pre-conditioned system, where the pre-conditioner is computed as in \cref{prp:approx_ns_kernels} using \\
\noindent{\bf 1.} $Q \geq \left(10 + 160 \Nyib^{\Hess}(\tilde{\la})\right)\log \frac{8 \btb\ka^2}{\tilde{\la} \delta}$ if using uniform sampling\\
\noindent{\bf 2.} $Q \geq \left(6 + 486 q^2 \Nyb^{\Hess}(\tilde{\la})\right)\log \frac{8 \btb\ka^2}{\tilde{\la} \delta}$ if using Nystr\"om leverage scores\\

Recall that $t$ denotes the number of approximate Newton steps performed at for each $\mu$ in Phase I and $T$ denotes the number of approximate Newton steps performed in Phase II, and that using \cref{thm:easy_first_phase}, $t = 2$ and $T = \lceil \log_2 \sqrt{1 \vee (\la \varepsilon^{-1}/\Rad^2)}   \rceil$. Moreover, recall that $K$ denotes the number of steps performed in Phase I. Define 
\[N_{ns}:= 2  \left\lfloor \left(3 + 11 \Rad \|\beta_{M,\la}\|\right)\log_2 (7 \Rad \|\nabla \hL_M(0)\|/\la)\right\rfloor +  \lceil \log_2 \sqrt{1 \vee (\la \varepsilon^{-1}/\Rad^2)}   \rceil.\]

Then with probability at least $(1-\delta)^{N_{ns}}$: 

\begin{itemize}
    \item The method presented in \cref{prp:approx_ns_kernels} returns a $1/7$- approximate Newton step at each time it is called in the algorithm.
    \item If $\beta$ denotes the result of the method, $\hnd_{M,\la}(\beta) \leq \sqrt{\varepsilon}$.
    \item The number of approximate Newton steps computed during the algorithm is bounded by $N_{ns}$; the complexity of the method is therefore of order $O(N_{ns}(M^2\max(M,Q)+ nM + \csla(\la)))$ in time and $O(M Q + M^2 +n)$ in memory, where $\csla(\la)$ is a bound on the complexity associated to the computing of leverage scores (see \cite{rudi2018fast} for details).
\end{itemize}

The algorithm is detailed in \cref{sec:algorithm}, in  \cref{alg:algo}. Note however that the notations are those of the main paper, which are slightly different from the ones used here.
\ep 
\bpr 
If we take the globalization scheme, using the parameters of \cref{thm:easy_first_phase}. Assume that all previous approximate Newton steps have been computed in a good way. Then the $\beta$ at which we are belongs to $\pthm_{\la,1/6\Rad}$. Thus, the hypotheses of this proposition imply that the hypothesis of \cref{prp:approx_ns_kernels} are satisfied; and hence, up to a $(1-\delta)$ probability factor, we can assume that the next approximate Newton step is performed correctly, continuing the globalization scheme in the right way. Thus, the globalization scheme converges as in \cref{thm:easy_first_phase}. 
\epr

\subsection{Statistical properties of the algorithm \label{app:main_kernels}}

The following theorem describes the computational and statistical behavior of our algorithm.

\bp[Behavior of an approximation to the projected problem] \label{prp:algo_stat_gen} Suppose that 
\cref{asm:iid1,asm:asmtech1,asm:minimizer1} are satisfied.\\
Let $n \in \N$, $\varepsilon > 0$, $\delta \in (0,1/2]$, $0 < \la \leq \Bts$. \\
Define $\tilde{\la} = \min(\la,\overline{\|\Hess\|})$ and assume $\frac{19 \btb \ka^2}{n}\log \frac{n}{2 \delta} < \tilde{\la}$, $n \geq 405 \btb \ka^2 \vee 67 \btb \ka^2 \log \frac{12 \btb \ka^2}{\delta}$, and $n \geq \triangle_1 \frac{\Bts}{\la}\log \frac{8 \square_1^2\Bts}{\la \delta}$. Assume
\[ \Cone \sqrt{ \frac{ \Deff_\lambda \vee (\Qs)^2}{n} ~\log \frac{2}{\delta}}\leq \frac{\la^{1/2}}{\Rad}, \qquad \Cone \bias_\la \leq \frac{\la^{1/2}}{\Rad}, \qquad 126 \sqrt{\varepsilon} \leq \frac{\la^{1/2}}{\Rad}.\]
Assume that the $M$ points $\tilde{x}_1,...,\tilde{x}_M$ are obtained through Nystr\"om sub-sampling using $\eta = \|\Cov\| \wedge \frac{\la \sqrt{2}}{1440 (\bts\vee 1)}$, with either\\
\noindent{\bf 1.} $M \geq \left(10 + 160 \Nyi^\Cov(\eta)\right)\log \frac{8 \ka^2}{\eta \delta}$ if using uniform sampling;\\
\noindent{ \bf 2.} 
$M \geq \left(6 + 486 q^2 \Ny^\Cov(\eta)\right)\log \frac{8 \ka^2}{\eta \delta}$ if using $q$-approximate leverage scores for $ \eta$, associated to the co-variance operator $\Cov$.

Assume we perform the globalization scheme as in \cref{prp:gen_scheme_optim}, i.e. with the parameters in \cref{thm:easy_first_phase},  where in order to compute any $\rho$ approximation of a regularized Newton step, we use a conjugate gradient descent on the pre-conditioned system, where the pre-conditioner is computed as in \cref{prp:approx_ns_kernels} using \\
\noindent{\bf 1.} $Q \geq \left(10 + 160 \Nyib^{\Hess}(\tilde{\la})\right)\log \frac{8 \btb\ka^2}{\tilde{\la} \delta}$ if using uniform sampling\\
\noindent{\bf 2.} $Q \geq \left(6 + 486 q^2 \Nyb^{\Hess}(\tilde{\la})\right)\log \frac{8 \btb\ka^2}{\tilde{\la} \delta}$ if using Nystr\"om leverage scores\\
Let $N_{ns}$ be defined as in \cref{prp:gen_scheme_optim}. Recall $N_{ns}$ is an upper bound for the number of approximate Newton steps performed in the algorithm. One can bound 
\[N_{ns} \leq 2  \left\lfloor \left(113 + 11 \Rad \|\fstar\|\right)\log_2 \frac{7 \Rad \|\nabla \hL_M(0)\|}{\la}\right\rfloor +  \left\lceil \log_2 \frac{\la^{1/2}}{\Rad \varepsilon}   \right\rceil. \]
Moreover, with probability at least $1 - (N_{ns} +2)\delta$, the following holds:

\[L(f_{\beta}) - L(\fstar) \leq \Cfone~ \bias_\la^2 + \Cftwo~ \frac{\Deff_\la \vee (\Qs)^2}{n}~\log \frac{2}{\delta} + \Cfthree~ \varepsilon .\]

where $\Cfone \leq 6.0\mathrm{e}{4}$, $\Cftwo \leq 6.0\mathrm{e}{6} $ and $\Cfthree \leq 810$, $\Cone$ is defined in \cref{lm:technical_before_proj}, and the other constants are defined in \cref{thm:general-result}. 
\ep 

\bpr 
This is a simple combination between \cref{prp:gen_scheme_optim,prp:nystrom_select,prp:proj_pb_stat_gen_kernels}. To bound the number of Newton steps $N_{ns}$, one simply uses the fact that under the conditions of the theorem, $\Rad\|\beta_{M,\la} \| \leq 10 + \Rad \|\fstar\|$. 
\epr

\br[Complexity] \label{rmk:complexity}
Let $L = \btb \ka^2$.  The complexity of the previous method using leverage scores computed for $\Cov$ for the Nystrom projections and for $\Hess(f_\beta)$ for choosing the $Q$ points at the different stages is the following. The total complexity in time will be of order:
\[O\left(N_{ns}\left(  n  \Ny^{\overline{\Hess}}(\la)\log(L \la^{-1}\delta^{-1})  + \btb^3 \Ny^{\Cov}(\la)^3 \log^3(L \la^{-1}\delta^{-1}) +  L / \la ~\btb^2 \Ny^{\Cov}(\la)^2 \right)\right).\]

The memory complexity can be bounded by
\[ O(\btb^2 \Ny^{\Cov}(\la)^2\log^2(L\la^{-1}\delta^{-1}) +n).\]
Here, we use the fact that $\Hess \leq \btb \Cov$.
\er 
We can now write down the previous proposition by classifying problems using \cref{asm:source,asm:deff} and in order to get optimal rates.

\bt[Performance of the scheme using pre-conditioning]\label{thm:main_kernels}
Let $\delta > 0$. 
Assume \cref{asm:iid1,asm:deff,asm:asmtech1,asm:minimizer1,asm:source} are satisfied. Let $n \geq \tilde{N}$, where $\tilde{N}$ is characterized in the proof,  $\lambda = \left( \left(\frac{\Qc}{\Lc}\right)^{2}~\frac{1}{n}\right)^{\frac{\alpha}{\alpha(1+2r) + 1}}$. \\\
Assume that the $M$ points $\tilde{x}_1,...,\tilde{x}_M$ are obtained through Nystr\"om sub-sampling using $\eta = \frac{\la \sqrt{2}}{1440 (\bts\vee 1)}$, with either\\
\noindent{\bf 1.} $M \geq \left(10 + 160 \Nyi^\Cov(\eta)\right)\log \frac{8 \ka^2}{\eta \delta}$ if using uniform sampling;\\
\noindent{ \bf 2.} 
$M \geq \left(6 + 486 q^2 \Ny^\Cov(\eta)\right)\log \frac{8 \ka^2}{\eta \delta}$ if using $q$-approximate leverage scores for $ \eta$, associated to the co-variance operator $\Cov$.

Assume we perform the globalization scheme as in \cref{prp:gen_scheme_optim}, i.e. with the parameters in \cref{thm:easy_first_phase},  where in order to compute any $\rho$ approximation of a regularized Newton step, we use a conjugate gradient descent on the pre-conditioned system, where the pre-conditioner is computed as in \cref{prp:approx_ns_kernels} using \\
\noindent{\bf 1.} $Q \geq \left(10 + 160 \Nyib^{\Hess}(\la)\right)\log \frac{8 \btb\ka^2}{\la\delta}$ if using uniform sampling\\
\noindent{\bf 2.} $Q \geq \left(6 + 486 q^2 \Nyb^{\Hess}(\la)\right)\log \frac{8 \btb\ka^2}{\la \delta}$ if using Nystr\"om leverage scores\\
Let $N_{ns}$ be defined as in \cref{prp:gen_scheme_optim}. Recall $N_{ns}$ is an upper bound for the number of approximate Newton steps performed in the algorithm. One can bound 
 \[ N_{ns} \leq  \left(227 + 22 \Rad \|\fstar\|\right)\left( \left\lceil \log_2 \left(7 \Rad \|\nabla \hL_M(0)\| \right)\right\rceil + \left\lceil \log_2\frac{n \Lc^2}{\Qc^2} \right\rceil +  \left\lceil \log_2 \frac{1}{\Rad \Lc} \right\rceil\right).\]
Moreover, with probability at least $1 - (N_{ns} +2)\delta$, the following holds:
\begin{itemize}
    \item all of the approximate Newton methods yield $\frac{1}{7}$-approximate Newton steps
    \item The scheme finishes, and the number of approximate Newton steps is bounded by $N_{ns}$. The total complexity of the method is therefore 
\[O((nM + M^3 + M^2Q + \csla)N_{ns}) \text{ in time },\qquad O(n + M^2) \text{ in memory} .\]
\item The returned $\beta$ is statistically optimal: 
\[L(f_\beta) - L(\fstar) \leq \Kfinal \left(\Qc^{\gamma} ~ 
    \Lc^{1-\gamma}\right)^2 ~~\frac{1}{n^{\gamma}} \log \frac{2}{\delta},\]
    where $\Kfinal$ is defined in \cref{thm:quantisource_kernels}.
\end{itemize}
\et

\bpr 
The proof consists mainly of combining \cref{prp:nystrom_select,prp:gen_scheme_optim,thm:quantisource_kernels}. 

Recall that we set $\la = \left(\frac{\Qc^2}{\Lc^2}\frac{1}{n}\right)^{\frac{\alpha}{\alpha(1+2r) + 1}}$. 
\paragraph{1.} Start by defining $\tilde{N}$ such that:
\begin{itemize}
    \item $\tilde{N} \geq N$ where $N$ is defined in \cref{thm:quantisource_kernels};
    \item $\forall n \geq \tilde{N},~ \la \leq \overline{\|\Hess\|}$. This is possible as $\frac{\alpha}{\alpha(1+2r) + 1}$ is a strictly positive exponent. 
    \item 
    $\forall n \geq \tilde{N},~\frac{19 \btb \vee 1 ~\ka^2}{n}\log \frac{n}{2 \delta} < \la$; this is possible as soon as $\frac{\alpha}{\alpha(1+2r) + 1} < 1$, i.e. this is satisfied since $r>0$;  
    \item $\tilde{N} \geq 405 \btb \vee 1 ~\ka^2 \vee 67 \btb \vee 1 ~\ka^2 \log \frac{12 \btb \vee 1~~\ka^2}{\delta}$;
    \item $\forall n \geq \tilde{N},~ \frac{\la \sqrt{2}}{1440 (\bts\vee 1)} \leq \|\Cov\|$.
\end{itemize}
We see that such a $\tilde{N}$ can be defined explicitly. 
\paragraph{2.} Combining the assumptions on $\tilde{N}$ with the ones on $M$, we see that all the assumptions of \cref{prp:nystrom_select} are satisfied and thus that with probability at least $1 - \delta$, all the hypotheses for \cref{thm:quantisource_kernels} are satisfied except the bound on $\hnd_{M,\la}(\beta)$.

\paragraph{3.} Applying \cref{prp:gen_scheme_optim}, taking $\sqrt{\varepsilon} =\Qc^{\gamma} ~ 
    \Lc^{1-\gamma} n^{-\gamma/2} $ and $\la = \left(\frac{\Qc^2}{\Lc^2}\frac{1}{n}\right)^{\frac{\alpha}{\alpha(1+2r) + 1}}$, we see that under these hypotheses, 
    \begin{align*}
    N_{ns}:= &2  \left\lfloor \left(3 + 11 \Rad \|\beta_{M,\la}\|\right)\log_2 \left(7 \Rad \|\nabla \hL_M(0)\|\left(\frac{n \Lc^2}{\Qc^2}\right)^{\frac{\alpha}{\alpha(1+2r) + 1}}\right)\right\rfloor \\
    &+  \left\lceil \log_2 \left(\frac{1}{\Rad \Lc}\left(\frac{n \Lc^2}{\Qc^2}\right)^{\frac{r \alpha}{\alpha(1+2r) + 1}}\right)  \right\rceil.
    \end{align*}
    Now we can bound this harshly:
    \[ N_{ns} \leq  \left(7 + 22 \Rad \|\beta_{M,\la}\|\right)\left( \left\lceil \log_2 \left(7 \Rad \|\nabla \hL_M(0)\| \right)\right\rceil + \left\lceil \log_2\frac{n \Lc^2}{\Qc^2} \right\rceil +  \left\lceil \log_2 \frac{1}{\Rad \Lc} \right\rceil\right).\]
    Now bounding $\Rad\|\beta_{M,\la}\| \leq 10 + \Rad \|\fstar\|$, we get
    \[ N_{ns} \leq  \left(227 + 22 \Rad \|\fstar\|\right)\left( \left\lceil \log_2 \left(7 \Rad \|\nabla \hL_M(0)\| \right)\right\rceil + \left\lceil \log_2\frac{n \Lc^2}{\Qc^2} \right\rceil +  \left\lceil \log_2 \frac{1}{\Rad \Lc} \right\rceil\right).\]

\paragraph{4. } Finally, we use a union bound to conclude.
\epr 

\newpage 

\section{Algorithm\label{sec:algorithm}}
\begin{algorithm}
\footnotesize
\fbox{
\begin{minipage}[t]{0.97\textwidth}
{\bf Input:} $(x_i,y_i)_{i=1}^n$, $n \in \N$, $\ell$ loss function, $k$ kernel function and $\la > 0$.\\
{\bf Return:} estimated function $\widehat{g}:\X \to \R$\\
Parameters: $Q, M, T \in \N$, $\mu_0 > 0$, $(q_k)_{k \in \N}$.\\
Fixed parameters: $t = 2$ from \cref{thm:main-thm}, $\tau = 3$ from \cref{prp:approx_ns_kernels} in \cref{app:ans_kernels}.\\
$(\bar{x}_j)_{j=1}^M \leftarrow \texttt{leverage-scores-sampling}((x_i)_{i=1}^n, M, \lambda, k)$\\
${\bf K} \leftarrow \texttt{kernel-matrix}((\bar{x}_j)_{j=1}^M,(\bar{x}_j)_{j=1}^M)$\\
${\bf T} \leftarrow \texttt{cholesky-upper-triangular}({\bf K})$\\
define the function $v(\cdot) = (k(\bar{x}_1,\cdot), \dots, k(\bar{x}_M, \cdot)) \in \R^M$\\
\fbox{
\begin{minipage}[t]{0.95\textwidth}
\textbf{define \texttt{compute-preconditioner}:}\\
{\bf Input:} $\alpha \in \R^M, \la > 0$\\
$c_i \leftarrow \sqrt{\ell^{(2)}(v(x_i)^\top {\bf T}^{-1} \alpha,y_i)}$ for all $i=1, \dots, n$\\
define the function $k'(\circ, \bullet)$ as 
$k'(\circ, \bullet):=  c_\circ \times c_\bullet \times k(x_\circ, x_\bullet) $ for $\circ, \bullet \in \{1,\dots,n\}$\\
$(h_s)_{s=1}^Q \leftarrow \texttt{leverage-scores-sampling}((i)_{i=1}^n, Q, \lambda, k')$\\
${\bf G} \leftarrow \texttt{kernel-matrix}((\bar{x}_j)_{i=1}^M, (x_{h_s})_{s=1}^Q, k)$\\
${\bf H} \leftarrow {\bf T}^{-\top} \times {\bf G} \times \text{diag}((c_{l_h}^2)_{h=1}^Q) \times {\bf G}^\top \times {\bf T}^{-1}$\\
${\bf B} \leftarrow \texttt{cholesky-upper-triangular}(\frac{1}{Q}{\bf H} + \la I)$\\
return $\bf B$
\end{minipage}
}
\fbox{
\begin{minipage}[t]{0.95\textwidth}
\textbf{define \texttt{preconditioned-conj-grad}:}\\
{\bf Input:} $\alpha \in \R^M, \mu >0, r \in \R^M, \tau \in \N, {\bf B} \in \R^{M \times M}$\\
$p \leftarrow r, s_0 \leftarrow \|r\|^2, \beta \leftarrow 0$\\
For $i=1,\dots, \tau$\\
${}\qquad z \leftarrow \mu {\bf B}^{-\top}{\bf B}^{-1}p + \frac{1}{n}\sum_{i=1}^n \ell^{(2)}(v(x_i)^\top {\bf T}^{-1} \alpha, y_i) ~ (v(x_i)^\top {\bf T}^{-1} {\bf B}^{-1} p) ~{\bf B}^{-\top}{\bf T}^{-\top} v(x_i)$ \\ 
${}\qquad a \leftarrow s_0/(p^\top z)$\\
${}\qquad \beta \leftarrow \beta + a p$\\
${}\qquad r \leftarrow r - a z,~~s_1 \leftarrow \|r\|^2$\\
${}\qquad p \leftarrow r + (s_1/s_0) p$\\
${}\qquad s_0 \leftarrow s_1$\\
return $\beta$
\end{minipage}
}
\fbox{
\begin{minipage}[t]{0.95\textwidth}
\textbf{define \texttt{appr-linear-solver}:}\\
{\bf Input:} $\alpha \in \R^M, \mu > 0, g \in \R^M$ \\
${\bf B} \leftarrow \texttt{compute-preconditioner}(\alpha, \mu)$\\
$u \leftarrow \texttt{preconditioned-conjugate-gradient}(\alpha, \mu, {\bf B}^{-\top}g, \tau = 3, {\bf B})$\\
return ${\bf B}^{-1} u$
\end{minipage}
}
\fbox{
\begin{minipage}[t]{0.95\textwidth}
\textbf{define \texttt{approximate-Newton}:}\\
Input: $\alpha_0 \in \R^M, \mu > 0, t\in \N$\\
For $j = 1,\dots, t$\\
${}\qquad g \leftarrow \mu \alpha_{j-1} + \frac{1}{n}\sum_{i=1}^n \ell^{(1)}(v(x_i)^\top {\bf T}^{-1}\alpha_{j-1}, y_i)~{\bf T}^{-\top} v(x_i)$\\
${}\qquad \alpha_j \leftarrow \alpha_{j-1} - \texttt{appr-linear-solver}(\alpha_{j-1}, \mu, g)$\\
return $\alpha_t$
\end{minipage}
}\\
$\alpha_0 \leftarrow 0$\\
For $k \in \N$\\
${}\qquad \alpha_{k+1} \leftarrow \texttt{approximate-Newton}(\alpha_k,\mu_k,t=2)$\\
${}\qquad \mu_{k+1} \leftarrow q_{k+1}\mu_{k}$\\
Stop when $\mu_{k+1} < \la$ and set $\alpha_{last}\leftarrow \alpha_k$\\
$\widehat{\alpha} \leftarrow \texttt{approximate-Newton}(\alpha_{last}, \la, T)$\\
return $\widehat{g}(\cdot):= v(\cdot)^\top {\bf T}^{-1} \widehat{\alpha}$
\end{minipage}
}
\caption{Algorithm efficient non-parametric learning for generalized self-concordant losses with optimal statistical guarantees discussed in \cref{sec:kernels} of the main paper.\label{alg:algo}}
\end{algorithm}
Let $N, M \in \N$ with 
$M \leq N$. In Alg.~\ref{alg:algo}, $\texttt{leverage-scores-sampling}((z_i)_{i=1}^N, M, k, \la)$ returns a subset of $(z_i)_{i=1}^N$ of cardinality $M$ sampled by using (approximate) leverage scores at scale $\la > 0$ and computed using the kernel $k$. An explicit example of an algorithm computing $\texttt{leverage-scores-sampling}$ is in \cite{rudi2018fast}. Moreover $\texttt{kernel-matrix}((x_i)_{i=1}^N, (x'_i)_{i=1}^M, k)$ computes the kernel matrix $K \in \R^{N\times M}$ where $K_{ij} = k(x_i, x'_j)$, with $N, M \in \N$.

\newpage

\section{Experiments\label{sec:exp-appendix}}

We present our algorithm's performance for logistic regression on two large scale data sets: Higgs and Susy. We have implemented our method using pytorch, and performed computations on one node of a Tesla P100-PCIE-16GB GPU. Recall that in the case of logistic regression, $\ell_{(x,y)}(t) = \log(1 + e^{-yt})$. 

In what follows, denote with $n$ the cardinality of the data set and $d$ the number of features of this data set. The error is measured in terms of classification error for both data sets. In both cases, we pre-process the data by substracting the mean and dividing by the standard deviation for each feature. The data sets are the following. 

\paragraph{Susy}($n = 5 \times 10^{6}$, $d = 18$, binary classification). We always use a Gaussian Kernel with $\sigma = 5$ for logistic loss (obtained through a grid search; note that in \cite{Rudi17}, $\sigma = 4$ is used for the square loss), and will always use $10^4$ Nystrom points.  

\paragraph{Higgs}($n = 1.1 \times 10^7$, $d = 28$, binary classification). We then apply a Gaussian Kernel with $\sigma = 5$, as in \cite{Rudi17} (we have also performed a grid search). \\\

For these data sets, we do not have a fixed test set, and thus set apart $20 \%$ of the data set at random to be the test set, and use the rest of the $80\%$ to train the classifier.\\\

In practice, we perform our globally convergent scheme with the following parameters.
\begin{itemize}
    \item We use $Q = M$ uniform random features to compute the pre-conditioner for each approximate Newton step;
    \item In the first phase, we decrease $\mu$ in a very fast way to $\la$ by starting at $\mu =1$ and dividing $\mu$ by 1000 after performing only a single approximate Newton step (using $2$ iterations of conjugate gradient descent);
    \item In the second phase, we perform $10$ approximate Newton steps (each ANS is computed using $8$ iterations of conjugate gradient descent).
\end{itemize}

\paragraph{Selection of $\la$}

In the introduction, we claim that in many a learning problem, the parameter $\la$ obtained through cross validation is often much smaller than the ones obtained in statistical bounds which are usually of order $\frac{1}{\sqrt{n}}$. This leads to very ill conditioned problems. 

\begin{figure}[ht]
\includegraphics[width=0.48\textwidth]{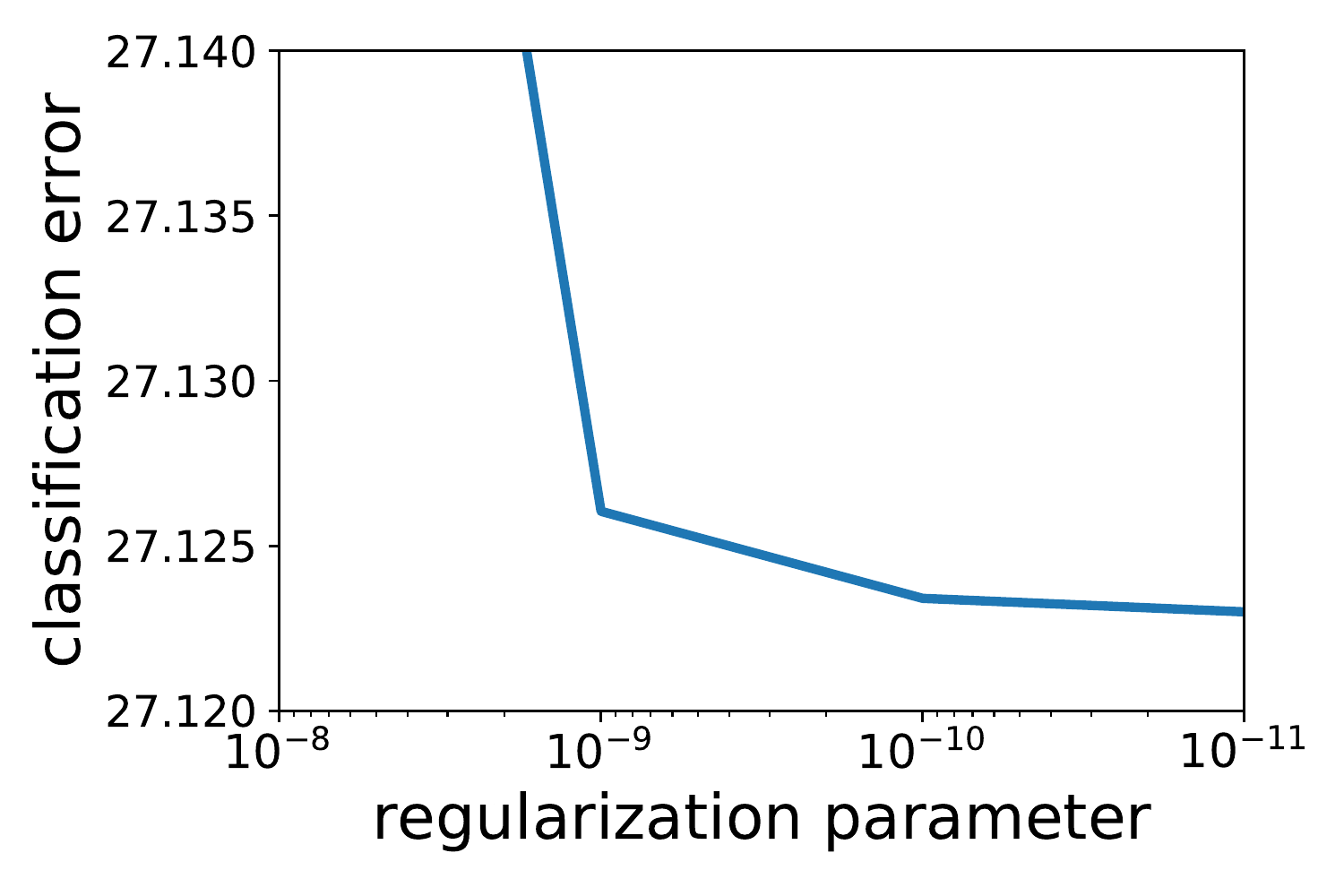}
\hspace{0.25cm}%
\includegraphics[width=0.48\textwidth]{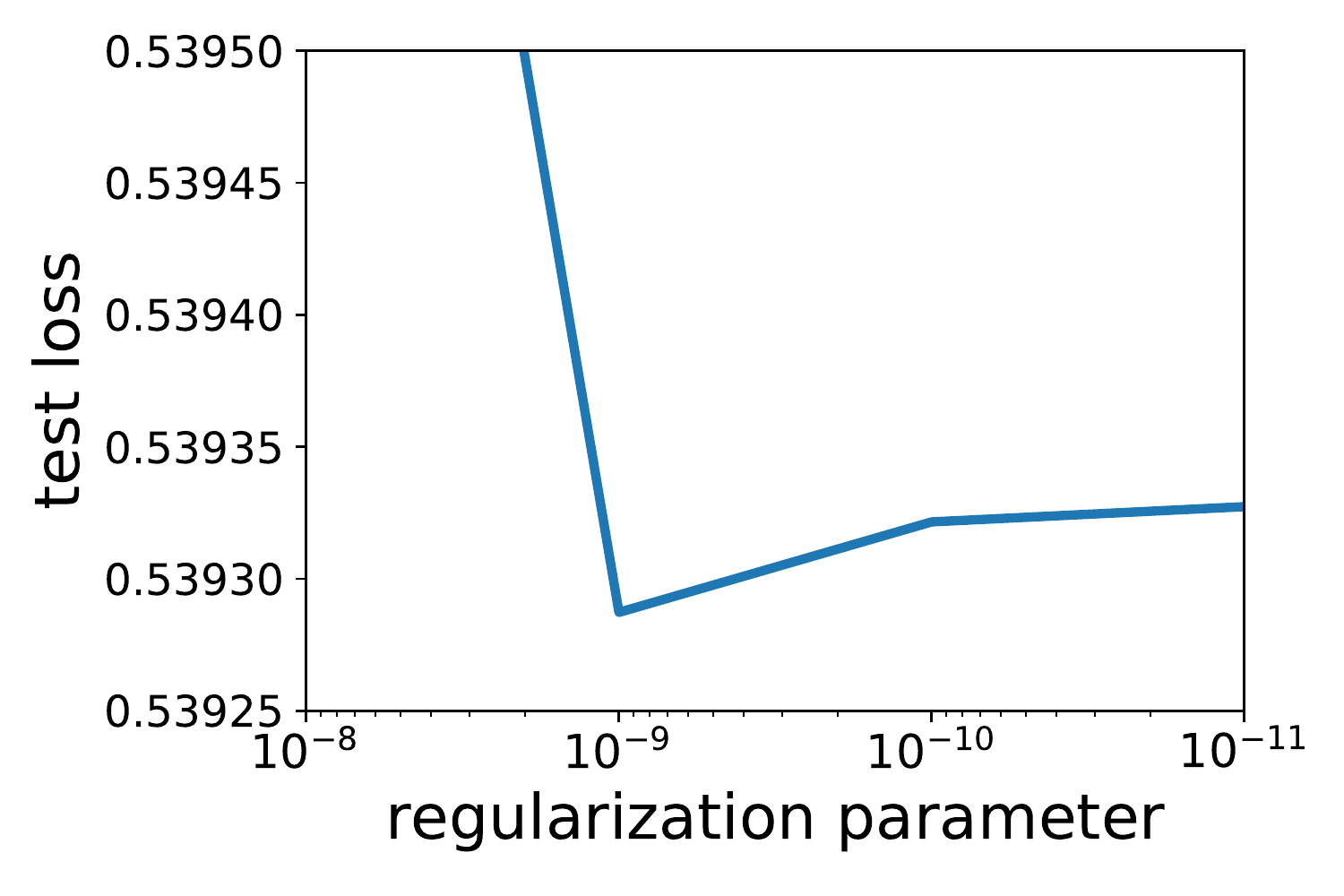} 
\caption{ \textbf{(Left)}  Classification error as a function of the regularization parameter and \textbf{(Right)} test loss as a function of the regularization parameter, when performing a logistic regression with $M = 2\times10^4$ Nystr\"om features on the entire Higgs data set; we select $\la = 10^{-9}$.}
\label{fig:test_higgs}
\end{figure}

\begin{figure}[ht]
\includegraphics[width=0.48\textwidth]{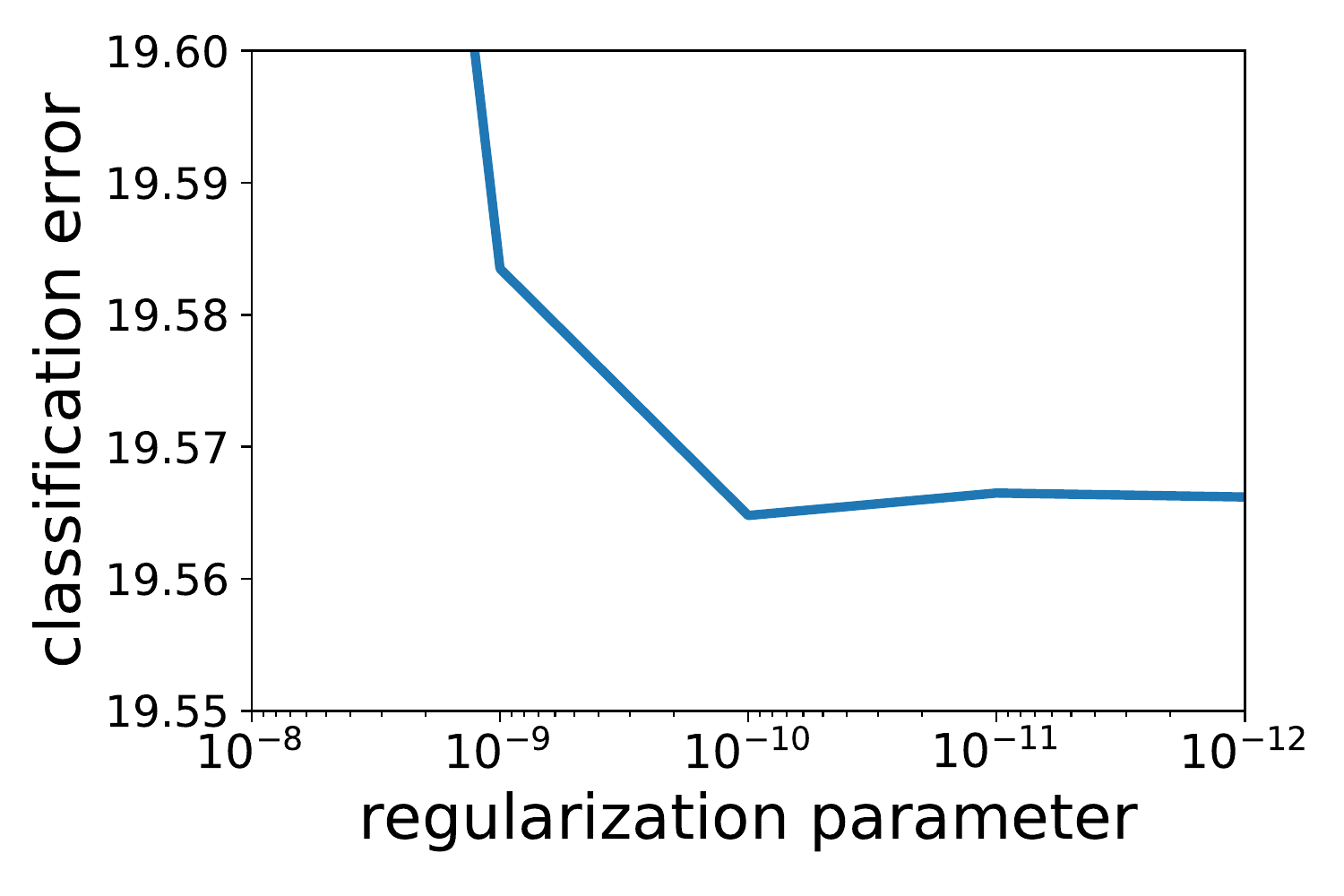}
\hspace{0.25cm}%
\includegraphics[width=0.48\textwidth]{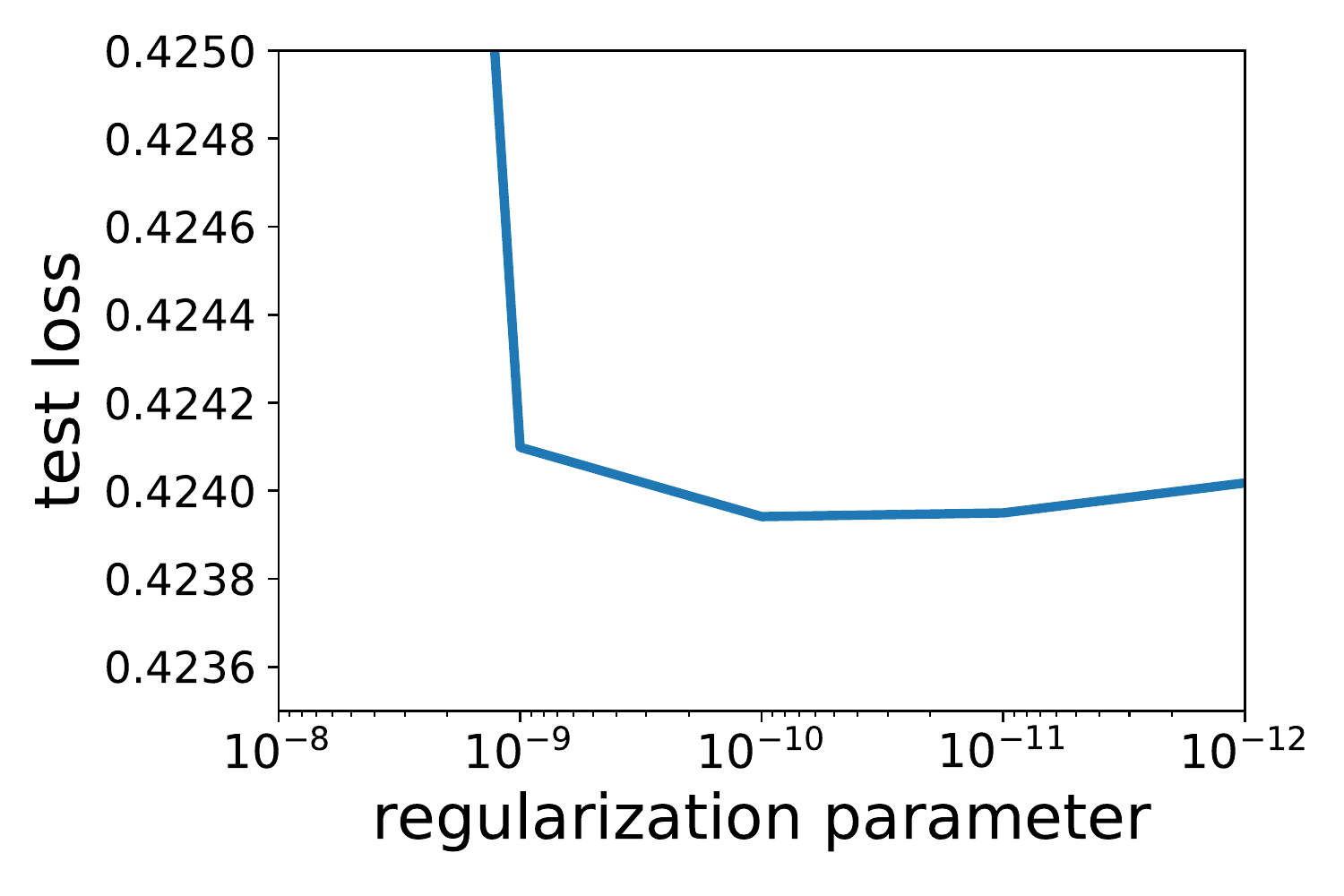} 
\caption{ \textbf{(Left)}  Classification error as a function of the regularization parameter and \textbf{(Right)} test loss as a function of the regularization parameter, when performing a logistic regression with $M = 10^4$ Nystr\"om features on the entire Susy data set; we select $\la = 10^{-10}$.}
\label{fig:test_susy}
\end{figure}

For both data sets, we select $\la$ (and $\sigma$, but we omit the double tables from this paper) by computing the test loss and classification errors for different values of $\la$, and report the evolution of these losses as a function of the parameter $\la$ in \cref{fig:test_higgs} for the Higgs data set, and \cref{fig:test_susy} for the Susy data set. 
We see that the optimal $\la$ yield strongly ill-conditioned problems. 

\paragraph{Comparison with accelerated methods}

Given the $M$ Nystrom points, our aims to minimize $\hL_{M,\la}$. 
From an optimization point of view, i.e. from a point of view where the aim is to minimize $\hL_{M,\la}$, we compare our method with a large mini-batch version of Katyusha accelerated SVRG (see \cite{Allen2017}). \\
Indeed, we perform this method using batch sizes of size $M$; the theoretical bounds provided in \cite{Allen2017} show that the algorithm has linear convergence, with a time complexity of order $O(nM + M^3 + M^2 \sqrt{\frac{L}{\la}})\log \frac{1}{\varepsilon}$ to reach precision $\varepsilon$. In the following plots, we compare both methods in terms of passes and time. \\\

By pass, we mean the following.
\begin{itemize}
    \item In the case of our second-order scheme, we define a pass on the data to be one step of the conjugate gradient descent used to compute approximate newton steps.
    \item In the case of Katyusha SVRG, we define a pass on the data to be either a full gradient computation or $n / M$ computations of the type $K_{\tau M}T^{-1} \beta $ where $T$ is an upper triangular matrix, and $K_{\tau M}$ is a $M \times M$ kernel matrix, associated to one batch gradient. 
\end{itemize}

\begin{figure}[ht]
\includegraphics[width=0.48\textwidth]{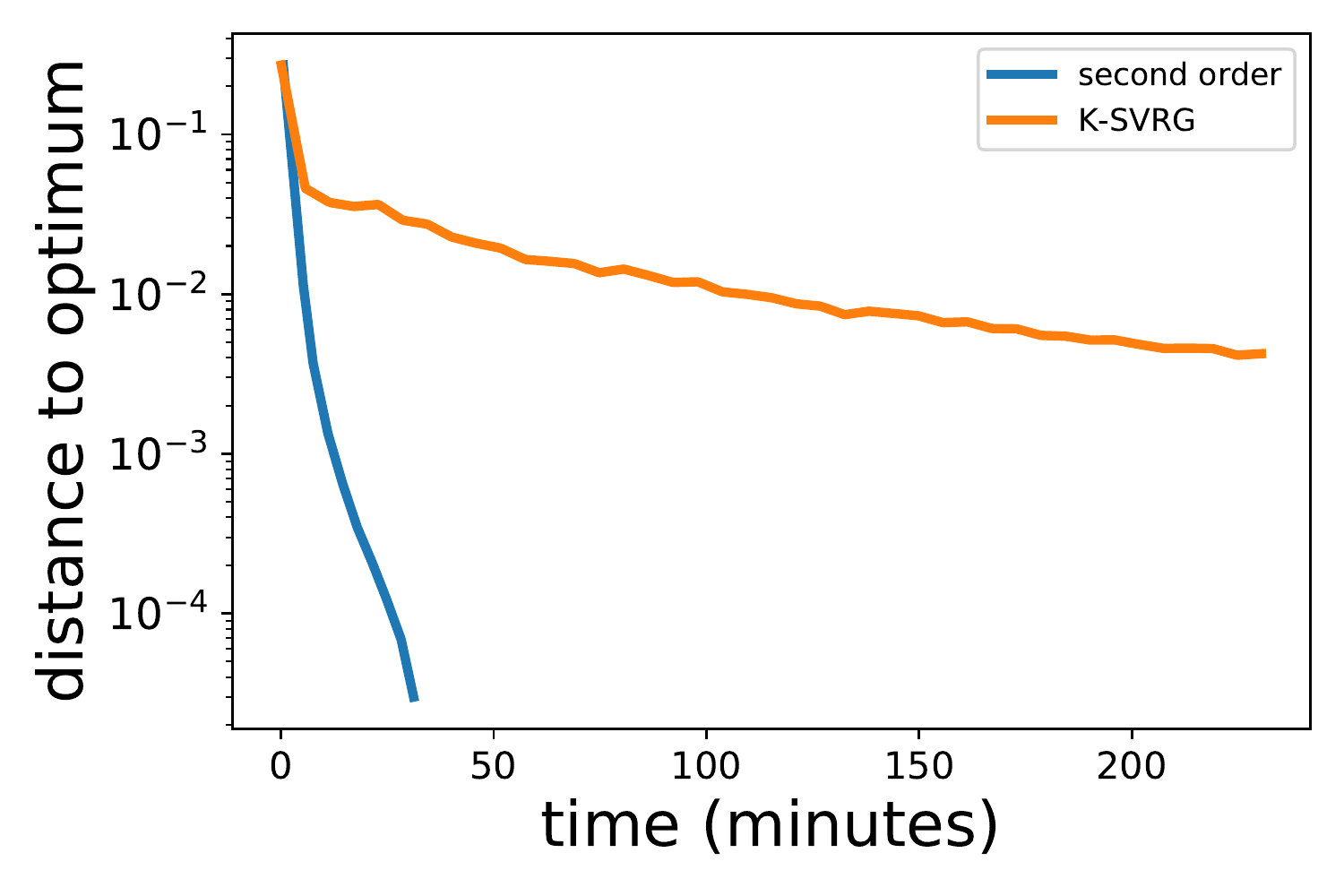}
\hspace{0.25cm}%
\includegraphics[width=0.48\textwidth]{susy_distance_class_small_arxiv.pdf} 
\caption{\textbf{(Left)} Distance to optimum as a function of  time and  \textbf{(Right)} distance to optimum and classification error as a function of the number of passes on the data when performing our second order scheme and K-SVRG to minimize the train loss on Susy, with $1.0\times10^4$ Nystr\"om points and $\la = 10^{-10}$.}
\label{fig:susy_compare_small}
\end{figure}

We use this notion to measure the speed of our method as they both correspond to natural $O(nM)$ operations, and incorporate the essential of the computing time. However, the second point is often much slower to compute than the first, due to the solving of the triangular system. Thus, the notion of passes is to take with precaution, as a pass for the accelerated SVRG algorithm takes much longer to run that a pass for our method. This is confirmed by the time plots (see \cref{fig:susy_compare_small} for in instance).

\begin{figure}[ht]
\includegraphics[width=0.48\textwidth]{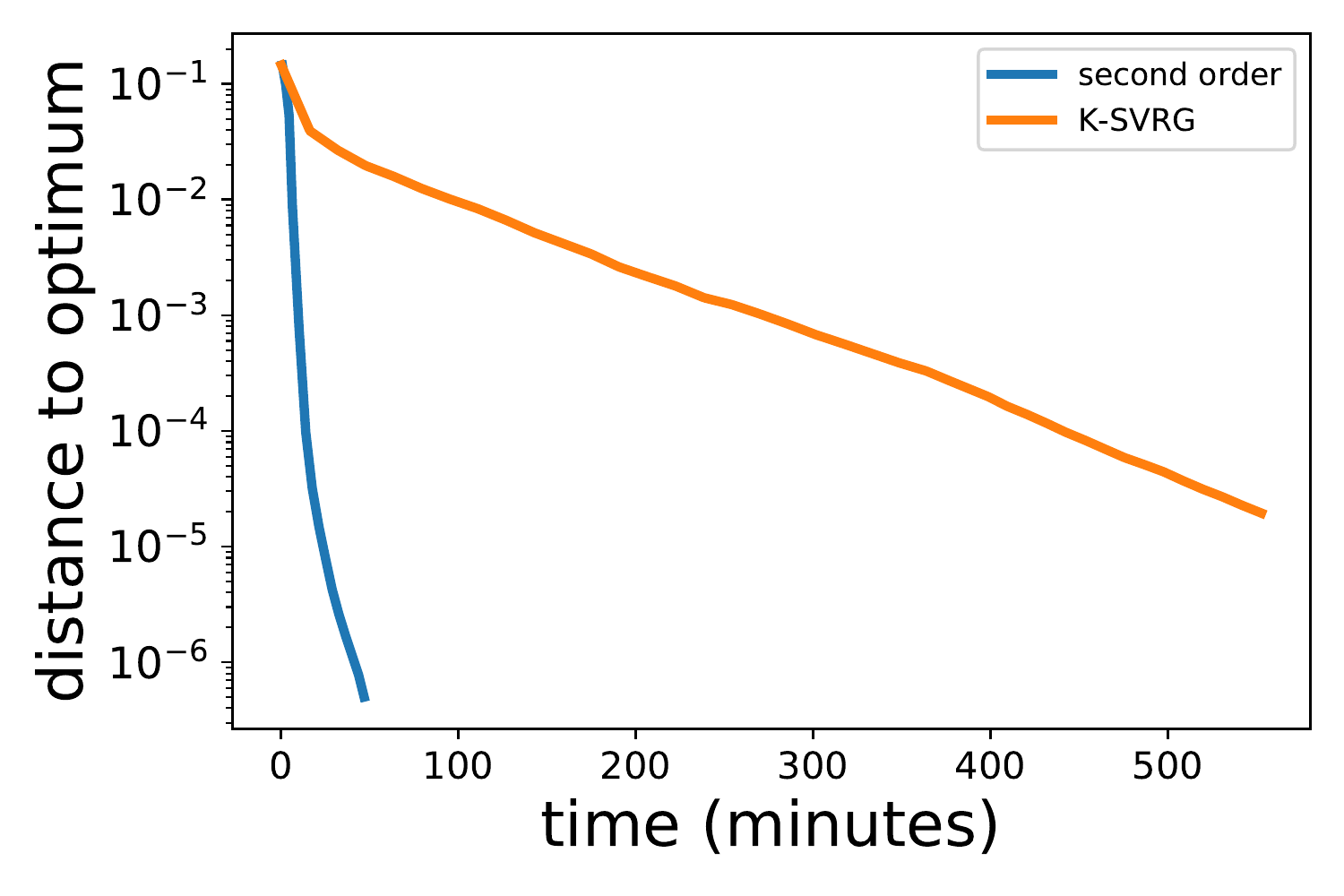}
\hspace{0.25cm}%
\includegraphics[width=0.48\textwidth]{higgs_distance_class_arxiv.pdf} 
\caption{\textbf{(Left)} Distance to optimum as a function of  time and  \textbf{(Right)} distance to optimum and classification error as a function of the number of passes on the data when performing our second order scheme and K-SVRG to minimize the train loss on Higgs, with $1.0\times10^4$ Nystr\"om points and $\la = 10^{-9}$.}
\label{fig:higgs_compare}
\end{figure}

\textit{Comparison between the two methods} - Due to the running time of K-SVRG, we compare both methods for $M = 10000$ Nystr\"om points for both data sets. We compare the performance of these two algorithm with respect to the distance to the optimum in function values as well as classification error \cref{fig:higgs_compare} for the Higgs data set, and in \cref{fig:susy_compare_small} for the Susy data set.

\textit{Note on the need for precise optimization} - As noted in the introduction, we see in both \cref{fig:susy_compare_small} and \cref{fig:higgs_compare} that precise optimization of the objective function is needed in order to get a good classification error. This justifies a posteriori the use of a second order method. In particular, in \cref{fig:susy_compare_small}, one notes the difference in behavior between the two methods : the second order method converges linearly in a fast way while the first order method slows down because of the condition number.

\begin{figure}[ht]
\includegraphics[width=0.48\textwidth]{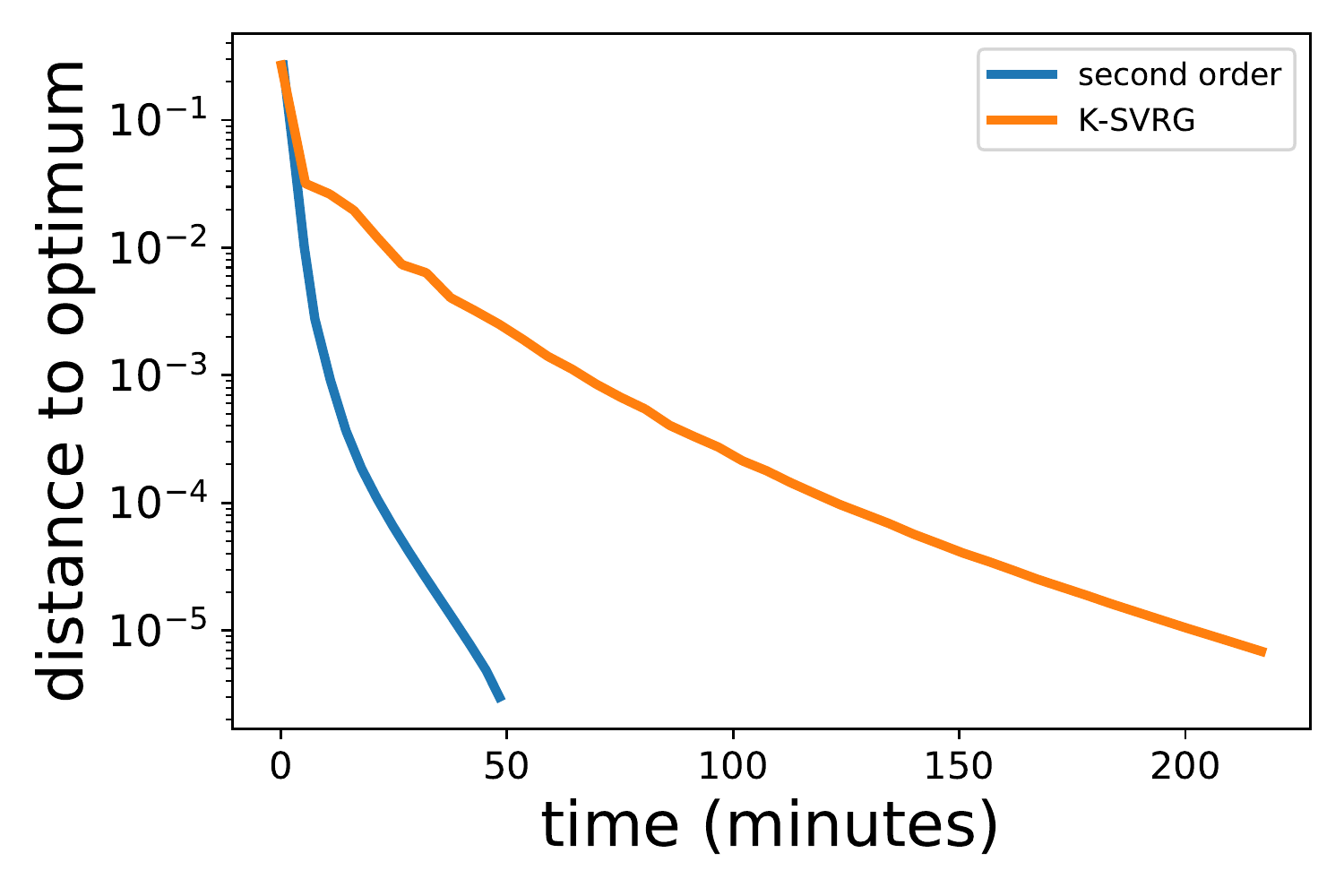}
\hspace{0.25cm}%
\includegraphics[width=0.48\textwidth]{susy_distance_arxiv.pdf} 
\caption{\textbf{(Left)} Distance to optimum as a function of  time and  \textbf{(Right)} distance to optimum and classification error as a function of the number of passes on the data when performing our second order scheme and K-SVRG to minimize the train loss on Susy, with $1.0\times10^4$ Nystr\"om points and $\la = 10^{-8}$.}
\label{fig:susy_compare}
\end{figure}

\textit{Note on ill-conditioning} - First note that in order to optimize test error, one gets very poorly conditioned problems. As predicted by the rates, we observe that K-SVRG is more sensible to ill-conditioning than our second order scheme. Indeed, in \cref{fig:susy_compare}, we have plotted the results for Susy for a smaller condition number with $\lambda = 10^{-8}$, compared to $\lambda = 10^{-10}$ to get optimal test error in \cref{fig:susy_compare_small}. We see that the difference in number of passes needed to reach a certain precision is much lower when $\la = 10^{-8}$ in \cref{fig:susy_compare}, confirming that K-SVRG behaves better when the condition number is smaller.

\paragraph{Performance of our method.}

In \cref{tab:performance},  we record the performance of the following methods, taking the $\la$ values we have obtained previously for the different data sets.  

For FALKON (see \cite{Rudi17}), we take the parameters suggested in the paper (except for the number of Nystr\"om points needed for Higgs, as our computational capacity is limited). 

\begin{table}[h!]
    \centering
    \begin{tabular}{|c|c|c|c|c|c|c|}
 \hline
    \multirow{2}{*}{Method} & \multicolumn{3}{c|}{Susy} &  \multicolumn{3}{c|}{Higgs}\\
 \cline{2-7}
 & c-error & $M$ & time(m) &c-error & $M$ & time(m)\\
 \hline
 Logistic regression with K-SVRG & 19.64\%  & $10^4$ &  230& 27.82 \% & $ 10^4$& 500 \\ 
 \hline 
 
 Logistic regression with our scheme & 19.5\% & $10^4$ & 15 & 26.9 \% & $2.5\times10^4$ & 65\\
 \hline
 Ridge Regression with FALKON (\cite{Rudi17}) & 19.7\% & $10^4$ & 5  & 27.16 \% & $2.5\times10^4$ & 60\\
 \hline
 
\end{tabular}
    \caption{Classification error of different methods}
    \label{tab:performance}
\end{table}

\newpage

\section{\label{app:proj}Solving a projected problem to reduce dimension}

\subsection{Introduction and notations}

In this section, we give ourselves a generalized self-concordant function $f$ whose associated subset we denote with $\G$. Once again, we will always omit the subscript $f$ in the notations associated to $f$.\\\

The aim of this section is the following. Given $f$ and $\la > 0$, computing an approximate solution to 
\[\xla  =\argmin_{x \in \hh}{f_\la(x)},\]
is often too costly. Instead, we look for a solution in a small subset of $\hh$ which we see as the image of a certain orthogonal projector $\Pj$ and which we denote $\hh_\Pj$. Usually, this subset will be finite dimensional and admit an easy parametrization. Thus we will compare an approximation of $\xla$ to an approximation of

\[\xlap = \argmin_{x \in \hh_\Pj}{f_\la(x)} = \argmin_{x \in \hh}{f(\Pj x) + \frac{\la}{2}\|x\|^2}.\]

Denote with $f_\Pj$ the mapping $x \in \hh \mapsto f(\Pj x)$. It is easy to see that, as $f$ is a generalized self-concordant function with $\G$, $f_\Pj$ is naturally a generalized self-concordant with $\G_\Pj:= \Pj \G = \left\{\Pj g ~:~ g \in \G\right\}$. Moreover, $\xlap = x^\star_{f_\Pj,\la}$. \\\

We will adopt the following notations for the quantities related to the generalized self-concordant function $f_\Pj$. Essentially, we always replace $f_\Pj$ simply by $\Pj$ from our definitions in appendix.
\begin{itemize}
    \item For the regularized function :
    \[\forall x \in \hh,~ \forall \la > 0,~ f_{\Pj,\la}(x) = f_\Pj(x) + \frac{\la}{2}\|x\|^2.\]
    \item For the Hessians 
    \[\forall x \in \hh, ~ \la >0,~ \Hess_{\Pj,\la}(x) = \Hess_{f_{\Pj},\la}(x) = \Pj \Hess(\Pj x  ) \Pj + \la \Id.\]
    \item $\forall h \in \hh,~ \tn_\Pj(h):= \tn_{f_\Pj}(h) = \tn(\Pj h)$.
    \item For the Newton decrement:
    \[\forall x \in \hh, ~ \la >0,~  \nd_{\Pj,\la}(x) = \nd_{f_\Pj,\la}(x) = \|\nabla f_{\Pj,\la}\|_{\Hess^{-1}_{\Pj,\la}(x)} = \|\Pj \nabla f(\Pj x) + \la x\|_{\Hess^{-1}_{\Pj,\la}(x)}.\]
    \item For the Dikin ellipsoid radius:
    $$\forall \la>0,~\forall x \in \hh,~ \radd_{\Pj,\la}(x):= \radd_{f_\Pj,\la}(x) = \frac{1}{\sup_{g \in  \G}\|\Pj g\|_{\Hess_{\la,\Pj}^{-1}(x)}};$$
    \item For the Dikin ellipsoid:
    \[\forall \la > 0,~ \forall \cc \geq 0,~ \dik_{\Pj,\la}(\cc):= \dik_{f_\Pj,\la}(\cc).\]
    
\end{itemize}

Note that for any $x \in \hh_\Pj$, $\rpla( x) \geq \rla(x)$.

We will now introduce the key quantities in order to compare an approximation of $\xlap$ to an approximation of $\xla$.

\bd[key quantities]
Define the following quantities
\begin{itemize}
    \item For any $\la >0$, the source term $\srce_\la:= \la\|\xla\|_{\Hess_\la^{-1}(\xla)} = \|\nabla f(\xla)\|_{\Hess_\la^{-1}(\xla)} $;
    \item Given an orthogonal projector $\Pj$, $\la >0$, and $x \in \hh$, the capacity of the projector $\cp(x,\la):= \frac{\|\Hess(x)^{1/2}(\Id - \Pj)\|^2}{\la}$.
\end{itemize}
\ed

\subsection{Relating the projected to the original problem\label{app:proj-to-ideal}}

Given $x \in \hh_\Pj$, our aim is to bound $\nd_\la(x)$ given $\nd_{\la,\Pj}(x)$ and $\srce_\la$. 

\bp \label{prp:link_proj}
Let $x \in \hh_\Pj$. If
\[\frac{\srce_\la}{\rla(\xla)} \leq \frac{1}{4},~ \cp(\xla,\la) \leq \frac{1}{120},~ \nd_{\Pj,\la}(x) \leq \frac{\rpla(x)}{2}, \]
Then it holds: 
\[\nd_\la(x) \leq 3(\nd_{\Pj,\la}(x) + \srce_\la).\]
Moreover, under these conditions, 
\begin{itemize}
    \item $  \|x - \xla\| \leq 7 \la^{-1/2}(\nd_{\Pj,\la}(x) + \srce_\la) $;
    \item $\la \|x\|_{\Hess^{-1}_{\Pj,\la}(x)} \leq 7 \nd_{\Pj,\la}(x) + 9 \srce_\la$.
\end{itemize}

\ep 

\bpr
In this proof, introduce the following auxiliary quantity:

\[\gamma_\la:= \frac{\srce_\la}{\rla(\xla)}.\]
\paragraph{1) Start by bounding $\tn(\Pj\xla - \xla)$.}

It holds:

\eqals{\tn(\Pj x - \xla) &= \sup_{g \in \G}|g \cdot (\Id - \Pj) \xla| \\
& \leq \frac{1}{\rla(\xla)} ~ \|(\Id - \Pj)\xla\|_{\Hess_\la(\xla)}\\
&\leq \frac{1}{\rla(\xla)}~\|\Hess_\la(\xla)^{1/2}(\Id - \Pj) \Hess_\la(\xla)^{1/2}\|~\|\Hess^{-1/2}_\la(\xla)\xla\|\\
& = \left(1 + \cp(\xla,\la)\right) ~\frac{\la \|\Hess^{-1/2}_\la(\xla)\xla\|}{\rla(\xla)}\\
& = (1 + \cp(\xla,\la)) ~\gamma_\la
.}

\paragraph{2) Then bound $\tn(\xlap - \Pj \xla)$}
First, bound $\nd_{\Pj,\la}(\Pj \xla)$: 

\eqals{\nd_{\Pj,\la}(\Pj \xla) & = \|\Pj \nabla f_\la(\Pj \xla) \|_{\Hess_{\la,\Pj}(\Pj \xla)^{-1}} \\
&\leq \| \nabla f_\la(\Pj \xla) \|_{\Hess_{\la}(\Pj \xla)^{-1}} 
.}

Using \cref{eq:hl}, we get $\| \nabla f_\la(\Pj \xla) \|_{\Hess_{\la}(\Pj \xla)^{-1}} \leq e^{\tn((\Id - \Pj)\xla)/2} \nd_\la(\Pj \xla)$. Using \cref{eq:go}, we can bound 
$$\nd_\la(\Pj \xla) \leq \dikins(\tn((\Id - \Pj)\xla))~\|(\Id - \Pj)\xla\|_{\Hess_\la(\xla)} \leq  \dikins(\tn((\Id - \Pj)\xla))~(1 + \cp(\xla,\la))\srce_\la.$$ 
Putting things together,

\[\nd_{\Pj,\la}(\Pj \xla) \leq e^{\tn((\Id - \Pj)\xla)/2} \dikins(\tn((\Id - \Pj)\xla))~(1 + \cp(\xla,\la))\srce_\la.\]
Now 
\[\frac{1}{\rpla(\Pj \xla)}\leq \frac{1}{\rla(\Pj \xla)} \leq  e^{\tn((\Id - \Pj)\xla)/2} \frac{1}{\rla(\xla)}.\]
Hence, 

\[\frac{\nd_{\Pj,\la}(\Pj \xla) }{\rpla(\Pj \xla)}\leq e^{\tilde{t}_\la} \dikins(\tilde{t}_\la)~\tilde{t}_\la,\qquad\tilde{t}_\la = (1 + \cp(\xla,\la))\gamma_\la. \]

Since $t \mapsto e^{t} \dikins(t)~t$ is an increasing function whose value in $0$ is $0$, we find numerically that for $t = \frac{3}{10}$, $e^{t} \dikins(t)~t \leq \frac{1}{2}$. Hence, if $(1 + \cp(\xla,\la))\gamma_\la \leq \frac{3}{10}$, then $\frac{\nd_{\Pj,\la}(\Pj \xla) }{\rpla(\Pj \xla)} \leq \frac{1}{2}$. Using \cref{lm:loc}, this shows that 
$$\tn_\Pj(\Pj \xla - \xlap) = \tn(\Pj \xla - \xlap) \leq \log 2 .$$ 

\paragraph{3) Getting a bound for $\tn(x - \xla)$.} To do so, combine the two previous bounds with the fact that if $\nd_{\Pj,\la}(x) \leq \frac{\rpla(x)}{2}$, then using \cref{lm:loc} with $f_{\Pj}$, $\tn_\Pj(x - \xlap) = \tn(x-\xlap) \leq \log 2$. Thus, if 

\[(1 + \cp(\xla,\la))\gamma_\la \leq \frac{3}{10},~\nd_{\Pj,\la}(x) \leq \frac{\rpla(x)}{2},\]
then it holds
\[\tn(x - \xla) \leq \frac{3}{10} + 2 \log 2.\]
\paragraph{4) A technical result to bound $\|\Hess_\la(x)^{-1/2} \Hess_{\Pj,\la}(x)^{1/2}\|$}. Using the fact that $\Pj x = x$, and \cref{lm:proj_hermit}, applied to $\Ab = \Hess(x)$, we get 
\[\|\Hess_\la(x)^{-1/2} \Hess_{\Pj,\la}(x)^{1/2}\| \leq  1 + \sqrt{\cp(x,\la)}.\]

Then, one can easily bound $\cp(x,\la) \leq e^{\tn(x-\xla)}\cp(\xla,\la)$. 
\paragraph{5) Let us now bound $\nd_\la(x)$.}

First, decompose the term 
\[\nd_\la(x) = \|\nabla f_\la(x)\|_{\Hess^{-1}_\la(x)} \leq  \|\Pj \nabla f_\la(x)\|_{\Hess^{-1}_\la(x)} + \|(\Id - \Pj)\nabla f(x)\|_{\Hess^{-1}_\la(x)} .\]

Since $x \in \hh_\Pj$,  $\|\Pj \nabla f_\la(x)\|_{\Hess^{-1}_\la(x)} = \|\nabla f_{\Pj,\la}(x)\|_{\Hess^{-1}_\la(x)}$, and using the previous point, we get

\[\|\Pj \nabla f_\la(x)\|_{\Hess^{-1}_\la(x)} \leq  \left(1 + e^{\tn(x - \xla)/2}\sqrt{\cp(\xla,\la)}\right) \nd_{\Pj,\la}(x).\]

Let us now bound the second term. We divide it into two terms: 
\[\|(\Id - \Pj)\nabla f(x)\|_{\Hess^{-1}_\la(x)} \leq \|(\Id - \Pj)\left(\nabla f(x) - \nabla f(\xla)\right)\|_{\Hess^{-1}_\la(x)} + \|(\Id - \Pj)\nabla f(\xla)\|_{\Hess^{-1}_\la(x)}.\]
The second term can be bounded in the following way: 
\[\|(\Id - \Pj)\nabla f(\xla)\|_{\Hess^{-1}_\la(x)} \leq \frac{1}{\sqrt{\la}} \|(\Id - \Pj)\Hess^{1/2}_\la(\xla)\|~\|\nabla f(\xla)\|_{\Hess^{-1}_\la(\xla)} \leq  \sqrt{1 + \cp(\xla,\la)} ~\srce_\la.\]
For the first term, we proceed in the following way. 
\eqals{
\|(\Id - \Pj)\left(\nabla f(x) - \nabla f(\xla)\right)\|_{\Hess^{-1}_\la(x)} & = \| \int_{0}^1{\Hess_\la^{-1/2}(x)(\Id - \Pj) \Hess(x_t)(x-\xla)~dt} \| \\
& \leq \frac{1}{\sqrt{\la}}~\int_{0}^1{\|(\Id - \Pj) \Hess^{1/2}(x_t)\|~ \|\Hess^{1/2}(x_t)(x-\xla)\|~dt}\\
& \leq \sqrt{ \cp(\xla,\la)}~\dikins(\tn(x - \xla))~\|x - \xla\|_{\Hess(\xla)}\\
& \leq \sqrt{ \cp(\xla,\la)} ~e^{\tn(x - \xla)} \nd_\la(x)
.}
Hence the final bound: 

\[\left(1 - \sqrt{ \cp(\xla,\la)} ~e^{\tn(x - \xla)}\right)\nd_\la(x) \leq \left(1 + e^{\tn(x - \xla)/2}\sqrt{\cp(\xla,\la)}\right) \nd_{\Pj,\la}(x) +  \sqrt{1 + \cp(\xla,\la)} ~\srce_\la.\]

Now if $\cp(\xla,\la) \leq \frac{1}{120}$, we see that $\sqrt{ \cp(\xla,\la)} ~e^{\tn(x - \xla)} \leq \frac{1}{2}$, and hence, using the bound on $\tn(x - \xla)$, 

\[\nd_\la(x) \leq 3(\nd_{\Pj,\la}(x) + \srce_\la).\]

\paragraph{6) Showing the last two points}. We leverage the fact that $\nd_\la(x) \leq 3(\nd_{\Pj,\la}(x) + \srce_\la)$ and $\tn(x - \xla) \leq \frac{3}{10} + 2 \log 2$.\\
To show the first bound, we plug in the previous results in the following equation: 
\[\|x - \xla\| \leq \la^{-1/2}\|x - \xla\|_{\Hess_\la(x)} \leq \frac{1}{\dikin(\tn(x - \xla))} ~ \la^{-1/2} \nd_\la(x).\]
The last inequality is obtained using \cref{eq:gl}.\\

To show the second point, we use the fact that $x \in \hh_\Pj$ to show that 
\[\la \|x\|_{\Hess_{\Pj,\la}^{-1}(x)} \leq  \la \|x\|_{\Hess^{-1}_{\la}(x)} \leq \la \|x - \xla\|_{\Hess_\la(x)} + \la\|\xla\|_{\Hess^{-1}_\la(x)}.\]
Then applying \cref{eq:hl} and \cref{eq:gl}:
\[\la \|x\|_{\Hess_{\Pj,\la}^{-1}(x)} \leq \frac{1}{\dikin(\tn(x - \xla))} ~ \nd_\la(x) + e^{\tn(x-\xla)/2} \srce_\la.\]
We then use the previous results to conclude.
\epr

\subsection{Finding a good projector}

\blm 
If for a certain $\eta \leq \la$ and for a certain constant $C$, $\|\Hess_\eta^{1/2}(x)(\Id-\Pj)\|^2 \leq C \eta$, then 
$$\cp(x,\la) \leq \frac{C\eta}{\la}.$$
\elm

\bpr 
This is completely direct, using the fact that $\Hess^{1/2}(x)\preceq \Hess_\eta^{1/2}(x)$.
\epr

\newpage 

\section{Relations between statistical problems and empirical problem.\label{app:proj_stat}}

In this section, we recall and reformulate the framework from \cite{marteau2019}.

\subsection{Statistical problem and ERM estimator}

Let $\Z$ be a Polish space and $Z$ be a random variable on $\Z$ with distribution $\rho$. 
Let $\hh$ be a separable Hilbert space, with norm $\|\cdot\|$, and let $(f_z)_{z \in \Z}$ be a family of functions on $\hh$.
Our goal is to minimize the \textit{expected risk} with respect to $x \in \hh$:

$$ \inf_{x \in \hh}~f(x):= \Exp{f_Z(x)}.$$ 
Given $(z_i)_{i=1}^n \in \Z^n$, we define the \textit{empirical risk}:

\[\hf(x):= \frac{1}{n}\sum_{i=1}^n f_{z_i}(x)  ,\]
and consider the following estimator based on regularized empirical risk minimization given $\la>0$ (note that the minimizer is unique in this case):

\[
\widehat{x}_\la^\star = \argmin_{x \in \hh} \hf_\la(x):= \hf(x)  + \frac{\la}{2} \|x\|^2,
\]
where we assume the following.
\begin{restatable}[i.i.d. data]{ass}{asmiid}
\label{asm:iid}
The samples $(z_i)_{1 \leq i \leq n}$ are independently and identically distributed according to $\rho$. 
\end{restatable}
We make the following assumption on the family $(f_z)$ (this is a reformulation of Assumption 8 in \cite{marteau2019})
\begin{restatable}[Generalized self-concordance]{ass}{asmgensc}
\label{asm:gen_sc}
For any $z \in \Z$, there exists an associated subset $\G_z \subset \hh$ such that $(f_z,\G_z)$ is generalized self-concordant in the sense of \cref{df:genscc}. 
\end{restatable}

Moreover we require the following technical assumption to guarantee that $f$ and and its derivatives are well defined for any $x \in \hh$ (this is a reformulation of Assumptions 3 and 4 in \cite{marteau2019}, and the necessary conditions to obtain \cref{prp:expectancy}).

\begin{restatable}[Technical assumptions]{ass}{asmtech}
\label{asm:asmtech}
The mapping $(z,x) \in \Z \times \hh \mapsto f_z(x)$ is measurable. Moreover, 
\begin{itemize}
    \item the random variables $\|f_Z(0)\|, \|\nabla f_Z(0)\|, \Tr(\nabla^2 f_Z(0))$ are are bounded;
    \item $\G:= \bigcup_{z \in \supp(Z)}\G_z$ is a bounded subset of $\hh$.
\end{itemize}
\end{restatable}
The assumptions above are usually easy to check in practice. In particular, if the support of~$\rho$ is bounded, the mappings $z \mapsto \ell_z(0),\nabla \ell_z(0), \Tr(\nabla^2\ell_z(0))$ are continuous, and $z \mapsto \G_z$ is uniformly bounded on bounded sets, then they hold. \\\

\bp \label{prp:scProblems}

Under \cref{asm:gen_sc,asm:asmtech}, the function $(f,\G)$ (or simply $f$) is generalized self-concordant. 

Moreover, under \cref{asm:iid}, define 
\[\hG:= \bigcup_{i=1}^n{\G_{z_i}}.\]
Then $(\hf,\hG)$ (or simply $\hf$) is generalized self-concordant. Moreover, note that $\hG \subset \G$. 

\ep

The main regularity assumption we make on our statistical problems follows (see Assumption 5 in \cite{marteau2019}).

\begin{restatable}[Existence of a minimizer]{ass}{asmminimizer}
\label{asm:minimizer}
There exists $x^\star \in \hh$ such that $ f(x^\star) = \inf_{x \in \hh} f(x).$
\end{restatable}

\paragraph{Notations}

We adopt all the notations from \cref{app:sc} for $f$ and $\hf$, which are generalized self-concordant functions with associated subsets given in \cref{prp:scProblems} with the following conventions:

\begin{itemize}
    \item For all quantities relating to $f$, we omit the subscript $f$ as usual;
    \item For all quantities relating to $\hf$, we omit the subscript $\hf$ and instead put a hat over all these quantities. For instance:
    \[\hHess(x):= \Hess_{\hf} (x) = \frac{1}{n}\sum_{i=1}^n{\nabla^2f_{z_i}(x)},~ \hrla(x):= \radd_{\hf,\la}(x) = \frac{1}{\sup_{g \in \hG}{\|g\|_{\hHess^{-1}_\la(x)}}}, \text{ etc...}\]
\end{itemize}

Recall the two main quantities introduced in \cite{marteau2019} to establish the quality of our estimator $\hxla$ (in \cite{marteau2019}, this is a mix between Proposition 2 and Definition 3).
\bp[Bias, degrees of freedom] \label{df:main_quant_stat} Suppose  \cref{asm:gen_sc,asm:asmtech,asm:minimizer} are satisfied. The following key quantities are well defined: 
\begin{itemize}
    \item the \emph{bias} $\bias_\la = \|\Hess_\la(\xo)^{-1/2} \nabla f_\la(\xo)\|$;
    \item the \emph{effective dimension} $\Deff_\la = \Exp{\|\Hess_\la(\xo)^{-1/2}\nabla f_Z(\xo)\|^2}$.
\end{itemize}
Moreover, we also introduce the following quantities:
\[\Bos:= \sup_{z \in \supp(Z)}{\|\nabla f_z(\xo)\|},\qquad \Bts:= \sup_{z \in \supp(Z)}{\Tr(\nabla^2 f_z(\xo))},\qquad \Qs = \frac{\Bos}{\sqrt{\Bts}}.\]
\ep

We can now recall the main theorem of \cite{marteau2019} (Theorem 4), which quantifies the behavior of the ERM estimator: 

\bt[Bound for the ERM estimator]\label[theorem]{thm:general-result}
Let $n \in \N$, $\delta \in (0,1/2]$, $0 < \la \leq \Bts$. Whenever
\[n \geq \triangle_1 \frac{\Bts}{\la}\log \frac{8 \square_1^2\Bts}{\la \delta},\qquad \sqrt{\triangle_2~  \frac{\Deff_\lambda \vee (\Qs)^2}{n} ~\log \frac{2}{\delta}}\leq \rla(\xo) ,\qquad 2\bias_\la \leq \rla(\xo),\]
 then with probability at least $1-2 \delta$, it holds
\eqal{\label{eq:general-result-bound}
f(\hxla) - f(\xo) \leq \mathsf{C}_{\textup{bias}}~\bias_\la^2 + \mathsf{C}_{\textup{var}}~\frac{\Deff_\la \vee (\Qs)^2}{n}~\log \frac{2}{\delta},
}
where $\mathsf{C}_{\textup{bias}},\mathsf{C}_{\textup{var}}, \square_1 \leq 414,~\triangle_1,~\triangle_2 \leq 5184$.
\et

\subsection{Link between a good approximation of $\hxla$ and $\xo$}

In this paper, we provide an algorithm which can effectively compute a good approximation of $\hxla$ (as it is a finite sum problem which can be solved). This algorithm will return a certain $x \in \hh$, whose precision with respect to the empirical problem will be characterized by $\hnd_\la(x)$. The aim of the following lemma is to see how this approximation $x$ behaves with respect to the statistical problem. 

\blm\label{lm:approx_emp_pb} Suppose the conditions for \cref{thm:general-result} are satisfied, i.e.
let $n \in \N$, $\delta \in (0,1/2]$, $0 < \la \leq \Bts$ and suppose
\[n \geq \triangle_1 \frac{\Bts}{\la}\log \frac{8 \square_1^2\Bts}{\la \delta},\qquad  \sqrt{\triangle_2~  \frac{\Deff_\lambda \vee (\Qs)^2}{n} ~\log \frac{2}{\delta}} \leq \rla(\xo), \qquad 2\bias_\la \leq \rla(\xo).\]
Let $x$ be an approximation of $\hxla$ characterized by its Newton decrement $\hnd_\la(x)$. If 
\[\hnd_\la(x) \leq \frac{\hrla(x)}{2},~ \hnd_\la(x) \leq \frac{\rla(\xo)}{2},\]
then with probability at least $1-2 \delta$, it holds

\[f(x) - f(\xo) \leq 14(f(\hxla) - f(\xo)) + 30\hnd_\la(x)^2.\]

\elm
\bpr 
Using \cref{eq:fv}, 
\eqals{
f(x) - f(\hxla) &\leq \langle \nabla f(\hxla),x - \hxla\rangle_{\hh} + \psi(\nm{x - \hxla})\|x - \hxla\|_{\Hess_\la(\hxla)}^2\\
& \leq \frac{1}{2}\| \nabla f(\hxla)\|^2_{\Hess^{-1}_\la(\hxla)} + \left(\psi(\nm{x - \hxla}) + \frac{1}{2}\right)\|x - \hxla\|^2_{\Hess_\la(\hxla)}
.}

\paragraph{1. Let us bound $\|\nabla f(\hxla)\|_{\Hess^{-1}_\la(\hxla)}$}

\eqals{
\|\nabla f(\hxla)\|_{\Hess^{-1}_\la(\hxla)} & \leq \int_{0}^1{\|\Hess^{-1/2}_\la(\hxla)\Hess(x_t)(\hxla - \xo)\|~dt}, &x_t = (1-t)\hxla + t\xo \\
&\leq \int_{0}^1{\|\Hess^{-1/2}_\la(\hxla)\Hess^{1/2}(x_t)\| ~\|\Hess^{1/2}(x_t)(\hxla - \xo)\|~dt}.}
Now using equation \cref{eq:h}
\[\Hess(x_t) \preceq e^{t \nm{\hxla - \xo}}\Hess(\hxla),\qquad \Hess(x_t) \preceq e^{(1-t)\nm{\hxla - \xo}} .\]
Thus:

\[\|\nabla f(\hxla)\|_{\Hess^{-1}_\la(\hxla)} \leq e^{\nm{\hxla - \xo}/2} ~ \|\hxla - \xo\|_{\Hess(\xo)} .\]

Finally, using equation \cref{eq:fv}
\[ \|\nabla f(\hxla)\|_{\Hess^{-1}_\la(\hxla)} \leq \frac{e^{\nm{\hxla - \xo}/2} }{\psi(-\nm{\hxla - \xo})^{1/2}}\left( f(\hxla) - f(\xo)\right)^{1/2}.
\]

\paragraph{2. Let us bound the terms involving $ \|x - \hxla\|_{\Hess_\la(\hxla)} $}

Note that using \cref{eq:gl} and \cref{eq:hl} applied to $\hf$, 

\[\|x - \hxla\|_{\Hess_\la(\hxla)}  \leq \|\Hess^{1/2}_\la(\hxla) \hHess^{-1/2}_\la(\hxla)\|~\frac{e^{\hnm{x - \hxla}/2}}{\dikin(\hnm{x - \hxla})} \hnd_\la(x).\]

This also leads to:

\eqals{
\nm{x - \hxla} &\leq \frac{1}{\rla(\hxla)}~\|\Hess^{1/2}_\la(\hxla) \hHess^{-1/2}_\la(\hxla)\|~\|x - \hxla\|_{\hHess_\la(\hxla)}\\
& \leq \frac{1}{\rla(\hxla)}~\|\Hess^{1/2}_\la(\hxla) \hHess^{-1/2}_\la(\hxla)\|~\frac{e^{\hnm{x - \hxla}/2}}{\dikin(\hnm{x - \hxla})} \hnd_\la(x).
}

\paragraph{3. Putting things together}

In the end, we get 

\eqals{f(x) - f(\xo) &\leq \left(1 + \frac{e^{\nm{\hxla - \xo}} }{\psi(-\nm{\hxla - \xo}) }\right)(f(\hxla)  - f(\xo))\\
& + \left(\psi(\nm{x - \hxla}) + \frac{1}{2}\right)\left(e^{\nm{\hxla - \xo_\la}/2}\|\Hess^{1/2}_\la(\xo_\la) \hHess^{-1/2}_\la(\xo_\la)\|\frac{e^{\hnm{x - \hxla}/2}}{\dikin(\hnm{x - \hxla})}\right)~\hnd_\la(x)^2
.}

Moreover, we bound 
\[\nm{x - \hxla} \leq e^{(\nm{\xo - \hxla} + \nm{\hxla - \xo_\la})/2}~\|\Hess^{1/2}_\la(\xo_\la) \hHess^{-1/2}_\la(\xo_\la)\|~\frac{e^{\hnm{x - \hxla}/2}}{\dikin(\hnm{x - \hxla})} \frac{\hnd_\la(x)}{\rla(\xo)}.\]

\paragraph{4. Plugging in previous results}

Under the assumptions of this lemma, which include the assumptions of Theorem 4. in \cite{marteau2019}, we get the following bounds.
\begin{itemize}
    \item In \cite{marteau2019},the assumptions of Theorem 4 imply that we can use Lemma 9, which uses Lemma 8 in which we show that with probability at least $1-\delta$, 
    \[\|\hHess_\la^{-1/2}(\xla)\Hess_\la(\xla)^{1/2}\|^2 \leq 2.\]
    \item Still using the assumptions of Theorem 4, we see in the proof of this theorem that the assumptions of Theorem 7 of \cite{marteau2019} are satisfied in the case where $\bias_\la \leq \frac{\rla(\xo)}{2}$, and thus that 
    $$\tn(\hxla - \xo_\la) \leq \log 2,~ \tn(\xo_\la - \xo) \leq \log 2.$$
    
\end{itemize}

Plugging in all these bounds, we get

\[\left(1 + \frac{e^{\nm{\hxla - \xo}} }{\psi(-\nm{\hxla - \xo}) }\right) \leq 14,~\nm{x - \hxla} \leq 6,\]
\[\left(\psi(\nm{x - \hxla}) + \frac{1}{2}\right)\left(e^{\nm{\hxla - \xo_\la}/2}\|\Hess^{1/2}_\la(\xo_\la) \hHess^{-1/2}_\la(\xo_\la)\|\frac{e^{\hnm{x - \hxla}/2}}{\dikin(\hnm{x - \hxla})}\right) \leq 30.\]
\epr

\subsection{Bounds when we solve a projected empirical problem \label{app:bounds_proj}}

In this section, we place ourselves in the setting of \cref{app:proj}. In this section, we had argued that for computational purposes, it was less costly to compute an approximate solution to a projected problem.\\\

In this section, we assume that we are going to project the regularized empirical problem, that is solve approximately

\[ x \approx \argmin_{x \in \hh}{\hf_{\Pj,\la}(x)} = \hf(\Pj x) + \frac{\la}{2}\|x\|^2.\]

for a given orthogonal projection $\Pj$. Recall from \cref{app:proj} that there is a natural way of seeing $\hf_\Pj$ as a generalized self-concordant function. We import all the notations from this section, keeping a $\widehat{\cdot}$ over all notations to mark the fact that we are projecting $\hf$ and not $f$. 

To quantify the quality of the approximation $x$, we will use the Newton decrement for the empirical projected problem $\hnd_{\Pj,\la}(x):=\nd_{\hf_\Pj,\la}(x)$.\\\

As we see in \cref{prp:link_proj}, under certain conditions, bounding $\hnd_\la(x)$ amounts to bounding two terms:

\begin{itemize}
    \item The empirical source $\hsrce_\la:= \la \|\hxla\|_{\hHess_\la^{-1}(\hxla)}$,
    \item The projected empirical Newton decrement $\hnd_{\Pj,\la}(x)$.
\end{itemize}

\paragraph{1. Bounding the empirical source term $\hsrce_\la$} Start by bounding the source empirical source term using quantities we know.

\blm[Empirical source]\label{lm:emp_source}
Let $n \in \N$, $\delta \in (0,1/2]$, $0 < \la \leq \Bts$. Whenever
\[n \geq \triangle_1 \frac{\Bts}{\la}\log \frac{8 \square_1^2\Bts}{\la \delta},\qquad  \sqrt{\triangle_2~  \frac{\Deff_\lambda \vee (\Qs)^2}{n} ~\log \frac{2}{\delta}} \leq \rla(\xo), \qquad 2\bias_\la \leq \rla(\xo).\]
The following holds, with probability at least $1 - 2 \delta$.

\[\hsrce_\la \leq 8~ \bias_\la + 80~ \sqrt{\frac{\Deff_\la \vee (\Qs)^2~\log \frac{2}{\delta}}{n}}.\]

Moreover, we also have the following bound :
\[\|\hxla - \xo\| \leq 3~\la^{-1/2}~\bias_\la + 8~\la^{-1/2} ~\sqrt{ \frac{\Deff_\la \vee (\Qs)^2~\log \frac{2}{\delta}}{n}}.\]
\elm

\bpr

We first decompose the source term into two terms, and then apply different bounds from \cite{marteau2019} to effectively bound it. We will use the following quantity:
\[\bvar_\la:=  \|\Hess_\la^{1/2}(\xla) \hHess_\la^{-1/2}(\xla)\|^2~\|\nabla \hf_\la(\xla)\|_{\Hess_\la^{-1}(\xla)}.\]

It is also defined in equation (23) in \cite{marteau2019}.

\paragraph{1. Dividing $\hsrce_\la$ into two controllable terms}. Decompose
\eqals{\hsrce_\la = \|\la \hxla\|_{\hHess_\la^{-1}(\hxla)} & \leq \|\hHess_\la^{-1/2}(\hxla)\Hess_\la^{1/2}(\hxla)\| ~ \|\la \hxla\|_{\Hess_\la^{-1}(\hxla)}\\
& \leq 
\|\hHess_\la^{-1/2}(\hxla)\Hess_\la^{1/2}(\hxla)\| ~  \left(\|\nabla f_\la(\hxla)\|_{\Hess^{-1}_\la(\hxla)} + \|\nabla f(\hxla)\|_{\Hess^{-1}_\la(\hxla)}\right)
.}

On the one hand, from the previous proof, we get

\eqals{\|\nabla f(\hxla)\|_{\Hess^{-1}_\la(\hxla)} &\leq e^{\nm{\hxla - \xo}/2} ~ \|\hxla - \xo\|_{\Hess(\xo)} \\
& \leq e^{\nm{\hxla - \xo}/2} ~\left( e^{\nm{\xo_\la - \xo}}\|\hxla - \xo_\la\|_{\Hess_\la(\xo_\la)} + \|\xo_\la - \xo \|_{\Hess_\la(\xo)}\right)\\
&\leq e^{\nm{\hxla - \xo}/2} ~\left( \frac{e^{\nm{\xo_\la - \xo}}}{\dikin(\nm{\hxla - \xo_\la})}\bvar_\la + \frac{1}{\dikin(\nm{\xo_\la - \xo})}\bias_\la\right)
.}

In the last line, we use the fact that 
$\|\hxla - \xo_\la\|_{\Hess_\la(\xo_\la)} \leq \| \Hess^{1/2}_\la(\xo_\la)\hHess^{-1/2}_\la(\xo_\la)\| ~\|\hxla - \xo_\la\|_{\hHess_\la(\xo_\la)}$ and then bound it using \cref{eq:gl} applied to $\hf$ to get 
\eqals{\|\hxla - \xo_\la\|_{\hHess_\la(\xo_\la)} &\leq 
\frac{1}{\dikin(\htn(\xla - \hxla))}\|\nabla \hf_\la(\xla)\|_{\hHess_\la^{-1}(\xla)} \\
& \leq \frac{1}{\dikin(\tn(\xla - \hxla))} \| \Hess^{1/2}_\la(\xo_\la)\hHess^{-1/2}_\la(\xo_\la)\|~\|\nabla \hf_\la(\xla)\|_{\Hess_\la^{-1}(\xla)}
.}

On the other hand, apply successively \cref{eq:gl} to $f$ and $\hf$ using the fact that $\htn \leq \tn$ to get 

\eqals{\|\nabla f_\la(\hxla)\|_{\Hess^{-1}_\la(\hxla)} &= \|\nabla f_\la(\hxla) - \nabla f_\la(\xo_\la)\|_{\Hess^{-1}_\la(\hxla)}\\
&\leq e^{\nm{\hxla - \xo_\la}/2} \dikins(\nm{\hxla - \xo_\la})~\|\hxla - \xo_\la\|_{\Hess_\la(\xo_\la)} \\
 &\leq e^{\nm{\hxla - \xo_\la}/2} \dikins(\nm{\hxla - \xo_\la})~\|\Hess^{1/2}_\la(\xo_\la)\hHess^{-1/2}_\la(\xo_\la)\|~\|\hxla - \xo_\la\|_{\hHess_\la(\xo_\la)}\\
 & \leq \frac{e^{\nm{\hxla - \xo_\la}/2}\dikins(\nm{\hxla - \xo_\la})}{\dikin(\nm{\hxla - \xo_\la})}~\|\Hess^{1/2}_\la(\xo_\la)\hHess^{-1/2}_\la(\xo_\la)\|^2~\|\nabla \hf_\la(\xo_\la)\|_{\Hess_\la(\xo_\la)}\\
 & = e^{3\nm{\hxla - \xo_\la}/2} \bvar_\la .
}

Putting things together: 

\[\hsrce_\la \leq \|\hHess_\la^{-1/2}(\hxla)\Hess_\la^{1/2}(\hxla)\| ~  \left(e^{3 \tn(\xla - \hxla)/2}\left(1 + \frac{1}{\dikin(\tn(\xla - \hxla))}\right) \bvar_\la + \frac{e^{\tn(\xla - \hxla)/2}}{\dikin(\tn(\xla - \xo))}\bias_\la \right). \]

\paragraph{2. We now import the results from \cite{marteau2019}}. The following hypotheses imply those of Thms 4 and 7 in \cite{marteau2019}:

Let $n \in \N$, $\delta \in (0,1/2]$, $0 < \la \leq \Bts$. Whenever
\[n \geq \triangle_1 \frac{\Bts}{\la}\log \frac{8 \square_1^2\Bts}{\la \delta},~~~n \geq \triangle_2~  \frac{\Deff_\lambda \vee (\Qs)^2}{\rla(\xo)^2} ~\log \frac{2}{\delta}, \bias_\la \leq \frac{\rla(\xo)}{2}.\]

In particular, they imply that with probability at least $1 - 2 \delta$: 

\begin{itemize}
    \item $\bvar_\la \leq \frac{1}{2}\bias_\la + 4 \square_1~\sqrt{ \frac{\Deff_\la \vee (\Qs)^2~\log \frac{2}{\delta}}{n}}$;
    \item $\|\Hess^{1/2}_\la(\xo_\la)\hHess^{-1/2}_\la(\xo_\la)\| \leq \sqrt{2}$;
    \item $\nm{\xo - \xo_\la} \leq \log 2$;
    \item $\nm{\hxla - \xo_\la} \leq \log 2$.
\end{itemize}

Hence, plugging these bounds in the previous equation, we get 

\[\hsrce_\la \leq 8 \bias_\la + 80 \sqrt{\frac{\Deff_\la \vee (\Qs)^2~\log \frac{2}{\delta}}{n}}.\]

\paragraph{3.} Note that in what has been done previously, we can bound: 

\[\|\hxla - \xla\|_{\Hess_\la(\xla)} \leq \frac{1}{\dikin(\tn(\xla - \hxla))} \bvar_\la \leq \bias_\la + 8 ~\sqrt{ \frac{\Deff_\la \vee (\Qs)^2~\log \frac{2}{\delta}}{n}}.\]

Moreover,

\[\|\xla - \xo\|_{\Hess_\la(\xo)} \leq \frac{1}{\dikin(\tn(\xla - \xo))} \|\nabla f_\la(\xo)\|_{\Hess^{-1}_\la(\xo)} \leq 2 \bias_\la.\]
Hence: 
\[\|\hxla - \xo\| \leq 3~\la^{-1/2}~\bias_\la + 8~\la^{-1/2} ~\sqrt{ \frac{\Deff_\la \vee (\Qs)^2~\log \frac{2}{\delta}}{n}}.\]
\epr

\paragraph{2. Final bound for the projected ERM approximation}

In this paragraph, denote with $\cp(x,\la)$ the quantity $\frac{\|\Hess^{1/2}(x)(\Id -\Pj)\|^2}{\la}$ and $\hcp(x,\la)$ the quantity $\frac{\|\hHess^{1/2}(x)(\Id -\Pj)\|^2}{\la}$

\blm \label{lm:technical_before_proj} 
Let $n \in \N$, $\delta \in (0,1/2]$, $0 < \la \leq \Bts$. Whenever
\[n \geq \triangle_1 \frac{\Bts}{\la}\log \frac{8 \square_1^2\Bts}{\la \delta},\qquad  \Cone \sqrt{ \frac{\Deff_\lambda \vee (\Qs)^2}{n} ~\log \frac{2}{\delta}} \leq \rla(\xo), \qquad \Cone \bias_\la \leq \rla(\xo),\]
if
\[\cp(\xo,\la) \leq \frac{\sqrt{2}}{480},~ \hnd_{\Pj,\la}(x) \leq \frac{\widehat{\radd}_{\Pj,\la}(x)}{2} ~\wedge~\frac{\rla(\xo)}{126},\]

the following holds, with probability at least $1 - 2 \delta$.

\[\hnd_\la(x) \leq \frac{\hrla(x)}{2},~ \hnd_\la(x) \leq \frac{\rla(\xo)}{2}.\]
Here, $\Cone = 1008$.
\elm

\bpr

Proceed in the following way.
\paragraph{1. } It is easy to see that the conditions of this lemma imply the conditions of \cref{thm:general-result}. Hence, as in the previous proofs, the following hold:

\begin{itemize}
    \item $\|\Hess^{1/2}_\la(\xo_\la)\hHess^{-1/2}_\la(\xo_\la)\| \leq \sqrt{2}$;
    \item $\nm{\xo - \xo_\la} \leq \log 2$;
    \item $\nm{\hxla - \xo_\la} \leq \log 2$.
\end{itemize}

\paragraph{2. } Let us now apply \cref{prp:link_proj} to $\hf$. 
If 
\[\frac{\hsrce_\la}{\hrla(\hxla)} \leq \frac{1}{4},~ \hcp(\hxla,\la) \leq \frac{1}{120},~ \hnd_{\Pj,\la}(x) \leq \frac{\widehat{\radd}_{\Pj,\la}(x)}{2}, \]
Then it holds: 
\eqal{\label{eq:res_proj_lemma}\hnd_\la(x) \leq 3(\hnd_{\Pj,\la}(x) + \hsrce_\la),\qquad \htn(x-\hxla) \leq \frac{3}{10} +2 \log 2.}
where the second bound is obtained in the proof of this proposition.
Now since 

\eqals{
\frac{1}{\hradd_{\la}(\hxla)} &\leq e^{\htn(\hxla - \xla)/2}\frac{1}{\hradd_\la(\xla)} & \cref{eq:hl}\\
& \leq e^{\htn(\hxla - \xla)/2}~\|\Hess^{1/2}_\la(\xo_\la)\hHess^{-1/2}_\la(\xo_\la)\|~\sup_{g \in \hG}{\|g\|_{\Hess^{-1}_\la(\xla)}} & \text{Def}\\
&\leq e^{\tn(\hxla - \xla)/2}~\|\Hess^{1/2}_\la(\xo_\la)\hHess^{-1/2}_\la(\xo_\la)\|~\sup_{g \in \G}{\|g\|_{\Hess^{-1}_\la(\xla)}} & \hG\subset \G\\
& =  e^{\tn(\hxla - \xla)/2}~\|\Hess^{1/2}_\la(\xo_\la)\hHess^{-1/2}_\la(\xo_\la)\|~\frac{1}{\rla(\xla)}& \text{Def}\\
& \leq e^{(\tn(\hxla - \xla) + \tn(\xla - \xo))/2}~\|\Hess^{1/2}_\la(\xo_\la)\hHess^{-1/2}_\la(\xo_\la)\|~\frac{1}{\rla(\xo)} & \cref{eq:hl} \\
& \leq \frac{2\sqrt{2}}{\rla(\xo)}. & \text{previous bounds}
}
In a similar way, we get $\hcp(\hxla,\la) \leq 2 \sqrt{2} \cp(\xo,\la)$. Thus, the conditions above are satisfied if the following conditions are satisfied:
\[\frac{\hsrce_\la}{\rla(\xo)} \leq \frac{\sqrt{2}}{16},~ \cp(\xo,\la) \leq \frac{\sqrt{2}}{480},~ \hnd_{\Pj,\la}(x) \leq \frac{\widehat{\radd}_{\Pj,\la}(x)}{2}.\]
Finally, note that under these conditions,
\eqal{\label{eq:techh}\frac{1}{\hradd_{\la}(x)} \leq \frac{e^{\htn(x - \hxla)/2}}{\hradd_{\la}(x)} \leq \frac{7}{\rla(\xo)}.}
using the previous bound and the bound on $\htn(x -\hxla)$. 
\paragraph{3.} Let us assume 
\[\frac{\hsrce_\la}{\rla(\xo)} \leq \frac{\sqrt{2}}{16},~ \cp(\xo,\la) \leq \frac{\sqrt{2}}{480},~ \hnd_{\Pj,\la}(x) \leq \frac{\widehat{\radd}_{\Pj,\la}(x)}{2}.\]
According to \cref{eq:techh}, and to \cref{eq:res_proj_lemma}, if
\[\hnd_{\Pj,\la}(x) + \hsrce_\la \leq \frac{\rla(\xo)}{42},\]
then it holds
\[\hnd_\la(x) \leq \frac{\hrla(x)}{2},~ \hnd_\la(x) \leq \frac{\rla(\xo)}{2}.\]
We simplify this condition as: 
\[\hnd_{\Pj,\la}(x) \leq \frac{\rla(\xo)}{126},\qquad \hsrce_\la \leq \frac{2 \rla(\xo)}{126}.\]
\paragraph{4.}
Now using the fact that under the conditions of this lemma, those of \cref{lm:emp_source} are satisfied:  
\[\hsrce_\la \leq 8~ \bias_\la + 80~ \sqrt{\frac{\Deff_\la \vee (\Qs)^2~\log \frac{2}{\delta}}{n}}.\]
Thus, $\hsrce_\la \leq \frac{2 \rla(\xo)}{126}$ holds, provided
\[\bias_\la \leq \frac{\rla(\xo)}{\Cone},~ n \geq \Cone^2 \frac{\Deff_\la \vee (\Qs)^2~\log \frac{2}{\delta}}{\rla(\xo)^2} ,\]
where $\Cone = 1008 $.
\epr

\bp[Behavior of an approximation to the projected problem] \label{prp:proj_pb_stat_gen}
Let $n \in \N$, $\delta \in (0,1/2]$, $0 < \la \leq \Bts$. Let $x \in \hh_\Pj$. Whenever
\[n \geq \triangle_1 \frac{\Bts}{\la}\log \frac{8 \square_1^2\Bts}{\la \delta},\qquad  \Cone \sqrt{ \frac{\Deff_\lambda \vee (\Qs)^2}{n} ~\log \frac{2}{\delta}} \leq \rla(\xo), \qquad \Cone \bias_\la \leq \rla(\xo),\]
if
\[\cp(\xo,\la) \leq \frac{\sqrt{2}}{480},~ \hnd_{\Pj,\la}(x) \leq \frac{\widehat{\radd}_{\Pj,\la}(x)}{2} ~\wedge~\frac{\rla(\xo)}{126}. \]

The following holds, with probability at least $1 - 2 \delta$.

\[f(x) - f(\xo) \leq \Cfone~ \bias_\la^2 + \Cftwo~ \frac{\Deff_\la \vee (\Qs)^2}{n}~\log \frac{2}{\delta} + \Cfthree~ \hnd^2_{\Pj,\la}(x), \]

where $\Cfone \leq 6.0\mathrm{e}{4}$, $\Cftwo \leq 6.0\mathrm{e}{6} $ and $\Cfthree \leq 810$, $\Cone$ are defined in \cref{lm:technical_before_proj}, and the other constants are defined in \cref{thm:general-result}. 
\ep 

\br[Constants]
In this result, absolutely huge constants are obtained. They are (of course) totally sub-optimal. Indeed, this analysis has been simplified by dividing the bound into blocks: error of the empirical risk minimization with regularization, error of the projection compared to this empirical risk minimizer. Going back and forth from empirical to statistical, from projected to non projected induces exponential explosion of the constants. There is a way of doing the analysis directly by projecting the statistical problem. However, in order to relate to our previous work \cite{marteau2019} and avoid re-doing all of our work we discarded this. If we were to perform this more direct analysis, we could keep the constants to a reasonable level, of order $10^2$. 

\er

\bpr 
We apply \cref{lm:approx_emp_pb}, using the previous lemma to guarantee the conditions.
\paragraph{1.} Under the conditions of this proposition, applying \cref{lm:technical_before_proj}, the conditions of \cref{lm:approx_emp_pb} are satisfied. Thus, 
\[f(x) - f(\xo) \leq 14(f(\hxla) - f(\xo)) + 30\hnd_\la(x)^2.\]
Moreover, from the previous proof, 
\[\hnd_\la(x) \leq 3(\hnd_{\Pj,\la}(x) + \hsrce_\la),\]
and seeing as \cref{lm:emp_source} is satisfied,
\[\hsrce_\la \leq 8~ \bias_\la + 80~ \sqrt{\frac{\Deff_\la \vee (\Qs)^2~\log \frac{2}{\delta}}{n}}.\]
This therefore yields:
\[\hnd_\la(x)^2 \leq 27 \hnd_{\Pj,\la}(x)^2 + 1726 \bias_\la^2 + 172600 \frac{\Deff_\la \vee (\Qs)^2~\log \frac{2}{\delta}}{n}.\]
\paragraph{2.} Moreover, from \cref{thm:general-result}, it holds:
\[f(\hxla) - f(\xo) \leq 414~\bias_\la^2 + 414~\frac{\Deff_\la \vee (\Qs)^2}{n}~\log \frac{2}{\delta}.\]
\paragraph{3. }
Putting things together:
\[f(x) - f(\xo) \leq \Cfone~ \bias_\la^2 + \Cftwo~ \frac{\Deff_\la \vee (\Qs)^2}{n}~\log \frac{2}{\delta} + \Cfthree~ \hnd^2_{\Pj,\la}(x).\]
We bound the constants in the theorem.

\epr

\blm

Under the conditions of the previous theorem, the following hold: 
\begin{itemize}
    \item $\frac{1}{\widehat{\radd}_{\Pj,\la}(x)} \leq \frac{8}{\rla(\xo)}$;
    \item $\la^{1/2}\|x - \xo\| \leq 7 \hnd_{\Pj,\la}(x) + 59 \bias_\la + 568\sqrt{\frac{\Deff_\la \vee (\Qs)^2~\log \frac{2}{\delta}}{n}}$;
    \item $\la \|x\|_{\hHess_{\Pj,\la}^{-1}(x)}  \leq 7 \hnd_{\Pj,\la}(x) + 72 \bias_\la + 720\sqrt{\frac{\Deff_\la \vee (\Qs)^2~\log \frac{2}{\delta}}{n}}$.
\end{itemize}
In particular,
$\frac{\la \|x\|_{\hHess_{\Pj,\la}^{-1}(x)}}{\widehat{\radd}_{\Pj,\la}(x)} \leq 11$.
\elm 

\bpr Let us prove the three statements. 

\paragraph{1.} Write $\frac{1}{\widehat{\radd}_{\Pj,\la}(x)} = \sup_{g \in \hG}{\|\Pj g\|_{\hHess^{-1}_{\Pj,\la}(x)}}$. Now
\[\sup_{g \in \hG}{\|\Pj g\|_{\hHess^{-1}_{\Pj,\la}(x)}} \leq \sup_{g \in \hG}{\| g\|_{\hHess^{-1}_{\la}(x)}} \leq e^{\htn(x - \hxla)/2}\sup_{g \in \hG}{\| g\|_{\hHess^{-1}_{\la}(\hxla)}}.\]

Now bound

\[\sup_{g \in \hG}{\| g\|_{\hHess^{-1}_{\la}(\hxla)}} \leq e^{\htn(\xla - \hxla)/2}\sup_{g \in \hG}{\| g\|_{\hHess^{-1}_{\la}(\xla)}} \leq e^{\htn(\xla - \hxla)/2}~\|\Hess^{1/2}_\la(\xla)\hHess^{-1/2}_\la(\xla)\|~\sup_{g \in \hG}{\| g\|_{\Hess^{-1}_{\la}(\xla)}}. \]

Finally bound 
\[\sup_{g \in \hG}{\| g\|_{\Hess^{-1}_{\la}(\xla)}} \leq e^{\tn(\xo - \xla)/2} \frac{1}{\rla(\xo)}.\]

Now using the fact that under the previous assumptions $\tn(\xo - \xla ), \tn(\xla - \hxla) \leq \log 2$, $\htn(x - \hxla) \leq \frac{3}{10} + 2 \log 2$ and $\|\Hess^{1/2}_\la(\xla)\hHess^{-1/2}_\la(\xla)\| \leq \sqrt{2}$, we get the first equation.

\paragraph{2. } In order to bound $\la^{1/2}\|x - \xo\|$, decompose
\[\la^{1/2}\|x - \xo\| \leq \la^{1/2}\|x - \hxla\| + \la^{1/2}\|\hxla - \xo\|.\]
Now use \cref{prp:link_proj} to bound
$  \la^{1/2}\|x - \hxla\| \leq 7 (\hnd_{\Pj,\la}(x) + \hsrce_\la) $. Using \cref{lm:emp_source}, under the conditions above, 

\[\hsrce_\la \leq 8~ \bias_\la + 80~ \sqrt{\frac{\Deff_\la \vee (\Qs)^2~\log \frac{2}{\delta}}{n}}.\]

Hence

\[ \la^{1/2}\|x - \hxla\| \leq 7 \hnd_{\Pj,\la}(x) + 56 \bias_\la + 560\sqrt{\frac{\Deff_\la \vee (\Qs)^2~\log \frac{2}{\delta}}{n}}. \]
Moreover, using again \cref{lm:emp_source}

\[\la^{1/2}\|\hxla - \xo\| \leq 3~\bias_\la + 8~ ~\sqrt{ \frac{\Deff_\la \vee (\Qs)^2~\log \frac{2}{\delta}}{n}}.\]
Combining these two inequalities, we get:

\[ \la^{1/2}\|x - \xo\| \leq 7 \hnd_{\Pj,\la}(x) + 59 \bias_\la + 568\sqrt{\frac{\Deff_\la \vee (\Qs)^2~\log \frac{2}{\delta}}{n}}.\]

\paragraph{3. } In order to bound $\la \|x\|_{\hHess_{\Pj,\la}^{-1}(x)}$, use \cref{prp:link_proj} to get  $\la \|x\|_{\hHess^{-1}_{\Pj,\la}(x)} \leq 7 \hnd_{\Pj,\la}(x) + 9 \hsrce_\la$.

Now using \cref{lm:emp_source}, the following bound holds:

\[ \la \|x\|_{\hHess_{\Pj,\la}^{-1}(x)}  \leq 7 \hnd_{\Pj,\la}(x) + 72 \bias_\la + 720\sqrt{\frac{\Deff_\la \vee (\Qs)^2~\log \frac{2}{\delta}}{n}}.\]

\epr

\bp[Simplification] 

\label{prp:proj_pb_stat_gen_part}
Let $n \in \N$, $\delta \in (0,1/2]$, $0 < \la \leq \Bts$. Let $x \in \hh_\Pj$. Whenever
\[n \geq \triangle_1 \frac{\Bts}{\la}\log \frac{8 \square_1^2\Bts}{\la \delta},\qquad  \Cone \sqrt{ \frac{\Deff_\lambda \vee (\Qs)^2}{n} ~\log \frac{2}{\delta}} \leq \frac{\sqrt{\la}}{\Rad}, \qquad \Cone \bias_\la \leq \frac{\sqrt{\la}}{\Rad},\]
if
\[\cp(\xo,\la) \leq \frac{\sqrt{2}}{480},~ \hnd_{\Pj,\la}(x) \leq \frac{\sqrt{\la}}{126\Rad},\]

then the following holds, with probability at least $1 - 2 \delta$.

\[f(x) - f(\xo) \leq \Cfone~ \bias_\la^2 + \Cftwo~ \frac{\Deff_\la \vee (\Qs)^2}{n}~\log \frac{2}{\delta} + \Cfthree~ \hnd^2_{\Pj,\la}(x), \]

where $\Cfone \leq 6.0\mathrm{e}{4}$, $\Cftwo \leq 6.0\mathrm{e}{6} $ and $\Cfthree \leq 810$, $\Cone$ are defined in \cref{lm:technical_before_proj}, and the other constants are defined in \cref{thm:general-result}. 

Moreover, in that case, $\Rad \|x - \xo\| \leq 10$. 

\ep

\subsection{Optimal choice of $\la$, specific source conditions \label{app:quantibounds}}

In this part, we continue to assume \cref{asm:gen_sc,asm:asmtech,asm:iid,asm:minimizer}. We present a classification of distributions $\rho$ and show that we can achieve better rates than the classical slow rates, as presented in Appendix F of \cite{marteau2019}. 

\subsubsection{Classification of distributions and statistical bounds for the ERM}

We use the following classification for distributions.

\bd[class of distributions]
Let $\alpha \in [1,+\infty]$ and $r \in [0,1/2]$. \\\
We denote with $ \clas_{\alpha,r}$ the set of probability distributions $\rho$ such that there exists $\Lc,\Qc \geq 0$, 
\begin{itemize}
    \item $\bias_\la \leq \Lc ~\lambda^{\frac{1 + 2r}{2}}$;
    \item $\Deff_\la \leq \Qc^2~ \la^{-1/\alpha}$;
\end{itemize}
where this holds for any $0<\la \leq 1$. 
For simplicity, if $\alpha = +\infty$, we assume that $\Qc \geq \Qs$.
\ed

 Note that given our assumptions, we always have 

\eqal{\label{eq:basic_class}\rho \in \clas_{1,0},~~~  \Lc = \|\xo\|,~\Qc = \Bos.}

We also define 
\eqal{ \label{eq:lambda1}\la_1 = \left(\frac{\Qc}{\Qs}\right)^{2\alpha} \wedge 1 , }
such that 
$$ \forall \la \leq \la_1,~ \Deff_\la \vee (\Qs)^2 \leq \frac{\Qc^2}{\la^{1/\alpha}} .$$

\paragraph{Interpretation of the classes}

\begin{itemize}
    \item The bias term $\bias_\la$ characterizes the regularity of the objective $\xo$. In a sense, if $r$ is big, then this means $\xo$ is very regular and will be easier to estimate. The following results reformulates this intuition.

\br[source condition]\label[remark]{rmk:source_cond}
Assume there exists $0 \leq r \leq 1/2$ and $v \in \hh$ such that 
\[\Pj_{\Hess(\xo)}\xo = \Hess(\xo)^r v.\]
Then it holds:  
\[\forall \la >0,~\bias_\la \leq \Lc~\la^{\frac{1 + 2r}{2}},~~~~\Lc = \|\Hess(\xo)^{-r}\xo\|.\]
\er

 \item The effective dimension $\Deff_\la$ characterizes the size of the space $\hh$ with respect to the problem. The higher $\alpha$, the smaller the space. If $\hh$ is finite dimensional for instance, $\alpha = +\infty$.
 
 \end{itemize}
 
 In this section, for any given pair $(\alpha,r)$ characterizing the regularity and size of the problem, we associate 
\[\beta = \frac{1}{1 + 2r + 1/\alpha},~~~~\gamma = \frac{\alpha(1+2r)}{\alpha(1+2r) + 1}.\]

 In \cite{marteau2019} (see corollary 3), explicit bounds are given for the performance of the regularized expected risk minimizer $\hxla$ depending on which class $\rho$ belongs to, i.e., as a function of $\alpha,r$.
 
 \bcor\label[corollary]{cor:bounds_source_capacity}
Let $\delta \in (0,1/2]$. Under \cref{asm:gen_sc,asm:asmtech,asm:iid,asm:minimizer}, if $\rho \in \clas_{\alpha,r}$ with $r>0$ , when $n \geq N$ and $\la = (C_0/n)^{\beta}$, then with probability at least $1-2\delta$,
\[f(\hxla) - f(\xo) \leq C_1 n^{-\gamma}\log \frac{2}{\delta},\]
with $C_0 = 256 (\Qc/\Lc)^2,~C_1 = 8 (256)^{\gamma}~(\Qc^{\gamma} ~ 
    \Lc^{1-\gamma})^2$ and $N$ defined in \cite{marteau2019}, and satisfying
$N = O(\poly(\Bos,\Bts,\Lc,\Qc,\Rad,\log(1/\delta)))$. 
\ecor

\subsubsection{Quantitative bounds for the projected problem}
In this part,  the aim is to show that if we approximately solve the projected problem up to a certain precision, then this approximation has the same statistical rates as the regularized ERM with the good choice of $\la$. For the sake of simplicity, we will assume that $r >0$.

In what follows, we define 

\eqal{\label{eq:Nvalue}N = \frac{ \Qc^2}{\Lc^2}\left(\Bts \wedge \la_0 \wedge \la_1\right)^{-1/\beta} ~~~\vee~~~\left(2.1e4 \frac{1}{1-\beta}  A \log  \left(1.4e6   \frac{1}{1-\beta}A^2 \frac{1}{\delta}\right)\right)^{1/(1-\beta)},}
where $A = \frac{\Bts \Lc^{2\beta} }{ \Qc^{2\beta} }$, $\la_0 = (\Cone \Lc \Rad \log \frac{2}{\delta})^{-1/r} \wedge 1$ and $\la_1 = \frac{\Qc^{2\alpha}}{(\Qs)^{2\alpha}}$.

\bt[Quantitative result with source $r > 0$]\label[theorem]{thm:quantisource}
Let $\rho \in \clas_{\alpha,r}$  and assume $r>0$. Let $\delta \in (0,\frac{1}{2}]$.\\\
Let $\Pj$ be an orthogonal projection, $x \in \hh$.
If 
    \[n \geq N,~~~\lambda = \left( \left(\frac{\Qc}{\Lc}\right)^{2}~\frac{1}{n}\right)^{\beta},~~~~\cp(\xo,\la) \leq \frac{\sqrt{2}}{480},~ \hnd_{\Pj,\la}(x) \leq \Qc^{\gamma} ~ 
    \Lc^{1-\gamma} n^{-\gamma/2}\]
    then with probability at least $1-2 \delta$,
    \[f(x) - f(\xo) \leq \Kfinal \left(\Qc^{\gamma} ~ 
    \Lc^{1-\gamma}\right)^2 ~~\frac{1}{n^{\gamma}} \log \frac{2}{\delta},\]
    where $N$ is defined in \cref{eq:Nvalue} and $\Kfinal \leq 7.0e6$. Moreover, $\Rad \|x - \xo\| \leq 10$.

\et

\bpr 

Using the definition of $\la_1$, as soon as $\la \leq \la_1$ ,it holds:  $\Deff_\la \vee (\Qs)^2 \leq \Qc^2 \la^{-1/\alpha}$.\\\

Let us formulate  \cref{prp:proj_pb_stat_gen_part} using the fact that $\rho \in \clas_{\alpha,r}$.\\\

Let $n \in \N$, $\delta \in (0,1/2]$, $0 < \la \leq \Bts$, $x \in \hh_\Pj$. Whenever
\[n \geq \triangle_1 \frac{\Bts}{\la}\log \frac{8 \square_1^2\Bts}{\la \delta},~~~\Cone \sqrt{ \frac{\Qc^2}{\la^{1/\alpha}n} ~\log \frac{2}{\delta}}\leq \frac{\la^{1/2}}{\Rad}, \Cone ~\Lc \la^{1/2 + r} \leq \frac{\la^{1/2}}{\Rad},\]
if
\[\cp(\xo,\la) \leq \frac{\sqrt{2}}{480},~ \hnd_{\Pj,\la}(x) \leq \Lc \la^{1/2  + r}, \]

The following holds, with probability at least $1 - 2 \delta$.

\[f(x) - f(\xo) \leq (\Cfone + \Cfthree) \Lc^2 \la^{1 + 2r} + \Cftwo~ \frac{\Qc^2}{\la^{1/\alpha} n}~\log \frac{2}{\delta}, \qquad \Rad \|x - \xo\|\leq 10,\]

where all constants are defined in \cref{prp:proj_pb_stat_gen_part}.

\paragraph{Assume that $r > 0$}.
Define 
\[\la_0 = (\Cone \Lc \Rad \log \frac{2}{\delta})^{-1/r} \wedge 1.\]
Then for any $\la \leq \la_0$:
\[ \Lc \la^{1/2 + r} \leq \frac{1}{\Cone}\frac{\sqrt{\la}}{\Rad}.\]

\noindent{\normalfont \bfseries 1)} First, we find a simple condition to guarantee 
$$\rla(\xo)^2\lambda^{ 1/\alpha} \geq \Ctwo~  \Qc^2
\frac{1}{n} ~\log \frac{2}{\delta}.$$

We see that if $\lambda \leq  \lambda_0$, then $\rla \geq \Cone \Lc \la^{1/2 + r} \log \frac{2}{\delta}$. Hence, this condition is satisfied if 
$$ \la \leq \la_0,~~~~\Cone^2 \Lc^2 \la^{1+ 2r + 1/\alpha} \geq \Ctwo~  \Qc^2
\frac{1}{n} .$$ 
Using the fact that $\Ctwo = \Cone^2$, we reformulate:
$$ \la \leq \la_0,~~~~\Lc^2 \la^{1+ 2r + 1/\alpha} \geq   \Qc^2
\frac{1}{n} .$$ 

\noindent{\normalfont \bfseries 2) } Now fix  
$$\lambda^{1 + 2r + 1/\alpha} =\frac{\Qc^2}{ \Lc^2} ~\frac{1}{n} \Longleftrightarrow \lambda = \left(\frac{\Qc^2}{ \Lc^2} ~\frac{1}{n}\right)^{\beta} .$$
where $\beta = 1/(1 + 2r + 1/\la) \in [1/2,1)$.

Using our restatement of \cref{prp:proj_pb_stat_gen}, with probability at least $1-2 \delta$,
$$\L(x) - \L(\xo) \leq 
  \left(\Cfone + \Cfthree + \Cftwo \log \frac{2}{\delta}\right)~\Lc^2 \lambda^{1 + 2r} \leq    \Kfinal~\log \frac{2}{\delta}~\Lc^2 \lambda^{1 + 2r},$$
  where $\Kfinal = \Cfone + \Cfthree + \Cftwo \leq 7.0e6$ (see \cref{prp:proj_pb_stat_gen}). \\\
  This result holds provided
     \eqal{\label{eq:cond_refor}0 < \lambda \leq \Bts \wedge \lambda_0 \wedge \la_1,~n \geq \triangle_1 \frac{\Bts}{\la}\log \frac{8 \square_1^2\Bts}{\la \delta}.}
     Indeed, it is shown in the previous point that the other conditions are satisfied.

\noindent{\normalfont \bfseries 3) }  Let us now work to guarantee the conditions in \cref{eq:cond_refor}.  \\\
First, to guarantee $n \geq \triangle_1 \frac{\Bts}{\la}\log \frac{8 \square_1^2\Bts}{\la \delta}$, bound
\[\frac{\Bts}{\la} = \frac{\Bts \Lc^{2\beta} n^{\beta}}{ \Qc^{2\beta} ~\log^{\beta} \frac{2}{\delta}} \leq 2~\frac{\Bts \Lc^{2\beta} }{ \Qc^{2\beta} } n^{\beta}.\]

Then apply lemma 15 from \cite{marteau2019} with $a_1 = 2\triangle_1$, $a_2 = 16 \square_1^2$, $A = \frac{\Bts \Lc^{2\beta} }{ \Qc^{2\beta} }$. Since $\beta \geq 1/2$, using the bounds in \cref{thm:general-result}, we find $a_1 \leq 10400$ and $a_2 \leq 64,$ hence the following sufficient condition:

\[ n \geq  \left(2.1e4 \frac{1}{1-\beta}  A \log  \left(1.4e6   \frac{1}{1-\beta}A^2 \frac{1}{\delta}\right)\right)^{1/(1-\beta)}.\]
  
 Then, to guarantee the condition 
\[ \lambda \leq \Bts \wedge \lambda_0 \wedge \la_1, \]
we simply need 
\[n \geq \frac{\Qc^2}{\Lc^2}\left(\Bts \wedge \la_0 \wedge \la_1\right)^{-1/\beta}.\]
Hence, defining 
$$N = \frac{ \Qc^2}{\Lc^2}\left(\Bts \wedge \la_0 \wedge \la_1\right)^{-1/\beta} ~~~\vee~~~\left(2.1e4 \frac{1}{1-\beta}  A \log  \left(1.4e6   \frac{1}{1-\beta}A^2 \frac{1}{\delta}\right)\right)^{1/(1-\beta)},$$ we see that as soon as $n \geq N$, \cref{eq:cond_refor} holds.

\epr

\newpage 

\section{Multiplicative approximations for Hermitian operators\label{sec:lemmas-operators}}

In this section, we put together useful tools for approximating linear operators and solving linear systems with regularization. 

In this section, $\Ab$ and $\mathbf{B}$ will always denote positive semi-definite Hermitian operators on a Hilbert space $\hh$, and $\Pj$ will denote an orthogonal projection operator. Moreover, given a positive semi-definite operator $\Ab$, and $\la >0$, $\Ab_\la$ will stand for the regularized operator $\Ab + \la \Id$.

\blm[Equivalence of Hermitian operators]\label{lm:equiv_operators}
Let $\mathbf{A}$ and $\mathbf{B}$ be two semi-definite Hermitian operators. Let $\la > 0$. Assume you have access to 
\[t:= \|\mathbf{A}_\la^{-1/2}(\mathbf{B} - \mathbf{A})\mathbf{A}^{-1/2}_\la\|.\]
It holds: 
\[\|\mathbf{A}_\la^{-1/2} \mathbf{B}_\la^{1/2}\|^2 \leq 1 + t \Leftrightarrow \mathbf{B}_\la \preceq (1+t)\mathbf{A}_\la.\]
Moreover, if $t<1$, 
\[\|\mathbf{B}_\la^{-1/2} \mathbf{A}_\la^{1/2}\|^2 \leq \frac{1}{1-t} \Leftrightarrow (1-t)\mathbf{A}_\la \preceq \mathbf{B}_\la .\]
\elm

\bpr 
For the first point, simply note that: 
\[\|\Ab_\la^{-1/2}\mathbf{B}_\la^{1/2}\|^2 = \|\Ab_\la^{-1/2} \mathbf{B}_\la \Ab_\la^{-1/2}\| = \|\Id +\Ab_\la^{-1/2} \left(\mathbf{B} - \Ab\right) \Ab_\la^{-1/2} \| \leq 1 + t.\]
For the second point,
\[\|\mathbf{B}_\la^{-1/2}\Ab_\la^{1/2}\|^2 = \|\left(\Ab_\la^{-1/2} \mathbf{B}_\la \Ab_\la^{-1/2}\right)^{-1}\| = \|\left(\Id +\Ab_\la^{-1/2} \left(\mathbf{B} - \Ab\right) \Ab_\la^{-1/2} \right)^{-1}\| .\]
Moreover, we know that if $\|\mathbf{H}\| < 1$ with $\mathbf{H}$ a Hermitian operator, then $\|(\Id + \mathbf{H})^{-1}\| \leq \frac{1}{1-\|\mathbf{H}\|}$. The result follows. 
\epr

We will now state a technical lemma which describes how combining approximation behaves. 

\blm[Combination of approximations]\label{lm:comb_approx}
Let $N \geq 1$. Let $(\Ab_i)_{1 \leq i \leq N+1}$ be a sequence of positive semi-definite Hermitian operators. Define
\[t_i:= \|\Ab_{i,\la}^{-1/2}(\Ab_{i+1} - \Ab_i)\Ab_{i,\la}^{-1/2}\|.\]
For any $1 \leq i,j \leq N+1$, define 
$$t_{i:j}:= \|\Ab_{i,\la}^{-1/2}(\Ab_{j} - \Ab_i)\Ab_{i,\la}^{-1/2}\|.$$
In particular, $t_i = t_{i:i+1}$. Then the following holds:

\[\forall 1 \leq i \leq j \leq N,~1 + t_{i:j} \leq \prod_{k=i}^{j-1}(1 + t_k)\]

Moreover, if $t_i < 1$, then it holds: 
\[\|\Ab_{i+1,\la}^{-1/2}(\Ab_{i} - \Ab_{i+1})\Ab_{i+1,\la}^{-1/2}\| \leq \frac{t_i}{1-t_i}\]
Hence, in that case
\[\forall 1 \leq j \leq i \leq N,~ 1 + t_{j:i} \leq \prod_{k=i}^{j-1}{\frac{1}{1-t_k}} \]

\elm 

\bpr
Let us prove everything for a sequence of three operators; the rest follows by induction. Let $\Ab_1,\Ab_2,\Ab_3$ be three positive semi-definite operators. 
\paragraph{1.} Bound
\eqals{t_{1:3} &= \|\Ab_{1,\la}^{-1/2}\left(\Ab_1 - \Ab_3\right)\Ab_{1,\la}^{-1/2}\| \\
&\leq \|\Ab_{1,\la}^{-1/2}\left(\Ab_1 - \Ab_2\right)\Ab_{1,\la}^{-1/2}\|  + \|\Ab_{1,\la}^{-1/2}\left(\Ab_2 - \Ab_3\right)\Ab_{1,\la}^{-1/2}\|\\
&\leq t_{1:2} + \|\Ab_{1,\la}^{-1/2}\Ab_{2,\la}^{1/2}\|^2 t_{2:3}\\
& \leq t_{1:2} + (1 + t_{1:2})t_{2:3}
.}
The last line comes from \cref{lm:equiv_operators}. Thus 
\[1 + t_{1:3} \leq 1 + t_{1:2} + t_{2:3} + t_{1:2} t_{2:3} = (1 + t_{1:2})(1 + t_{2:3}).\]
\paragraph{2.} Let us now bound $t_{2:1}$ knowing $t_{1:2}$. This will imply the rest of the lemma. 

Indeed, simply note that 
\[t_{2:1} = \|\Ab^{-1/2}_{2,\la}(\Ab_2 - \Ab_1)\Ab^{-1/2}_{2,\la} \| \leq \|\Ab^{-1/2}_{2,\la}\Ab^{1/2}_{1,\la}\|^2 ~t_{1:2} .\]
Using \cref{lm:equiv_operators}, if $t_{1:2} < 1$, $\|\Ab^{-1/2}_{2,\la}\Ab^{1/2}_{1,\la}\|^2 \leq \frac{1}{1-t_{1:2}}$, hence
\[t_{2:1} \leq \frac{t_{1:2}}{1-t_{1:2}}.\]
\epr

\blm[Projection of Hermitian operators]\label{lm:proj_hermit}
For any Hermitian operator $\Ab$ and orthogonal projection $\Pj$, the following holds:
\[\|\Ab_\la^{-1/2}(\Ab - \Pj \Ab \Pj )\Ab_\la^{-1/2}\| \leq \left(1 + \frac{\|\Ab^{1/2}(\Id - \Pj)\|}{\sqrt{\la}}\right)^2 - 1.\]
In particular, 
\[\|\Ab^{-1/2}_\la \left(\Pj \Ab \Pj + \la \Id\right)^{1/2}\| \leq 1 + \frac{\|\Ab^{1/2}(\Id - \Pj)\|}{\sqrt{\la}}.\]
Moreover, if 
\[\frac{\|\Ab^{1/2}(\Id - \Pj)\|}{\sqrt{\la}} < \sqrt{2}-1,\]
then it holds 
\[\|\Ab^{1/2}_\la \left(\Pj \Ab \Pj + \la \Id\right)^{-1/2}\|^2 \leq \frac{1}{2 - \left(1 + \frac{\|\Ab^{1/2}(\Id - \Pj)\|}{\sqrt{\la}}\right)^2}.\]
We also always have:
\[\| \left(\Pj \Ab \Pj + \la \Id\right)^{-1/2} \Pj \Ab^{1/2}_\la\|^2 \leq 1 .\]
\elm

\bpr 
For any Hermitian operator $\Ab$, the following computation holds:
\eqals{\|\Ab_\la^{-1/2}(\Ab - \Pj \Ab \Pj )\Ab_\la^{-1/2}\| &=
\|\Ab_\la^{-1/2}(\Ab - (\Id - (\Id - \Pj)) \Ab (\Id - (\Id - \Pj ))\Ab_\la^{-1/2}\|\\
& \leq 2\|\Ab_\la^{-1/2}(\Id-\Pj)\Ab \Ab_\la^{-1/2}\| + \|\Ab_\la^{-1/2}(\Id-\Pj)\Ab (\Id - \Pj)\Ab_\la^{-1/2}\|\\
& \leq \frac{2\|\Ab^{1/2}(\Id - \Pj)\|}{\sqrt{\la}} + \frac{\|\Ab^{1/2}(\Id-\Pj)\|^2}{\la}\\
& = \left(1 + \frac{\|\Ab^{1/2}(\Id - \Pj)\|}{\sqrt{\la}}\right)^2 - 1
.}
\epr

\blm[Relationship between approximations]\label{lm:relation_ship}
Let $\Ab$ and $\Bb$ be two positive semi-definite hermitian operators. Let $\la > 0$, $b \in \hh$ and $\rho > 0$.  If
\[\|\Ab_\la^{-1/2}(\Bb - \Ab)\Ab_\la^{-1/2}\| \leq \frac{1}{2} \wedge \frac{\err}{4},\qquad \nsa \in \lso(\Bb_\la,b,\err/4) ,\]
then it holds:
$$\nsa \in \lso(\Ab_\la,b,\err).$$

\elm

\bpr  Assume $\nsa \in \lso(\Bb_\la,b,\terr/4)$ for a certain $\terr$. 
Decompose 
\eqals{
\|\Ab_\la^{-1}b - \nsa\|_{\Ab_\la} &\leq \|\Ab_\la^{-1}b - \Bb_\la^{-1}b\|_{\Ab_\la}  + \|\Bb_\la^{-1}b - \nsa\|_{\Ab_\la} \\
&\leq \|\Ab_\la^{1/2}(\Ab_\la^{-1} - \Bb_\la^{-1})\Ab_\la^{1/2}\|~ \|b\|_{\Ab_\la^{-1}} + \|\Ab_\la^{1/2}\Bb_\la^{-1/2}\| ~\|\Bb_\la^{-1}b - \nsa\|_{\Bb_\la}.}

Now using the fact that $\Ab_\la^{-1} - \Bb_\la^{-1} = \Bb_\la^{-1}(\Bb - \Ab)\Ab_\la^{-1}$, 

\begin{align*}\|\Ab_\la^{1/2}(\Ab_\la^{-1} - \Bb_\la^{-1})\Ab_\la^{1/2}\| &\leq \|\Ab_\la^{-1/2}(\Bb-\Ab)\Ab_\la^{-1/2}\|~\|\Ab_\la^{1/2}\Bb_\la^{-1} \Ab_\la^{1/2}\| \\
&=  \|\Ab_\la^{-1/2}(\Bb-\Ab)\Ab_\la^{-1/2}\|~\|\Ab_\la^{1/2}\Bb_\la^{-1/2} \|^2.\end{align*}
Moreover, 
\[\|\Bb_\la^{-1}b - \nsa\|_{\Bb_\la} \leq \terr \|b\|_{\Bb_\la^{-1}} \leq \|\Ab^{1/2}\Bb^{-1/2}\|~\|b\|_{\Ab_\la^{-1}}.\]
Putting things together, and noting that from \cref{lm:equiv_operators}, $\|\Ab^{1/2}\Bb^{-1/2}\|^2 \leq \frac{1}{1-\|\Ab_\la^{-1/2}(\Bb - \Ab)\Ab_\la^{-1/2}\|}$ as soon as $\|\Ab_\la^{-1/2}(\Bb - \Ab)\Ab_\la^{-1/2}\| < 1$, it holds: 
$$\nsa \in \lso(\Ab_\la,b,\err),~\err = \frac{\terr + \|\Ab_\la^{-1/2}(\Bb - \Ab)\Ab_\la^{-1/2}\|}{1 - \|\Ab_\la^{-1/2}(\Bb - \Ab)\Ab_\la^{-1/2}\|} .$$

Choosing the right values for $\terr$ and $\|\Ab_\la^{-1/2}(\Bb - \Ab)\Ab_\la^{-1/2}\|$ yields the result.
\epr 
\subsection{Results for Nystrom sub-sampling\label{app:sub_nystrom}}

Recall the notations from \cref{app:kernels}. 

We write without proof the following lemmas, which are just restatements of lemmas 9 and 10 of \cite{Rudi17}. 

\blm[Uniform sampling]
Let $\delta >0$. If $\left\{\tilde{z}_1,...,\tilde{z}_m\right\}$ are sampled uniformly, then if $0 < \la < \|\Ab\|$, $m \leq n$ and 
\[m \geq \left(10 + 160 \Nyi^\Ab(\la)\right)\log \frac{8 \|v\|_{L^{\infty}(Z)}^2}{\la \delta}.\]
Then it holds, with probability at least $1-\delta$:
\[\|\Ab_\la^{-1/2}(\hAb - \Ab)\Ab_\la^{-1/2}\| \leq \frac{1}{2},\qquad\|\hAb_{m,\la}^{-1/2}(\hAb - \hAb_m)\hAb_{m,\la}^{-1/2}\| \leq \frac{1}{2} .\]
\elm

\blm[Nystrom sampling]
Let $\delta >0$. If $\left\{\tilde{z}_1,...,\tilde{z}_m\right\}$ are sampled using $q$-approximate leverage scores for $t = \la$, then if $t_0 \vee \frac{19 \|v\|_{L^{\infty}(Z)}^2}{n}\log \frac{n}{2 \delta} < \la < \|\Ab\|$, and $n \geq 405 \|v\|_{L^{\infty}(Z)}^2 \vee 67 \|v\|_{L^{\infty}(Z)}^2 \log \frac{12 \|v\|_{L^{\infty}(Z)}^2}{\delta}$, if
\[m \geq \left(6 + 486 q^2 \Ny^\Ab(\la)\right)\log \frac{8 \|v\|_{L^{\infty}(Z)}^2}{\la \delta}.\]
Then it holds, with probability at least $1-\delta$:
\[\|\Ab_\la^{-1/2}(\hAb - \Ab)\Ab_\la^{-1/2}\| \leq \frac{1}{2},\qquad\|\hAb_{m,\la}^{-1/2}(\hAb - \hAb_m)\hAb_{m,\la}^{-1/2}\| \leq \frac{1}{2} .\]
\elm

\blm Let $\la >0$. 
Assume: 
\[\|\Ab_\la^{-1/2}(\hAb - \Ab)\Ab_\la^{-1/2}\| \leq \frac{1}{2},\qquad\|\hAb_{m,\la}^{-1/2}(\hAb - \hAb_m)\hAb_{m,\la}^{-1/2}\| \leq \frac{1}{2} .\]
Denote with $P_m$ the projection on $\lspan(v_{\tilde{z}_j})_{1\leq j \leq m}$. Then the following holds:
\[\|\Ab^{1/2}_\la(\Id - \Pj_m)\|^2 \leq 3 \la,\]
and for any partial isometry $V$,
\[\frac{1}{2}\left(V^* \hAb_m V + \la \Id \right)\preceq  V^* \hAb V + \la \Id \preceq \frac{3}{2}\left( V^* \hAb_m V + \la \Id\right).\]

\elm

\bpr 
For the first point, use the well known fact that 
\[\Id - \Pj_m \leq \la \hAb_{m,\la}^{-1},\]
since the range of $\Pj_m$ contains that of $\hAb_m$. 
Thus, 
\[\|\Ab_\la^{1/2}(\Id - \Pj_m)\|^2 \leq \la \|\Ab_\la^{1/2} \hAb_{m,\la}^{-1/2}\|^2.\]
Now using \cref{lm:comb_approx},
\[\|\Ab_\la^{-1/2}(\hAb - \Ab)\Ab_\la^{-1/2}\| \leq \frac{1}{2} \implies \|\hAb_\la^{-1/2}(\hAb - \Ab)\hAb_\la^{-1/2}\| \leq 1.\]
Hence, again using \cref{lm:comb_approx}, 
\[\|\hAb_{m,\la}^{-1/2}(\hAb_m - \Ab)\hAb_{m,\la}^{-1/2}\| \leq 2,\]
and therefore, using \cref{lm:equiv_operators},
\[\|\Ab_\la^{1/2} \hAb_{m,\la}^{-1/2}\|^2 \leq 3.\]
For the second point, this is only a consequence of \cref{lm:equiv_operators}. 
\epr 
Now state two results which show that 

\blm[Uniform sampling yielding $\rho$-approximation]\label{lm:sub_unif}
Let $0 < \rho \leq 1$ and $\delta >0$. Let $b \in \hh$. If $\left\{\tilde{z}_1,...,\tilde{z}_m\right\}$ are sampled uniformly, $0 < \la < \|\Ab\|$, $m \leq n$ and 
\[m \geq \left(2 + \frac{48}{\err} + \frac{5000}{\rho^2} \Nyi^\Ab(\la)\right)\log \frac{8 \|v\|_{L^{\infty}(Z)}^2}{\la \delta}.\]
Then it holds, with probability at least $1-\delta$:
\[x \in \lso(\hAb_{m,\la},b,\err/4) \implies x \in \lso(\Ab_\la,b,\err)  .\]
In particular, with probability $1-\delta$, 
\[\hAb_{m,\la}^{-1}b \in \lso(\Ab_\la,b,\err) . \]
\elm

\bpr 
Apply Lemma 9 from \cite{Rudi17} with $\eta = \frac{\rho}{12} < \frac{1}{2}$. We find that under the conditions above, with probability at least $1-\delta$, 
\[\|\Ab_\la^{-1/2}(\hAb - \Ab)\Ab_\la^{-1/2}\| \leq \eta ,\qquad\|\hAb_{m,\la}^{-1/2}(\hAb - \hAb_m)\hAb_{m,\la}^{-1/2}\| \leq \eta .\]
Now use \cref{lm:comb_approx} to see that
\[\|\Ab_\la^{-1/2}(\hAb_m - \Ab)\Ab_\la^{-1/2}\| \leq (1 + \eta^2) - 1 \leq 3 \eta \leq \err/4.\]
Thus, we can apply \cref{lm:relation_ship} to get the desired result. 
\epr

\blm[Leverage scores Nystrom sampling yielding $\rho$-approximation]\label{lm:sub_ny}
Let $\delta >0$. If $\left\{\tilde{z}_1,...,\tilde{z}_m\right\}$ are sampled using $q$-approximate leverage scores for $t = \la$, then if $t_0 \vee \frac{19 \|v\|_{L^{\infty}(Z)}^2}{n}\log \frac{n}{2 \delta} < \la < \|\Ab\|$, and $n \geq 405 \|v\|_{L^{\infty}(Z)}^2 \vee 67 \|v\|_{L^{\infty}(Z)}^2 \log \frac{12 \|v\|_{L^{\infty}(Z)}^2}{\delta}$, if
\[m \geq \left(2 + \frac{24}{\err} + \frac{13000 q^2}{\rho^2} \Ny^\Ab(\la)\right)\log \frac{8 \|v\|_{L^{\infty}(Z)}^2}{\la \delta}.\]
Then it holds, with probability at least $1-\delta$:
\[x \in \lso(\hAb_{m,\la},b,\err/4) \implies x \in \lso(\Ab_\la,b,\err)  .\]
In particular, with probability $1-\delta$, 
\[\hAb_{m,\la}^{-1}b \in \lso(\Ab_\la,b,\err) . \]
\elm

\bpr 
The proof is exactly the same as that of the previous lemma, using Lemma 10 instead of Lemma 9 in \cite{Rudi17}. 
\epr 

\end{document}